\documentclass[]{informs3} 

\OneAndAHalfSpacedXI 

\usepackage[normalem]{ulem}  
\useunder{\uline}{\ul}{}
\usepackage{hyperref}
\hypersetup{
    colorlinks,
    citecolor=black,
    filecolor=black,
    linkcolor=black,
    urlcolor=black
}

\usepackage[numbered,open]{bookmark} 

\usepackage{fancyhdr}
\usepackage{booktabs}
\usepackage{multirow}
\usepackage[oldcommands,ruled,vlined,linesnumbered]{algorithm2e}

\usepackage{natbib}
 \bibpunct[, ]{(}{)}{,}{a}{}{,}%
 \def\bibfont{\small}%
 \def\BIBand{and}%

\usepackage{float}
\usepackage{graphicx}

\usepackage{amssymb}
\usepackage{amsmath}
\usepackage{tikz}
\usepackage{float}
\usepackage{comment}

\usepackage{lscape}
\usepackage{datetime}
\usepackage[english]{babel}
\usepackage{flafter} 
\usepackage{placeins} 
\usepackage{setspace} 
\usepackage{multirow}
\usepackage{subfigure}

\usepackage{colortbl}

\usepackage[inline]{enumitem}
\newlist{inlinelist}{enumerate*}{1}
\setlist[inlinelist]{label=(\roman*)}

\TheoremsNumberedThrough     

\EquationsNumberedThrough    

\begin{document}

\TITLE{ Service Time Window Design in Last-Mile Delivery }

\ARTICLEAUTHORS{%

\AUTHOR{S. Davod Hosseini}
\AFF{Sobey School of Business, Saint Mary's University, Halifax, Canada\\ \EMAIL{davod.hosseini@smu.ca}}

\AUTHOR{Borzou Rostami}
\AFF{Alberta School of Business, University of Alberta, Edmonton, Canada\\  \EMAIL{borzou@ualberta.ca}}

\AUTHOR{Mojtaba Araghi}
\AFF{Lazaridis School of Business and Economics, Wilfrid Laurier University, Waterloo, Canada\\ \EMAIL{maraghi@wlu.ca}}
}

\ABSTRACT{%
Our study focuses on designing reliable service time windows for customers in a last-mile delivery system to boost dependability and enhance customer satisfaction. To construct time windows for a pre-determined route (e.g., provided by commercial routing software), we introduce two criteria that balance window length and the risk of violation. The service provider can allocate different penalties reflecting risk tolerances to each criterion, resulting in various time windows with varying levels of service guarantee. Depending on the degree of information available about the travel time distribution, we develop two modeling frameworks based on stochastic and distributionally robust optimization. In each setting, we derive closed-form solutions for the optimal time windows, which are functions of risk preferences and the sequence of visits. We further investigate fixed-width time windows, which standardize service intervals, and the use of a policy that allows vehicles arriving before the lower bounds to wait rather than incur a penalty. Next, we integrate service time window design with routing optimization into a unified framework that simultaneously determines optimal routing and time window allocations. We demonstrate the efficacy of our models on a rich collection of instances from well-known datasets. While a small portion of the time windows designed by the stochastic model was violated in out-of-sample tests, the distributionally robust model consistently delivered routes and time windows within the service provider’s risk tolerance. Our proposed frameworks are readily compatible with existing routing solutions, enabling service providers to design time windows aligned with their risk preferences. It can also be leveraged to produce the most efficient routes with narrow time windows that meet operational constraints at controlled levels of service guarantee.
}

\KEYWORDS{On-time last-mile delivery, uncertain travel times, correlation, time window assignment, risk analysis, routing optimization}


\maketitle

\section{Introduction}
\label{sec:intro}

The exponential growth of e-commerce has substantially amplified the volume of business-to-consumer deliveries, witnessing a 54\% surge in domestic parcel delivery volumes across the European Union between 2016 and 2020 \citep{Euro2022}. In the fulfillment of online orders, the so-called \emph{last mile}, the delivery of the order from the carrier to the customer's doorstep, is arguably the least efficient and most expensive stage in the delivery process \citep{macioszek2018first}. This, along with the consistent growth of urban population, could cause more traffic congestion in a city center, creating a negative impact on the well-being of the city economically, socially, and environmentally \citep{deng2021urban}.
On the other hand, however, research has shown that customers prioritize timely delivery and reliability when choosing an online retailer, and delivery performance is a significant factor in customer satisfaction and loyalty \citep{salari2022real}. 
Unlike on-demand delivery (such as food deliveries), where the product should be delivered as soon as possible, most online purchases are scheduled to be delivered during a time window provided by the retailer’s or a third-party delivery system. In spite of that, Deloitte reports a failure rate of 10\%-15\% on the first attempt at home delivery by carrier companies in Spain \citep{Deloitte2020last}. Similar numbers are reported for US, Germany, and UK with an average cost of 15 USD per failed delivery \citep{loqate2022fixing}. Failed deliveries also give rise to enormous amounts of extra emissions \citep{van2015comparative}, a deterioration of service levels \citep{mangiaracina2019innovative}, and collection and delivery point expenses \citep{liu2019assessing}.

Along with inaccurate delivery information, the absence of recipients is one of the major causes of failed deliveries. Some research efforts have focused on predicting the probability of failed delivery attempts using machine learning algorithms \citep{lim2023right}. Many innovations have also been used in recent years to reduce the failed delivery incidents such as decoupling the delivery and pickup by using smart lockers \citep{lyu2022last}. However, such a strategy is not suitable for many deliveries, e.g., perishable products such as flowers, bulky items such as furniture, or packages requiring a signature. For instance, 58\% of big and bulky last mile deliveries were rescheduled \citep{DispatchTrack2022big}. Leaving the item at the door, even when  possible, is not an effective solution as some areas have witnessed an increase in parcel thefts, prompting the police to recommend scheduling deliveries for when someone is at home \citep{WaterlooPolice2023}. Thus, the need to provide customers with reliable delivery time frame still exists. 

A promised delivery time window helps reduce last-mile operations costs as well as customers’ uncertainty and the inconvenience of waiting, and so becomes a lever to manage customer expectations and improve customer satisfaction \citep{cui2020sooner}. However, most businesses that entail pick-up/delivery (including the front runners in this market such as Amazon and FedEx) or service (home services such as installations, repairs, and maintenance as well as home healthcare services such as nursing and physical therapy) only provide their customers with an arrival day or a wide time window within the day. For example, on average less than 10\% of packages are assigned delivery time windows in the dataset released for 2021 Amazon Last Mile Routing Research Challenge \citep{merchan20222021}. Even when the two- to four-hour “Estimated Delivery Window”s are provided by Amazon, “they are not guaranteed and may be subject to change. Deliveries can arrive before or after estimated windows’’ \citep{Amazon2023}.

To optimize resources and streamline operations in last-mile delivery, companies often use approaches to consolidate multiple orders within a specific geographic area. As a common practice in the industry, a single vehicle makes multiple stops to pick up or drop off goods or products from/at different locations along a predetermined route. Although this approach enhances efficiency and significantly reduces the shipping cost of smaller lots, it increases the coordination complexity.
Specifically, it complicates the estimation of delivery time window, as a single disruption along the route will affect all the remaining delivery promises and the error will be amplified as one moves further in the network. Moreover, a late delivery to a customer cannot be offset by an early delivery to another customer as customers expect \emph{on-time} delivery; nearly one in three customers considers early deliveries to be just as bad as late deliveries \citep{DispatchTrack2022big}.

Our study is motivated by such widespread and growing applications of this type of last-mile operation in both pick-up/delivery businesses and the service industry. 
Service providers, visiting multiple customers in each delivery route, desire to design the route and narrow but reliable visit time windows that are guaranteed with some level of confidence.
However, even for a given route, designing such reliable time windows is challenging in real-life situations due to the stochasticity of arrival times to customers. Moreover, travel times among road segments are correlated rather than statistically independent, as for instance congestion on one road is prone to cause congestion on nearby roads
\citep[see, for instance,][]{seshadri2012algorithm,Nicholson201514,Letchford2015steiner,rostami2021branch}. By neglecting correlations among arc travel times, the forecasted travel time variance may be underestimated by up to 75\% \citep{parent2010spatial}. Therefore, addressing the stochasticity of arc travel times and accounting for the correlation among different road segments' travel times are critical in designing any routes and service time windows. 

In this paper, we introduce a novel approach for optimizing last-mile delivery with time window assignment in a network characterized by stochastic and potentially correlated travel times. Our focus is on a service provider tasked with efficiently delivering goods to a set of customers within a predetermined time frame (time budget).
The goal of our study is twofold. First, we aim to design reliable service time windows for a pre-determined route (e.g., provided by any commercial routing software) that accommodates the variability in travel times. This entails designing time windows that minimize possible violations, considering both early and late arrivals at each customer location, and factoring in the service provider's risk preferences/tolerance.
We introduce two distinct modeling frameworks, grounded in stochastic and distributionally robust optimization principles depending on the degree of information available about travel time distribution within the network. 
Second, we extend our two modeling frameworks to optimize routing decisions and time window designs concurrently. Depending on the degree of information available for the underlying travel time distribution, we utilize the previously derived time window characteristics for a given route and develop tractable formulations and sophisticated algorithms for obtaining optimal routes.

\begin{figure}[t]
    \centering
    \includegraphics[scale=.5]{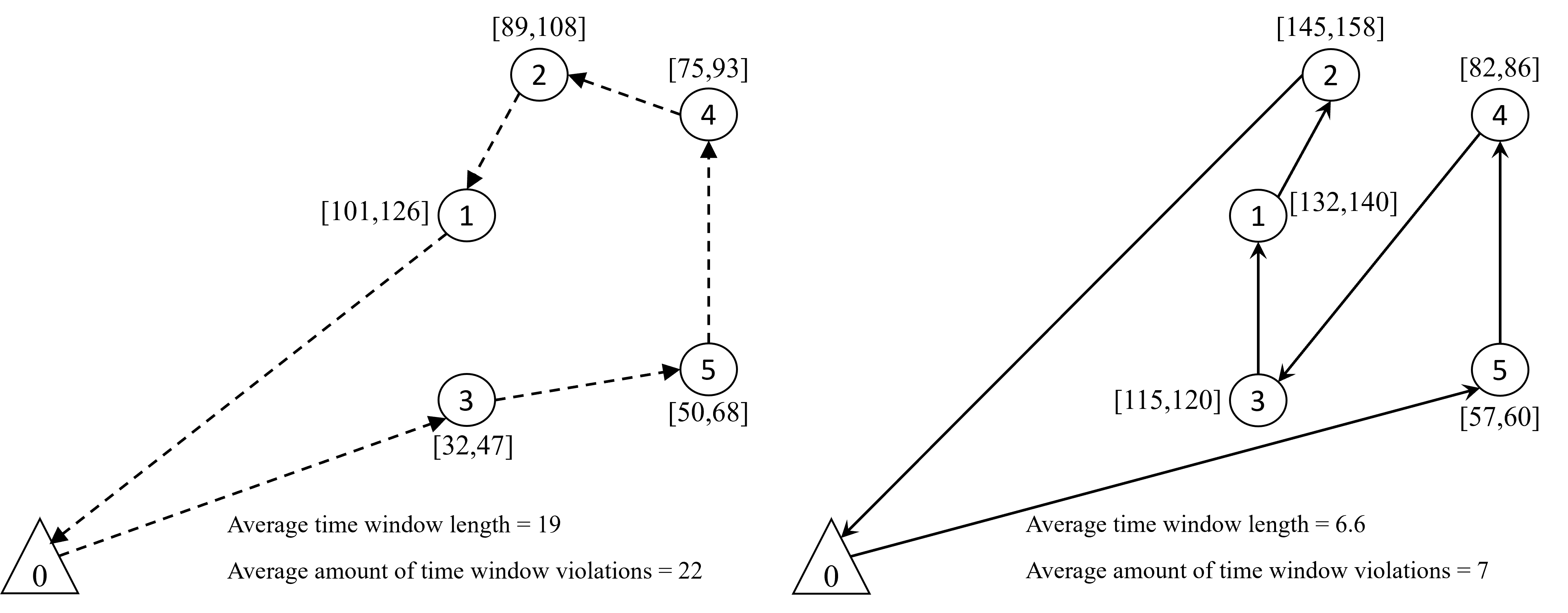}
    \caption{Different time window characteristics (in 1000 training samples) for two example routes}
    \label{fig:two_routes}
\end{figure}

We now present an example to illustrate how our tow goals, time window design and route selection, interact under uncertainty.
Figure 1 compares two ways of serving five customers within a 200-minute time budget, based on a network derived from \cite{solomon1987algorithms}. The \emph{dotted} route is first obtained by optimizing solely for minimal travel time (e.g., using route planning software). We then apply our time window design approach—using a 95\% service guarantee under a stochastic setting with 1000 samples drawn from a Normal distribution—to assign time windows to each customer. Although this route respects the time budget, it yields an average time window length of 19 minutes and an average amount of violations of 22 minutes.
If, instead, we \emph{also} optimize the route choice when assigning time windows, we arrive at the \emph{solid} route. This solution still meets the 200-minute budget but achieves a much shorter average time window length (6.6 minutes) and a much lower average amount of violations (7 minutes). However, improving reliability in this manner comes at the cost of higher overall travel time. 
These two routes illustrate a “tortoise versus hare” trade-off: if minimizing total travel time is the primary goal, one can use existing software to solve the routing problem first and then apply our time window model. Conversely, if tighter windows and higher reliability are desired, our integrated approach to route and time window design is preferable—even though it may slightly increase total travel time.

\subsection{Related Literature}
\label{sec:related_work}
The realm of last-mile delivery encompasses various concepts, challenges, and research opportunities, as explored in \cite{Savelsbergh2016city} and \cite{boysen2021last}. Decision problems in last-mile delivery can be classified into three levels: (i) infrastructure design or setup, (ii) fleet sizing and staffing, and (iii) routing and scheduling. These levels represent a continuum from long-term strategic planning to short-term operational tasks. This study specifically focuses on the third level, which has received extensive attention in the transportation and operations research literature through various formulations of traveling salesman problems (TSPs) and their extensions, such as vehicle routing problems (VRPs) \citep[see, e.g.,][]{laporte2010concise, gendreau2014chapter}. 
However, our study aims to address the unique challenges associated with the explicit design of service times, traditionally considered as given inputs to these problems, thereby extending existing knowledge and methodologies in the field. By incorporating scheduling aspects alongside routing optimization, we strive to uncover novel insights and develop innovative solutions that advance the state-of-the-art in last-mile delivery.

It is worth noting that recent advancements in the field have introduced innovative approaches to address last-mile delivery challenges. For example, shared mobility \citep{qi2018shared}, the integration of predictions with order assignments \citep{liu2021time}, and the utilization of crowdshipping \citep{dayarian2020crowdshipping} or crowdsourcing \citep{fatehi2022crowdsourcing} have all demonstrated the potential to optimize routing and scheduling for efficient customer delivery. However, in what follows, we conduct an extensive literature review on service time design in the context of last-mile delivery. Furthermore, we delve into research studies that explicitly address the correlation between travel times in the design of service times and routing decisions. By synthesizing and analyzing these related works, our goal is to establish a clear understanding of the developments in this field and identify gaps that our research aims to address.

\paragraph{\textbf{Time Window Assignment.}} While routing optimization with given time windows has been extensively studied in the literature \cite[see][for a detailed recent review]{zhang2024generalized},
the concept of assigning time windows to customers in last-mile delivery has emerged more recently and is still an area with limited literature. In these problems, time windows are no longer treated as inputs but become an integral part of the decision-making process.

The first group of papers in this domain focuses on the selection of an endogenous time window, of fixed width, from an exogenous time frame for each customer. \cite{spliet2015time} address this problem in the context of a retail distribution network with demand uncertainty. \cite{spliet2015discrete} study the discrete variant of this problem, where a time window for each customer needs to be selected from a finite set of candidate time windows. \cite{spliet2018time} extend these works by incorporating time-dependent travel times. These papers assume known probability distributions of travel times and propose heuristic solution methods based on the branch-price-and-cut algorithm to solve the routing optimization problem. \cite{subramanyam2018scenario} generalize the work of \cite{spliet2015time} and study problems with scenario-based models of uncertainty in which any operational parameter may be uncertain and  the endogenous time windows may be chosen from either continuous or discrete sets. To handle cases where estimating probability distributions is challenging, \cite{hoogeboom2021robust} propose a robust formulation based on a risk measure and develop a branch-and-cut framework to solve the problem exactly. \cite{martins2019product} extend this type of problem to a product-oriented time window assignment problem, where multi-compartment vehicles are routed to transport products with different temperature requirements to grocery stores within their preferred time windows. They designed an adaptive large neighborhood search method to solve the proposed problem.

The second group of papers focuses on assigning time windows where customers do not impose exogenous time frames. In \cite{jabali2015self}, delivery time windows for customers are determined by the service provider, considering travel time uncertainty modeled by disruption scenarios. However, only one arc is allowed to be disrupted along a route, and the duration of each disruption is assumed to be a discrete random variable with a known probability distribution. They develop a hybrid two-stage tabu search algorithm to find good solutions. \cite{vareias2019assessing} extend this work by allowing multiple arcs to be disrupted simultaneously, with the duration of each disruption considered a continuous random variable. The arrival time at each customer becomes a continuous random variable, depending on the arcs where disruptions have occurred. They propose an adaptive large neighborhood search algorithm, iteratively solving the routing problem and the time window assignment problem. \cite{yu2023time} extend the prior works by considering multiple sources of uncertainties including travel time, service durations, and customer cancellations as well as handling both static and dynamic models by leveraging a rolling horizon approach. In a recent related work, \cite{ulmer2024optimal} consider the situation where a time window is communicated when customers request service during a booking period, which are all served at a later date. In their approach, time window decisions are decoupled from the final service plan: the final routing is determined independently from the assigned time windows. The goal is to minimize the expected time window size across all customers while a chance constraint ensures a high percentage of time windows are satisfied.

Despite the progress in this research domain, none of the aforementioned papers provide time windows with a guarantee of being respected. Furthermore, the existing focus primarily revolves around designing routes and a time window simultaneously. However, our approach takes practical steps by designing time windows that come with a certain level of guarantee provided by the service provider for any predetermined route obtained from any source (e.g., delivery routing software). This allows our approach to be applicable to any delivery company that already has access to routing optimization software. By utilizing our approach, these companies can provide their customers with reliable time windows for deliveries. Furthermore, we seamlessly integrate this time window design with the routing optimization process, ensuring the generation of an optimal route that adheres to the designed time windows.

\paragraph{\textbf{Correlated Travel Times.}} In all the above papers in which assigning time windows to customers has been studied, the travel times of the arcs are assumed to be either deterministic or stochastic and  independent, and, hence, the correlation among them is not considered explicitly. In general, the correlation between travel times has received little attention in the literature of stochastic routing optimization. Over the last four decades, a rich body of stochastic programming models has been developed in the literature to address several variants of routing optimization problems under uncertainty \cite[for an overview see][]{gendreau2014chapter}. Most papers, however, have assumed that uncertainties are independently distributed in order to avoid the tremendous increase in computational complexity. The case of uncertain arc travel times is no exception; most of the works assume the independence of travel times for the sake of simplicity \citep[a recent literature review is provided in][]{rostami2021branch,rajabi2019reliable,bakach2021solving}. However, it contradicts real-life contexts, where for example, the existence of a traffic jam on one road is likely to cause a traffic jam on nearby roads, or poor weather conditions may cause delays on all roads in a certain area \citep{agrawal2012price}. In what follows, we briefly review related papers that address the travel time correlation in routing optimization. However, non of them address the time window assignment, which is the main focus of our study. 

\cite{lecluyse2009vehicle} extend the VRP with time-dependent travel times by adding the standard deviation of the travel time to the objective function to address the variability of the travel times, whose distribution is assumed to be log-normal. They demonstrate the trade-off between the expected travel time and its standard deviation using  simulation, and conclude that as more weight is given to the variability component, the resulting optimal route will take a slightly longer travel time, but will be more reliable. \cite{Letchford2015steiner} study an extension of the Steiner TSP with correlated costs, which follow a multi-variate distribution whose first and second moments are known. Four different integer programming formulations, two quadratic and two linear, were presented to find the efficient tours, in which there is a trade-off between minimizing the expected cost of the tour and minimizing the variance of the cost. 
\cite{rostami2021branch} study the capacitated VRP with stochastic and correlated arc travel times, where the first and the second moments of the travel time probability distributions are assumed to be known, and the correlations are represented by a variance-covariance matrix. Similar to the previous two works in VRP, they seek a trade-off between the expected travel time and its variance (as a measure of the travel time reliability) by adopting a mean-variance approach. The problem is modeled as a binary quadratic program and solved by branch-price-and-cut algorithms. 
They demonstrated that their models can yield routes with a total expected travel time slightly larger than the one of the routes found by the standard VRP, but with significantly less variance.
\cite{bakach2021solving} study a VRP with a makespan objective and stochastic and correlated travel times. The authors present an approach that approximates the expected makespan and the standard deviation based on extreme-value theory. They demonstrate the impact of different correlation patterns and levels of correlation on route planning and report that cost savings of up to 13.76\% can be obtained by considering correlation.

\subsection{Our Contributions}
\label{sec:contrib}
The paper's contributions are summarized as follows.
\begin{inlinelist}
\noindent 
    \item We propose a new approach for designing arrival time windows under uncertainty, using two criteria—window length and on-time arrivals. By factoring in travel-time variability and risk preferences, we minimize violations. We introduce two frameworks—stochastic and distributionally robust optimization—fitted to the available travel time data. This enables service providers to leverage any predetermined route in last-mile delivery and offer more accurate and reliable arrival time windows that enhance customer satisfaction.
    \item For any given routing decision in the stochastic model, we derive closed-form solutions for the time windows and extend our analysis to encompass fixed-width time windows and waiting policies for early arrivals, demonstrating the relationship between these solutions and the service provider’s risk tolerance in terms of specific levels of service guarantee.
    \item To address the correlation between the arcs' travel times in an explicit manner and incorporate distributional ambiguity, we propose a distributionally robust optimization model in which partial distributional information on mean and covariance are used within an ambiguity set. Similar to the stochastic model, for any routing decision, we compute closed-form solutions for the optimal time windows, which allows the service provider to derive managerial insights.
    \item We extend our two modeling frameworks to optimize both routing decisions and time window design concurrently. Depending on the underlying travel time distribution, we utilize the previously derived time window characteristics to develop tractable formulations and sophisticated algorithms for obtaining optimal solutions. 
     \item We conduct extensive computational experiments on benchmark instances to assess our models' potential, robustness, and efficiency. Our findings offer managerial insights on generating reliable time windows tailored to risk tolerances and balancing violation rates with window lengths.
\end{inlinelist}

\subsection{Paper Structure}
\label{sec:paper_structure}
The remainder of this paper is organized as follows. In Section \ref{sec:model}, we describe the proposed criteria to design the service time windows followed by the directions to find their optimal values and characterize their structural properties under both fully and partially known joint distribution of travel times. Section \ref{sec:integrated} presents a combined approach to integrating time windows' design and routing decisions along with decomposition methods to efficiently solve the generated stochastic and distributionally robust models. Section \ref{sec:comput} is allotted to computational study, managerial insights, and analyses. Conclusion and future research directions are ultimately provided in Section \ref{sec:conclusion}.
\section{Service Time Windows Design} 
\label{sec:model}

In this section, we present the service time window design for visiting a set of customers denoted by $V_0$ located at locations $1,2,\ldots, n$. Let location 0 indicate the depot (origin) where drivers are initially deployed and $V=\{0,1,2,…,n\}$ the set of locations that a driver must visit during a time period (usually not more than eight hours). We show by $A$ the set of links between these locations and by $\tilde{t}_{ij}$, $(i,j)\in A$ the uncertain random travel times for traversing them. We assume that the service time at location $j$ is part of $\tilde{t}_{ij}$ for all $(i,j) \in A$.

The aim is to design a priori route with a time window for visiting each customer such that the eventual a posteriori route over any random link travel time is optimal with a certain level of guarantee of not violating the time windows. This can be formally stated as a two-stage procedure: (i) For any pre-determined route (e.g., provided by any commercial routing software), considering the variability in travel times, design reliable service time windows $[\ell^k,u^k]$  for each customer $k\in V_0$ and inform them about the potential time windows of the visit. This entails designing narrowest time windows that minimize possible violations, considering both early and late arrivals at each customer location, and factoring in the service provider's risk preferences/tolerance. (ii) Acknowledging the significance of route selection in the accuracy and reliability of assigned time windows, integrate routing decisions and time window design concurrently.

To design the visit time windows, we need to compute the random arrival time at each customer $k\in V_0$. To this end, we define $y_{ij}^k$ as a binary decision that is equal to 1 if link $(i,j) \in A$ is on the route from depot 0 to customer $k$. Given a routing decision $\boldsymbol{y}^k$, the random arrival time at customer $k$ is then determined as
\begin{flalign*}
\label{eq:arrival_time}
     \tau^k(\boldsymbol{y}^k,\boldsymbol{\tilde{t}})\triangleq {\boldsymbol{y}^k}^\top\boldsymbol{\tilde{t}}=\sum_{(i,j)\in A}y_{ij}^k \tilde{t}_{ij},
\end{flalign*}
taking into account all random link travel times that are part of the route to customer $k$. Based on this random arrival time, we design service time windows $[\ell^k, u^k]$ for each customer $k$ that satisfy the following two criteria C1 and C2:
\begin{align*}
    &\text{C}1: \quad \ell^k\leq \tau^k(\boldsymbol{y}^k,\boldsymbol{\tilde{t}}) \leq u^k,\\
    &\text{C}2: \quad (u^k-\ell^k) \text{ is minimized}.
\end{align*}
The first criterion states that the arrival time to customer $k$ must lie in the proposed time window. A trivial solution that satisfies C1 is to set $\ell^k=0$ and $u^k=+\infty$. However, the service provider is interested in providing its customers with  tight service time windows, in which $u^k-\ell^k$ values are as small as possible. The second criterion addresses this concern. Note that while we allow each customer \(k\) to have a potentially different time window length \(\bigl(u^k - \ell^k\bigr)\) in this section, Extension \ref{subsec:stochastic_fixed} considers an alternative setting in which all time windows share a common, fixed width \(w\), addressing situations where consistent service intervals are required for operational or contractual reasons.

Considering the randomness of the arrival time $\tau^k\bigl(\boldsymbol{y}^k,\boldsymbol{\tilde{t}}\bigr)$ at customer $k$, we address criterion C1 by defining two \emph{random} on-time performance metrics 
\begin{align*}
    & h_\ell^k(\boldsymbol{y}^k, \boldsymbol{\tilde{t}}, \ell^k) = \max\bigl\{\ell^k \;-\; \tau^k\bigl(\boldsymbol{y}^k,\boldsymbol{\tilde{t}}\bigr),\, 0\bigr\}, \quad \text{and} \\
    & h_u^k(\boldsymbol{y}^k, \boldsymbol{\tilde{t}}, u^k) = \max\bigl\{\tau^k\bigl(\boldsymbol{y}^k,\boldsymbol{\tilde{t}}\bigr) \;-\; u^k,\, 0\bigr\}.
\end{align*}
These capture \emph{earliness} (if the arrival time falls below~$\ell^k$) and \emph{tardiness} (if it exceeds~$u^k$), respectively, for a service time window $[\ell^k, u^k]$. 
We next explain how these random metrics $h_\ell^k$ and $h_u^k$ can be \emph{aggregated} into a single, deterministic cost measure under two different modeling paradigms:
\begin{itemize}
    \item Stochastic Programming (Section \ref{sec:stochastic}): In this approach, we assume the travel-time distribution is fully known. Consequently, we take the \emph{expected} values of the random earliness and tardiness metrics $h_\ell^k$ and $h_u^k$ with respect to that known distribution. This yields $\mathcal{H}_\ell^k\bigl(\boldsymbol{y}^k,\ell^k\bigr)$ and $\mathcal{H}_u^k\bigl(\boldsymbol{y}^k,u^k\bigr)$, which summarize \emph{average} early- and late-arrival behavior, respectively, for customer $k$.
    \item Distributionally Robust Optimization (Section \ref{sec:robust}): Here, only partial information (e.g., the mean and covariance) of the travel-time distribution is available. Instead of a single expectation, we adopt a \emph{worst-case} viewpoint and evaluate $h_\ell^k$ and $h_u^k$ under the most adverse distribution consistent with that partial information. The resulting $\mathcal{H}_\ell^k\bigl(\boldsymbol{y}^k,\ell^k\bigr)$ and $\mathcal{H}_u^k\bigl(\boldsymbol{y}^k,u^k\bigr)$ thus reflect \emph{worst-case} on-time performance for customer $k$.
\end{itemize}
\noindent
By combining these measures with the time window width penalty (criterion C2), we obtain the overall service time window cost for customer $k$ as
\begin{equation}
\label{eq:Hfunction}
    \mathcal{H}^k(\boldsymbol{y}^k,\ell^k, u^k) = a_w^k(u^k-\ell^k) + a_{\ell}^k\,\mathcal{H}_{\ell}^k(\boldsymbol{y}^k, \ell^k) + a_u^k\mathcal{H}_u^k(\boldsymbol{y}^k, u^k), 
\end{equation}
whose specific form will depend on whether we evaluate $\mathcal{H}_\ell^k$ and $\mathcal{H}_u^k$ via expectation (Section \ref{sec:stochastic}) or via a robust supremum (Section \ref{sec:robust}).
The weights $a_w^k, a_u^k, a_{\ell}^k \in (0,1]$, $\forall k\in V_0$ in \eqref{eq:Hfunction} are real penalty parameters representing the importance of each component  for the service provider. More precisely, $a_w^k$ is the penalty associated with the length (width) of the time window, corresponding to criterion C2, while  $a_{\ell}^k$ and $a_u^k$ are penalties associated with the earliness and tardiness metrics, respectively, corresponding to criterion C1. 

In the definition of cost function \eqref{eq:Hfunction}, we assume that if the vehicle arrives earlier than the assigned start time, the service provider will not wait and will start the service upon the arrival. This case is suitable for routing in dense urban areas where parking spaces are extremely limited \citep{jaillet2016routing} or 
when operating on a tight schedule or under significant travel‐time uncertainty. By allowing early arrivals to begin service immediately, we preserve slack in the schedule that can be used to absorb delays elsewhere. If we were to eliminate earliness (i.e., always wait until the assigned start time), we would lose this buffer, increasing the risk (and cost) of tardiness at subsequent customers—particularly in contexts with limited time budgets and high penalties for running late.
However, this assumption can be relaxed to allow waiting if a vehicle arrives before $\ell^k$, as discussed in Extension \ref{sec:waiting_allowed}.

Given a routing decision $\boldsymbol{y}^k$, the primary objective of the service provider is to design a service time window $[\ell^k, u^k]$ for each customer $k$ that minimizes the service time window design cost $\mathcal{H}^k$ in \eqref{eq:Hfunction}. This can be achieved by solving the following optimization problem for each customer $k$:
\begin{align} 
    \label{main_subproblem}\text{SP$^k(\boldsymbol{y}^k$)}:\quad \min_{\ell^k, u^k} &~ \Big\{\mathcal{H}^k(\boldsymbol{y}^k, \ell^k, u^k):~
   0\leq \ell^k \leq u^k \Big\}.
\end{align}
In the subsequent sections, \ref{sec:stochastic} and \ref{sec:robust}, we present the exact formulations to compute $\mathcal{H}_\ell^k$, $\mathcal{H}_u^k$, and consequently $\mathcal{H}^k$ in each setting of a fully known distribution and a partially known distribution of random travel times, respectively, for a given routing decision $\boldsymbol{y}^k$.

\subsection{Design under Fully Known Distribution}
\label{sec:stochastic}

Consider the random vector $\boldsymbol{\tilde{t}}$ of link travel times with continuous probability density function $p(\boldsymbol{\tilde{t}})$ that follows a continuous distribution $\mathrm{P}$. For a given routing decision $\boldsymbol{y}^k$, the arrival time $\tau^k(\boldsymbol{y}^k,\boldsymbol{\tilde{t}})$ at each customer $k \in V_0$ is a random variable with the distribution induced by that of $\boldsymbol{\tilde{t}}$.
We show by  $F^k(\boldsymbol{y}^k, \epsilon^k)$ the cumulative distribution function of arriving time at customer $k$, i.e.,  $F^k(\boldsymbol{y}^k, \epsilon^k)=\text{Pr}(\tau^k(\boldsymbol{y}^k,\boldsymbol{\tilde{t}})\leq \epsilon^k)$ with $\epsilon^k$ being a positive real number.
Knowing the probability distribution of $\tau^k(\boldsymbol{y}^k,\boldsymbol{\tilde{t}})$, we define $\mathcal{H}_{\ell}^k$ and $\mathcal{H}_{u}^k$  to quantify the expected earliness and tardiness at customer $k$, respectively, as follows 

\begin{align}
    &\mathcal{H}_{\ell}^k(\boldsymbol{y}^k,\ell^k) \triangleq \mathbb{E}_{\mathrm{P}}\left[\left(\ell^k-\tau^k(\boldsymbol{y}^k,\boldsymbol{\tilde{t}})\right)^+\right]  \label{eq:h_l_stochastic}, \\  
    &\mathcal{H}_u^k(\boldsymbol{y}^k, u^k) \triangleq \mathbb{E}_{\mathrm{P}}\left[\left(\tau^k(\boldsymbol{y}^k,\boldsymbol{\tilde{t}})-u^k\right)^+\right], \label{eq:h_u_stochastic}
\end{align}
where 
$\left(.\right)^+ = \max \{., 0\}$. 
Plugging \eqref{eq:h_l_stochastic} and \eqref{eq:h_u_stochastic} into the cost function \eqref{eq:Hfunction}, we can derive an optimal time window solution for each customer $k \in V_0$ with a given routing decision $\boldsymbol{y}^k$ by solving the stochastic optimization problem gained in \eqref{main_subproblem}.
    
\begin{proposition}
\label{Proposition1}
With \(\mathcal{H}_{\ell}^k\) and \(\mathcal{H}_{u}^k\) defined as in \eqref{eq:h_l_stochastic} and \eqref{eq:h_u_stochastic}, the function \(\mathcal{H}^k\) in \eqref{eq:Hfunction} is convex and continuously differentiable with respect to \(\ell^k\) and \(u^k\). Moreover, we have:
\begin{align*}
    &\frac{\partial}{\partial \ell^{k}} \mathcal{H}^k(\boldsymbol{y}^k, \ell^k, u^k) = -a_w^k + a_\ell^k \,F^k\left(\boldsymbol{y}^{k}, \ell^{k}\right), \quad \text{and} \\
    &\frac{\partial}{\partial u^{k}} \mathcal{H}^k(\boldsymbol{y}^k, \ell^k, u^k) = a_w^k + a_u^k \Big(F^k\left(\boldsymbol{y}^{k}, u^{k}\right) - 1\Big).
\end{align*}
\end{proposition}
\proof{Proof.} See Appendix A.
\hfill $\square$

Using Proposition \ref{Proposition1} and considering the fact that the feasible set of model \eqref{main_subproblem} is linear, this model is a convex optimization problem, and, hence, a local minimum is the global one. In what follows, we use the results of Proposition \ref{Proposition1} to characterize the penalty parameters $\boldsymbol{a}^k$ and derive a closed-form solution for the optimal time windows. To this end, we consider model \eqref{main_subproblem} and ignore the condition $0\leq \ell^k \leq u^k$ for now. Later, we show the optimal solution we construct will satisfy this condition. The stationarity conditions of $\mathcal{H}^k$ with respect to $\ell^k$ and $u^k$ 
force the local minimizer $(\bar{\ell}^k, \bar{u}^k)$ of $\mathcal{H}^k$ (and hence the global minimizer because of the convexity) to satisfy the following equations:
\begin{align}
  \label{eq:condition1}  &F^k\left(\boldsymbol{y}^{k}, \bar {\ell}^{k}\right) = \text{Pr}\big(\tau^k(\boldsymbol{y}^k,\boldsymbol{\tilde{t}})\leq \bar {\ell}^{k}\big)= \frac{a_w^k}{a_\ell^k}, \quad \text{and}\\
\label{eq:condition2}    &F^k\left(\boldsymbol{y}^{k}, \bar{u}^{k}\right) = \text{Pr}\big(\tau^k(\boldsymbol{y}^k,\boldsymbol{\tilde{t}})\leq \bar {u}^{k}\big)=1-\frac{a_w^k}{a_u^k}.
\end{align}

\begin{figure}[h]
    \includegraphics[scale=.7]{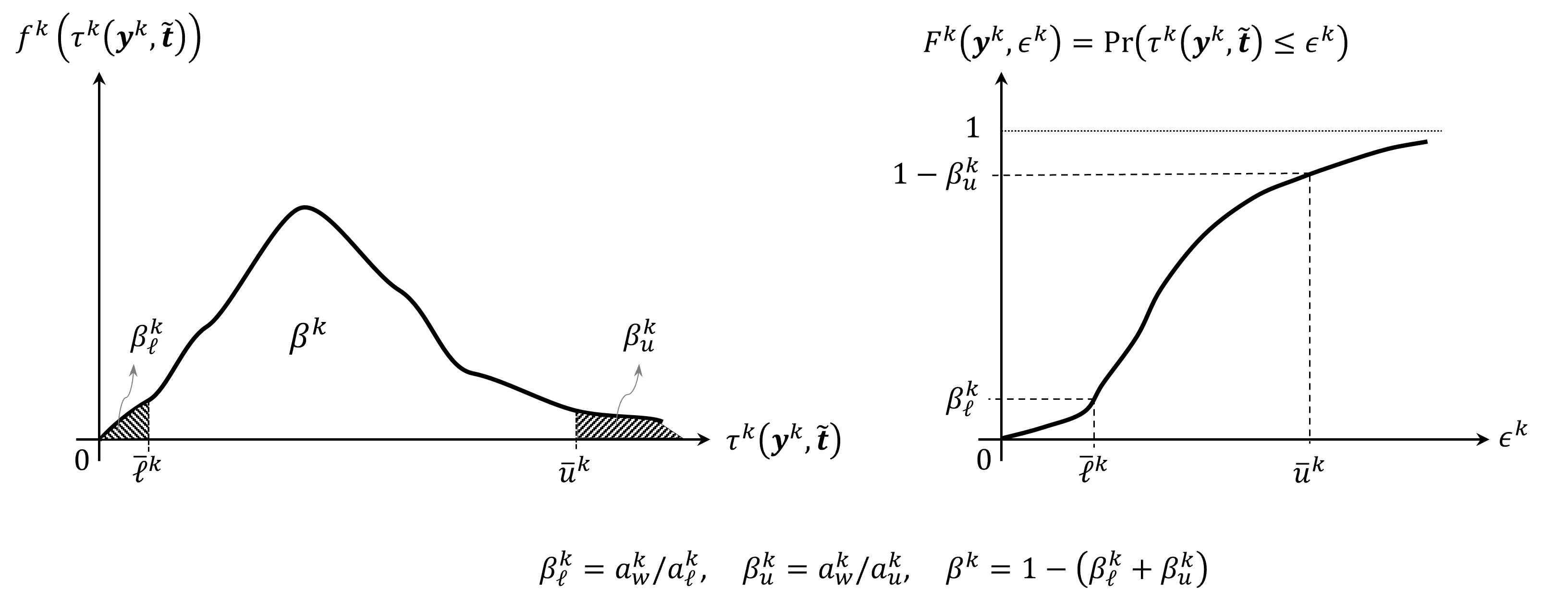}
    \caption{Confidence levels in designing the optimal service time window $[\bar{\ell}^k, \bar{u}^k]$ for customer $k \in V_0$}
    \label{fig:confidence_level}
\end{figure}

These equations reveal several important insights about the optimal time windows. In particular, if $a_w^k \ll min\{a_\ell^k, a_u^k\}$ (e.g., $a_w^k \rightarrow 0$), the optimal time window $[\bar{\ell}^k, \bar{u}^k]$ is equal to $[0,+\infty]$. This is true because $F\left(\boldsymbol{y}^{{k}}, \epsilon^{k}\right)$ is continuous and non-decreasing in $\epsilon^{k}$ with limit 1 as $\epsilon^{k} \rightarrow+\infty$ and limit 0 as $\epsilon^{k} \rightarrow 0$. This immediate result confirms our initial observation of the necessity of both conditions C1 and C2. More importantly,
   because $F\left(\boldsymbol{y}^{{k}}, \epsilon^{k}\right)$ is continuous and non-decreasing in $\epsilon^{k}$, we have
\begin{align*}
    &\text{Pr}\big(\tau^k(\boldsymbol{y}^k,\boldsymbol{\tilde{t}})\leq \ell ^{k}\big)\leq \frac{a_w^k}{a_\ell^k} & \text{for all}~\ell^{k}\leq \bar {\ell}^{k}, \\
    &\text{Pr}\big(\tau^k(\boldsymbol{y}^k,\boldsymbol{\tilde{t}})\geq {u}^{k}\big)\leq \frac{a_w^k}{a_u^k} & \text{for all}~u^{k}\geq \bar {u}^{k},
\end{align*}
meaning that, $\bar{\ell}^k$ is the largest one among the lower bounds for which $\text{Pr}(\tau^k(\boldsymbol{y}^k,\boldsymbol{\tilde{t}})\leq \ell ^{k})\leq \frac{a_w^k}{a_\ell^k}$, and $\bar{u}^k$ is the smallest one among the upper bounds for which $\text{Pr}(\tau^k(\boldsymbol{y}^k,\boldsymbol{\tilde{t}})\geq {u}^{k})\leq \frac{a_w^k}{a_u^k}$ (see Figure \ref{fig:confidence_level}). 
Furthermore, if we assume $a_w^k/a_\ell^k + a_w^k/a_u^k \leq 1$; then, equations \eqref{eq:condition1} and \eqref{eq:condition2} imply that $F^k\left(\boldsymbol{y}^{k}, \bar {\ell}^{k}\right) \leq F^k\left(\boldsymbol{y}^{k}, \bar {u}^{k}\right)$, which in turn yields $\bar{\ell}^k \leq \bar{u}^k$.
This assumption stands in contrast to another extreme case (in addition to the case $a_w^k \to 0$ discussed previously). If \(a_w^k\) grows very large relative to \(a_\ell^k\) and \(a_u^k\) (i.e., \(a_w^k \to +\infty\)), then any nonzero window width becomes prohibitively expensive, forcing the window to shrink toward a single point and thus incurring earliness or tardiness for all arrivals. By ensuring \(a_w^k / a_\ell^k + a_w^k / a_u^k \le 1\), we obtain a moderate penalty structure that avoids this degenerate outcome, which along with C2 yields a nontrivial, finite time window in the optimal solution.
The following proposition summarizes the main results. 

\begin{proposition}
\label{Proposition2}
For each customer $k\in V_0$ with fixed routing decision $\boldsymbol{y}^k$, assuming that $a_w^k/a_\ell^k + a_w^k/a_u^k \leq 1$, the optimal service time window $[\bar{\ell}^k, \bar{u}^k]$ of problem \eqref{main_subproblem} is given as:

\begin{align*}
     \bar{\ell}^k = \max \Big\{ \ell^k : ~\operatorname{\text{Pr}}\left(\ell^k \leq \tau^k(\boldsymbol{y}^k, \boldsymbol{\tilde{t}}) \leq \bar{u}^k\right) \geq 1-\frac{a_w^k}{a_{\ell}^k}-\frac{a_w^k}{a_u^k} \Big\}, \\
     \bar{u}^k = \min \Big \{u^k : ~\operatorname{\text{Pr}}\left(\bar{\ell}^k \leq \tau^k(\boldsymbol{y}^k, \boldsymbol{\tilde{t}}) \leq u^k\right) \geq 1-\frac{a_w^k}{a_{\ell}^k}-\frac{a_w^k}{a_u^k} \Big\}.
\end{align*}

\end{proposition} \vspace{-1em}
As illustrated in Figure \ref{fig:confidence_level}, the immediate corollary of Proposition \ref{Proposition2} is that, by providing the service time window $[\bar{\ell}^k, \bar{u}^k]$ to customer $k$, the service provider has a confidence level of $1 - \beta_\ell^k=1- a_w^k/a_\ell^k$ in arriving at the customer's location after $\bar{\ell}^k$, and a confidence level of $1 - \beta_u^k=1-a_w^k/a_u^k$ in arriving at the customer's location before $\bar{u}^k$. This will yield a joint confidence level $\beta^k = 1 - (\beta_\ell^k + \beta_u^k)$ in arriving at the customer's location within the time window $[\bar{\ell}^k, \bar{u}^k]$. Hence, the choice of penalty parameters $\boldsymbol{a}^k$ expresses how the service provider is concerned about the length of the time window as well as the too early and/or too late arrivals at the customers. For example, having $a_\ell^k < a_u^k$, the late violation rate is expected to be less than the early violation rate, i.e, $\beta_\ell^k > \beta_u^k$. Note that in real-world applications, decision makers desire a high confidence level (usually $\beta^k \geq 90\%$). Hence, it is realistic to set the penalty parameters $\boldsymbol{a}^k$ in the rest of the paper such that $\frac{a_w^k}{a_u^k}< 0.5$ and $\frac{a_w^k}{a_{\ell}^k}  < 0.5$.

\paragraph{\textbf{Structural Properties under Sample Average Approximation.}}
Although Proposition \ref{Proposition2} provides the structure of the optimal time window for each customer with a certain level of guarantee, deriving optimal time windows requires complete knowledge of the joint distribution 
$\mathrm{P}$ and involves performing integration operations. To address this, one may employ the sample average approximation (SAA) scheme.
Given a set of $Q$ samples $\boldsymbol{t}^{[1]}, \boldsymbol{t}^{[2]}, ...,\boldsymbol{t}^{[Q]}$ of travel time vectors generated from the probability distribution of $\tilde{\boldsymbol{t}}$ with density $p(\tilde{\boldsymbol{t}})$, we can approximate the arrival time at each customer $k$ for each sample $q$ as $\tilde \tau^k\left(\boldsymbol{y}^{k}, \boldsymbol{t}^{[q]}\right)=\sum_{(i, j) \in A}t_{ij}^{[q]} y_{i j}^{k}$ which ultimately yields the following approximations \eqref{eq:approx_l} and \eqref{eq:approx_u} for the  $\mathcal{H}_{\ell}^k(\boldsymbol{y}^k, \ell^k)$ and $\mathcal{H}_{u}^k(\boldsymbol{y}^k, u^k)$ given in \eqref{eq:h_l_stochastic} and \eqref{eq:h_u_stochastic}, respectively: 
 \begin{align}
 \label{eq:approx_l}
      &\tilde{\mathcal{H}}_{\ell}^k(\boldsymbol{y}^k, \ell^k) = \frac{1}{Q} \sum_{q=1}^{Q}\left(\ell^{k}-\tilde \tau^k(\boldsymbol{y}^{k}, \boldsymbol{t}^{[q]})\right)^{+}, \text{ and} \\
 \label{eq:approx_u}
    &\tilde{\mathcal{H}}_{u}^k(\boldsymbol{y}^k, u^k) = \frac{1}{Q} \sum_{q=1}^{Q}\left(\tilde \tau^k(\boldsymbol{y}^{k}, \boldsymbol{t}^{[q]})-u^{k}\right)^{+}.
 \end{align}

We can linearize  the nonlinear terms in $\tilde{\mathcal{H}}_{\ell}^k$ and $\tilde{\mathcal{H}}_{u}^k$ through introducing nonnegative auxiliary variables $\boldsymbol{v}_1$ and $\boldsymbol{v}_2$, respectively, and adding the following constraints for each sample $q$:
\begin{align} 
    \label{eq:Hl_v1} & \ell^{k}-\sum_{(i, j) \in A}\left(t_{i j}^{[q]} y_{i j}^{k}\right) \leq v_{1}^{k[q]} \\
    \label{eq:Hu_v2} & \sum_{(i, j) \in A}\left(t_{i j}^{[q]} y_{i j}^{k}\right)-u^{k} \leq v_{2}^{k[q]} \\
    \label{eq:luv_sign} & 0\leq \ell^{k}\leq u^{k}, \quad v_{1}^{k[q]}, v_{2}^{k[q]} \geq 0,
\end{align}
which leads to the approximation of service time window design problem \eqref{main_subproblem}:

\begin{equation}
\label{eq:SP(SM)}
  \tilde{\text{SP}}^k(\boldsymbol y^k): \quad \min_{\ell^k, u^k, \boldsymbol{v_{1}}, \boldsymbol{v_{2}}} \Big\{a_w^k(u^k - \ell^k) + \frac{a_\ell^k}{Q} \sum_{q=1}^{Q} v_{1}^{k[q]} + \frac{a_u^k}{Q} \sum_{q=1}^{Q} v_{2}^{k[q]}: ~\eqref{eq:Hl_v1},\eqref{eq:Hu_v2}, \eqref{eq:luv_sign}\Big\}.
\end{equation}

Here, we show how to derive  a closed-form for the optimal time window $[\bar{\ell}^k, \bar{u}^k]$ in $\tilde{\text{SP}}^k(\boldsymbol y^k)$. To do so, we sort all observed arrival times to customer $k$, $\tilde \tau^{k[q]}=\sum_{(i, j) \in A} t_{i j}^{[q]} {y}_{i j}^{k}$, to obtain a permutation $\Lambda$  such that 
 $ \tilde \tau^{k[\Lambda_1]} \leq  \tilde \tau^{k[\Lambda_2]} \leq \ldots \leq   \tilde \tau^{k[\Lambda_Q]}$. For ease of presentation, we let $\Lambda=(1,2,\ldots Q)$. We then define critical sample indices $P_1^{k}, P_2^k\in \{1, 2, \ldots, Q\}$ such that
\begin{align}
\label{eq:Qorder1}
   & \frac{1}{Q}\,\sum_{q=1}^{P_1^{k}-1} a_\ell^k < a_w^k \leq   \frac{1}{Q}\,\sum_{q=1}^{P_1^{k}} a_\ell^k, \quad \text{ and} \\
    \label{eq:Qorder2}
    & \frac{1}{Q}\,\sum_{q=P_2^{k}+1}^{Q} a_u^k < a_w^k \leq   \frac{1}{Q}\,\sum_{q=P_2^{k}}^{Q} a_u^k.
\end{align}

 We are now ready to present the structural properties of the approximate optimal time windows. 
 
\begin{proposition}
\label{Proposition3}
For each customer $k \in V_0$, the optimal time window for the approximate model $\tilde{\text{SP}}^k(\boldsymbol y^k)$ in \eqref{eq:SP(SM)} is given as
\begin{align*}
&{\bar \ell^k}= \tilde \tau^{k[P_1^{k}]}=\sum_{(i, j) \in A} t_{ij}^{[P_1^{k}]} {y}_{ij}^{k},
&{\bar u^k}= \tilde \tau^{k[P_2^{k}]}=\sum_{(i, j) \in A} t_{ij}^{[P_2^{k}]} {y}_{ij}^{k}.
\end{align*}
Moreover, the optimal values for the auxiliary variables  $\boldsymbol{v}_1$ and $\boldsymbol{v}_2$ are
\begin{align*}
    &\bar v_1^{k[q]} = \begin{cases}
 \tilde \tau^{k[P_1^{k}]}-\tilde \tau^{k[q]}  & q=1,\ldots, P_1^{k}-1;\\
0& q=P_1^{k}, \ldots, Q.
 \end{cases}
    &\bar v_2^{k[q]} = \begin{cases}
0& q=1, \ldots, P_2^{k};\\
\tilde \tau^{k[q]}-\tilde \tau^{k[P_2^{k}]}  & q=P_2^{k}+1,\ldots, Q.
 \end{cases}
\end{align*}
\end{proposition}

\proof{Proof.} See Appendix B.
\hfill $\square$

Proposition \ref{Proposition3} produces an important insight into the approximate time windows.
The approximate time window solutions imply that samples $q=1,\ldots, P_1^{k}-1$ and samples $q=P_2^{k}+1,\ldots, Q$ violate the assigned time window $[{\bar{\ell}^k}, {\bar{u}^k}]$ by arriving too early or too late at customer $k$, respectively. That is, we can derive the violation rates $\frac{P_1^{k}-1}{Q}$ and $\frac{Q-P_2^{k}}{Q}$ for early and late arrivals, respectively. These violations are represented as penalties through the optimal values of $\boldsymbol{v}_1$ and $\boldsymbol{v}_2$ variables. Whereas there is no penalty for the samples in which arrival times occur within the time window (i.e., ${\bar{v}_1^k}=0$ for samples $q=P_1^{k}, \ldots, Q$ arriving after the earliest time ${\bar{\ell}^k}$, and  ${\bar{v}_2^k}=0$ for samples $q=1, \ldots, P_2^{k}$ arriving before the latest time ${\bar{u}^k}$), the penalty for the samples with time window violation is the difference between the arrival time ($\tilde \tau^{k[q]}$) and the lower bound ($\tilde \tau^{k[P_1^{k}]}$) or the upper bound ($\tilde \tau^{k[P_2^{k}]}$).

More importantly, using \eqref{eq:Qorder1} and \eqref{eq:Qorder2}, we can derive $\frac{P_1^{k}-1}{Q} < \frac{a_w^k}{a_\ell^k}$ and $\frac{Q-P_2^{k}}{Q} < \frac{a_w^k}{a_u^k}$, respectively. 
That is, $\beta_\ell^k=\frac{a_w^k}{a_\ell^k}$ and $\beta_u^k=\frac{a_w^k}{a_u^k}$ can be interpreted as the risk tolerance of the service provider on either side of the time window. Therefore, the samples' early violation rate $\frac{P_1^{k}-1}{Q}$ and late violation rate $\frac{Q-P_2^{k}}{Q}$ are less than the service provider's maximum acceptable violation rates $\beta_\ell^k$ and $\beta_u^k$, respectively.  We investigate these insights via numerical experiments in Section \ref{sec:comput}.

Before studying the time window design under partially known distribution, we examine how our assumptions on variable-length time windows and the decision not to wait when arriving before the time windows' lower bounds affect the structure of our model and the corresponding results in the following two subsections, respectively.

\subsubsection{Extension 1: Fixed-Length Time Windows.}
\label{subsec:stochastic_fixed}

To address scenarios where service intervals must remain consistent for operational or contractual reasons, this extension modifies the original model of Section \ref{sec:stochastic} by enforcing a fixed length \(w\) for all customer time windows. Each customer \(k \in V_0\) is now assigned a window \([ \ell^k, \ell^k + w ]\), where \(w \geq 0\) is a common decision variable. The revised cost function combines the fixed-width penalty, expected earliness, and expected tardiness:  
\begin{align}
    & \mathcal{H}^k(\boldsymbol{y}^k, \ell^k, w) = a_w w + a_{\ell}^k~\mathbb{E}_{\mathrm{P}}\left[\left(\ell^k-\tau^k(\boldsymbol{y}^k,\boldsymbol{\tilde{t}})\right)^+\right] + a_u^k~\mathbb{E}_{\mathrm{P}}\left[\left(\tau^k(\boldsymbol{y}^k,\boldsymbol{\tilde{t}})- (\ell^k + w)\right)^+\right]. \label{eq:Hfunction_appx}
\end{align} 
The optimization problem \eqref{main_subproblem_appx} then jointly determines \(\ell^k\) and \(w\) to minimize this cost for a given routing decision $\boldsymbol{y}^k$:  
\begin{align}
    & \text{SP$^k(\boldsymbol{y}^k$)}:\quad \min_{\ell^k, w}~\Big\{\mathcal{H}^k(\boldsymbol{y}^k, \ell^k, w):~ \ell^k \geq 0, w \geq 0 \Big\}. \label{main_subproblem_appx}   
\end{align}
Analogous to the discussions in Propositions \ref{Proposition1} and \ref{Proposition2}, one can show that the first‐order conditions at the local optimum $(\bar{\ell}^k, \bar{w})$—and hence the global optimum because $\mathcal{H}^k(\boldsymbol{y}^k, \ell^k, w)$ is again convex in $(\ell^k, w)$—link the cumulative distribution function $F^k$ of random arrival times to the penalty parameters $\boldsymbol{a}$:
\begin{align}
    \label{eq:condition1_appx_2}  &F^k\left(\boldsymbol{y}^{k}, \bar {\ell}^{k}\right) = \text{Pr}\left(\tau^k(\boldsymbol{y}^k,\boldsymbol{\tilde{t}})\leq \bar {\ell}^{k}\right) = \frac{a_w}{a_\ell^k}, \\
    \label{eq:condition2_appx_2}  &F^k\left(\boldsymbol{y}^{k}, \bar{\ell}^{k}+ \bar{w} \right) = \text{Pr}\left(\tau^k(\boldsymbol{y}^k,\boldsymbol{\tilde{t}})\leq \bar{\ell}^{k} + \bar{w}\right) = 1 - \frac{a_w}{a_u^k}.
\end{align}
This ensures $\bar{\ell}^k$ and $\bar{w}$ balance the trade-off between window length and violation probabilities. In fact, by assigning the time window $[\bar{\ell}^k, \bar{\ell}^k + \bar{w}]$ to each customer $k$, the service provider attains a joint confidence level of $1 - (a_w/a_\ell^k + a_w/a_u^k)$ for on-time arrival, provided $a_w/a_\ell^k + a_w/a_u^k \leq 1$.

To bypass the computational challenges of integrating over the full travel-time distribution \(\mathrm{P}\), one may employ SAA. Given \(Q\) travel-time samples \(\{\boldsymbol{t}^{[q]}\}_{q=1}^Q\), the expected earliness and tardiness penalties are approximated by $\frac{1}{Q} \sum_{q=1}^Q \left( \ell^k - \tilde{\tau}^k (\boldsymbol{y}^k, \boldsymbol{t}^{[q]}) \right)^+$ and $\frac{1}{Q} \sum_{q=1}^Q \left( \tilde{\tau}^k (\boldsymbol{y}^k, \boldsymbol{t}^{[q]}) - (\ell^k + w) \right)^+$, respectively, where each $(.)^+$ term turns into a linear function once the auxiliary variables \(v_1^{k[q]}, v_2^{k[q]}\) are introduced for each sample $q$, respectively, enforcing:  
\begin{align} 
    \label{eq:Hl_v1_appx_2} & \ell^{k}-\sum_{(i, j) \in A}\left(t_{i j}^{[q]} y_{i j}^{k}\right) \leq v_{1}^{k[q]}, \\
    \label{eq:Hu_v2_appx_2} & \sum_{(i, j) \in A}\left(t_{i j}^{[q]} y_{i j}^{k}\right)- (\ell^{k}+w) \leq v_{2}^{k[q]}, \\
    \label{eq:luv_sign_appx_2} & \ell^{k}, w, v_{1}^{k[q]}, v_{2}^{k[q]} \geq 0.
\end{align}
This will transform problem \eqref{main_subproblem_appx} into a tractable linear program:  
\begin{equation}
\label{eq:SP(SM)_appx_2}
  \tilde{\text{SP}}^k(\boldsymbol y^k): \quad \min_{\ell^k, w, \boldsymbol{v_{1}}, \boldsymbol{v_{2}}} \Big\{ a_w w + \frac{a_\ell^k}{Q} \sum_{q=1}^{Q} v_{1}^{k[q]} + \frac{a_u^k}{Q} \sum_{q=1}^{Q} v_{2}^{k[q]}: ~\eqref{eq:Hl_v1_appx_2},\eqref{eq:Hu_v2_appx_2}, \eqref{eq:luv_sign_appx_2}\Big\}.
\end{equation}

\subsubsection{Extension 2: Waiting before Time Windows’ Lower Bounds.}
\label{sec:waiting_allowed}

To accommodate waiting at each customer $k \in V_0$ when the vehicle arrives before $\ell^k$ (the nominal lower bound), let  $T^k$ represent the \emph{actual} service start time at customer $k$, which can be calculated recursively as
\begin{equation}
\label{eq:recursive_1}
    T^k = \max \left\{T^{i_{k-1}} + \tilde{t}_{i_{k-1}k}, ~\ell^k \right\},
\end{equation}
where $i_{k-1}$ is the customer served right before customer $k$ ($y^k_{i_{k-1}k} = 1$). Therefore, if the vehicle's arrival time at customer $k$ is earlier than $\ell^k$, the service will be postponed and start at $\ell^k$. Otherwise, service commences immediately upon arrival at time $T^{i_{k-1}} + \tilde{t}_{i_{k-1}k}$. Assuming the partial route to reach customer $k$ is in the sequence $\{i_0, i_1, i_2, ..., i_{k-1}, k\}$, where $i_0 = 0$, the service start times at the customers served before customer $k$ can be determined recursively using \eqref{eq:recursive_1} as follows
\begin{align*}
    & T^{i_0} = 0, \\
    & T^{i_1} = \max \{ \tilde{t}_{0i_1}, \ell^{i_1} \}, \\
    & T^{i_2} = \max \{T^{i_1} + \tilde{t}_{i_1i_2}, \ell^{i_2} \}, \\
    & \vdots \\
    & T^{i_{k-1}} = \max \{T^{i_{k-2}} + \tilde{t}_{i_{k-2}i_{k-1}}, \ell^{i_{k-1}} \}.
\end{align*}
That is, any waiting or delay at earlier stops propagates forward, which results in
\begin{align*}
    T^k = \max \Bigg\{ \max \bigg\{ \cdots \max \Big\{ \max \{ \tilde{t}_{0i_1}, \ell^{i_1} \} + \tilde{t}_{i_1i_2}, \ell^{i_2} \Big\} \cdots, \ell^{i_{k-1}} \bigg\} + \tilde{t}_{i_{k-1}k}, ~\ell^k \Bigg\},
\end{align*}
whose nested sequence of maxima can be unrolled into a single maximum over all nodes along the route:
\begin{align}
\label{eq:recursive_2}
    T^k = \max_{r \in \{i_0, i_1, i_2, \dots, i_{k-1}, k\}} 
    \left\{
    \ell^r + \sum_{a \in \big\{(r, i_{\hat{r}+1}), (i_{\hat{r}+1}, i_{\hat{r}+2}), \dots, (i_{k-1}, k)\big\}} \tilde{t}_a
    \right\},
\end{align}
where $r = i_{\hat{r}}$ indicates a node along the partial route from node 0 to customer $k$.

This way, the time window design cost associated with customer $k$ in \eqref{eq:Hfunction} will transform to 
\begin{equation}
\label{eq:Hfunction_ext}
   \mathcal{H}^k(\boldsymbol{y}^k, \ell^k, u^k) = a_w^k(u^k - \ell^k) + a_u^k ~ \mathbb{E}_{\mathrm{P}}\left[\left(T^k - u^k\right)^+\right], 
\end{equation}
where the second term penalizes the risk of arriving late (tardiness). When the vehicle arrives early, a larger $\ell^k$ forces service to begin later—narrowing the effective window $(u^k - T^k)$ and thereby reducing the associated width penalty. However, a narrower effective window provides less slack to absorb travel time variations, which can increase the risk (and cost) of tardiness. In contrast, setting $\ell^k$ very small allows the service to begin earlier if the vehicle arrives early, thus widening the effective window and incurring a higher window‐width penalty. However, that wider window reduces tardiness risk. This trade-off requires balancing two objectives: minimizing the window length and mitigating tardiness risk. The optimal choice of \(\ell^k\) and \(u^k\) are thus determined by weighing the penalty for excessive slack against the risk of late arrivals, as governed by the penalty parameters \(a_w^k\) and \(a_u^k\), receptively.

Using SAA, one can reformulate \eqref{eq:Hfunction_ext} by approximating the expected tardiness via $\frac{1}{Q} \sum_{q=1}^{Q}\left(T^{k[q]} - u^{k}\right)^{+}$, where the realization of $T^k$ is denoted by $T^{k[q]}$ in the $q$-th sample of travel times, $q \in \{1, 2, ..., Q\}$. In order to linearize the reformulation, in addition to the calculation of $T^k$ as presented in \eqref{eq:recursive_2}, one can define the auxiliary variable $v^{k[q]}$ to develop a tractable time window assignment optimization problem for a given routing decision $\boldsymbol{y}$ as follows
\begin{subequations}
\label{eq:SM_ext}
\begin{align} 
    \label{eq:SM_obj_ext} \min_{
    \boldsymbol{\ell}, \boldsymbol{u}, \boldsymbol{v}} & \sum_{k \in V_{0}}\left( a_w^k(u^k - \ell^k) + \frac{a_u^k}{Q} \sum_{q=1}^{Q} v^{k[q]}\right)
    \\
    \mbox{s.t.} \quad
    \label{eq:SM_v2_ext} & T^{k[q]} - u^{k} \leq v^{k[q]} \quad &&\forall k \in V_{0}, \forall q \in\{1,2, \ldots, Q\} \\
    & T^{k[q]} \geq \ell^r + \sum_{a \in \big\{(r, i_{\hat{r}+1}), (i_{\hat{r}+1}, i_{\hat{r}+2}), \dots, (i_{k-1}, k)\big\}} t^{[q]}_a \quad &&\forall k \in V_{0}, \forall q \in\{1,2, \ldots, Q\}, \forall r \in \{i_0, i_1, i_2, \dots, i_{k-1}, k\} \\
    \label{eq:SM_sign_ext} & 0\leq \ell^{k} \leq u^{k}, \quad v^{k[q]} \geq 0 \quad &&\forall k \in V_{0}, \forall q \in\{1,2, \ldots, Q\}
\end{align}
\end{subequations}

While the extended model \eqref{eq:SM_ext} accommodates waiting before the lower bounds, this feature disrupts the original probabilistic linkage between penalty parameters (\(a_w^k, a_u^k\)) and service guarantees. In particular, this modification renders the confidence levels \(1 - a_w^k/a_u^k\) (derived in Propositions \ref{Proposition1}–\ref{Proposition2}) inapplicable to tardiness violations. Restoring those guarantees would require redefining the penalty structure or imposing additional constraints on \(T^k\), fundamentally altering the original framework’s risk-tolerance interpretation. Consequently, we defer a full computational validation of this extension to future work that addresses these theoretical gaps.
\subsection{Design under Partially Known Distribution}
\label{sec:robust} 

An insufficient number of data samples or the unreliability of data samples makes the underlying probability distribution of travel times uncertain. Therefore, the earliness and tardiness measures defined in Section \ref{sec:stochastic} can be affected by the misspecification of the underlying arrival time distribution. Distributionally robust optimization (DRO) is an alternative approach that utilizes limited distributional information. 
The main idea is to embrace the fact that the distribution P is known to belong to an ambiguity set $\mathbb{D}$. This approach has recently become increasingly popular \citep[see][]{delage2010distributionally,wiesemann2014distributionally} and has been applied to routing optimization under uncertainty \citep[see, for instance,][]{carlsson2013robust,mohajerin2018data}.

Any DRO model's tractability and solution performance strongly depends on the limited distributional information and hence the choice of the ambiguity set. To address the correlation between the links' travel times, we assume that the joint distribution P of travel times $\tilde{\boldsymbol{t}}$ belongs to the ambiguity set $\mathbb{D}$ with a given set of information on the mean vector and the covariance matrix. 
Specifically, we assume that the service provider does not have access to the full empirical distribution of travel times (through the full evolving history of travel time observations) but instead relies on the sample estimates of the mean vector, $\boldsymbol{\hat{\mu}}$, and covariance matrix, $\hat{\mathcal{C}}$, which define the ambiguity set $\mathbb{D}$. This way, the DRO model accounts for the uncertainty in these estimates by bounding deviations from the sample mean and restricting the discrepancy between the estimated and true covariance matrix.
To measure the degree of ambiguity about the estimates of mean and covariance, we define $\mathbb{D}$ as

\begin{equation} \label{eq:ambiguity_set1}
\mathbb{D} \triangleq 
\left\{
\begin{array}{l|l}
\mathrm{P} \in \mathcal{M}_{+} & 
\begin{array}{ll}
\mathrm{P}(\tilde{\boldsymbol{t}} \in \mathbb{R}^{|A|})=1 & \text{(a)} \\[1ex]
\left( \mathbb{E}_{\mathrm{P}}(\tilde{\boldsymbol{t}})-\boldsymbol{\hat{\mu}} \right)^\top \boldsymbol{\hat{\mathcal{C}}}^{-1} \left( \mathbb{E}_{\mathrm{P}}(\tilde{\boldsymbol{t}})-\boldsymbol{\hat{\mu}} \right) \leq \alpha_1 & \text{(b)} \\[1ex]
\left\| \operatorname{Cov}_{\mathrm{P}}(\tilde{\boldsymbol{t}})-\boldsymbol{\hat{\mathcal{C}}} \right\|_F \leq \alpha_2,\quad \operatorname{Cov}_{\mathrm{P}}(\tilde{\boldsymbol{t}}) \succeq 0 & \text{(c)}
\end{array}
\end{array}
\right\}
\end{equation}
where $\mathcal{M}_{+}$ is the set of all probability measures on the measurable space $(\mathbb{R}^{|A|}, \mathfrak{B})$ with the $\sigma$-algebra $\mathfrak{B}$ on $\mathbb{R}^{|A|}$. (\ref{eq:ambiguity_set1}b) assumes that the true mean of $\tilde{\boldsymbol{t}}$ lies in an ellipsoid of size $\alpha_1$ centered at the estimate $\boldsymbol{\hat{\mu}}$, and  (\ref{eq:ambiguity_set1}c) forces the Frobenius norm of difference between the estimate $\boldsymbol{\hat{\mathcal{C}}}$ and the true covariance matrix of $\tilde{\boldsymbol{t}}$ to lie in size $\alpha_2$.

We use the ambiguity set $\mathbb{D}$ because it captures our concerns on correlated arc travel times and leads to a tractable optimization model to derive the optimal time windows described below. Using $\mathbb{D}$, we redefine the on-time performance measures used in \eqref{eq:Hfunction} as
\begin{align}
     &\mathcal{H}_{\ell}^k(\boldsymbol{y}^k, \ell^k) \triangleq \sup_{\mathrm{P} \in \mathbb{D}} \mathbb{E}_{\mathrm{P}}\left[\left(\ell^k-\tau^k(\boldsymbol{y}^k,\boldsymbol{\tilde{t}})\right)^+\right] \label{eq:h_l_robust}, \quad \text{and} \\  
     &\mathcal{H}_u^k(\boldsymbol{y}^k, u^k) \triangleq \sup_{\mathrm{P} \in \mathbb{D}} \mathbb{E}_{\mathrm{P}}\left[\left(\tau^k(\boldsymbol{y}^k,\boldsymbol{\tilde{t}})-u^k\right)^+\right] \label{eq:h_u_robust},
 \end{align}
which represent the worst-case earliness and tardiness over the ambiguity set $\mathbb{D}$. Plugging these into the service time window design problem  \eqref{main_subproblem}, we can derive an optimal time window solution for each customer. To cope with the problem's difficulty, we first consider a particular case of $\mathbb{D}$ where we assume the random travel time $\tilde{\boldsymbol{t}}$ with \emph{known} mean $\boldsymbol{\Bar{\mu}}$ and covariance $\boldsymbol{\Bar{\mathcal{C}}} \succeq 0$ follows a family of distributions defined as $\mathcal{F}_{(\boldsymbol{\bar{\mu}}, \boldsymbol{\bar{\mathcal{C}}})} \triangleq \{\mathrm{P} \in \mathcal{M}_+:\mathrm{P}\left(\tilde{\boldsymbol{t}} \in \mathbb{R}^{|A|} \right) = 1, \; \mathbb{E}_\mathrm{P}\left(\tilde{\boldsymbol{t}}\right) = \boldsymbol{\Bar{\mu}}, \; \operatorname{Cov}_\mathrm{P}\left(\tilde{\boldsymbol{t}}\right) = \boldsymbol{\Bar{\mathcal{C}}} \succeq 0\}$. Then, we extend our results for the general case $\mathbb{D}$ in Section \ref{sec:DRO-OM}.

Let us consider the on-time performance measures defined in  \eqref{eq:h_l_robust} and \eqref{eq:h_u_robust}. Given the definition of $\mathcal{F}_{(\boldsymbol{\bar{\mu}}, \boldsymbol{\bar{\mathcal{C}}})}$, the expected value of the random arrival time, $\mathbb{E}_{\mathrm{P}}[\tau^k(\boldsymbol{y}^k,\boldsymbol{\tilde{t}})]$, and the standard deviation of the random arrival time, $\text{STD}_{\mathrm{P}}[\tau^k(\boldsymbol{y}^k,\boldsymbol{\tilde{t}})]$, at customer $k$ can be stated as 
${\boldsymbol{y}^k}^\top \boldsymbol{\bar{\mu}}$ and  $\sqrt{{\boldsymbol{y}^{k}}^{\top} \boldsymbol{\bar{\mathcal{C}}} {\boldsymbol{y}^{k}}}$, respectively. Therefore, we can compute the supremums utilizing the Jensen's inequality  \citep[see][]{scarf1958min} as follows:
\noindent
\begin{align*}
  &   \sup_{\mathrm{P} \in \mathcal{F}_{(\boldsymbol{\bar{\mu}}, \boldsymbol{\bar{\mathcal{C}}})}} \mathbb{E}_{\mathrm{P}}\left[\left(\ell^k-\tau^k(\boldsymbol{y}^k,\boldsymbol{\tilde{t}})\right)^+\right] = 1/2\left(\ell^{k}-{\boldsymbol{y}^{{k}}}^{\top} \boldsymbol{\bar{\mu}} + \sqrt{{{\boldsymbol{y}^{k}} ^ {\top}} \boldsymbol{\bar{\mathcal{C}}} \boldsymbol{y}^{{k}} + \left( \ell^{k}-{\boldsymbol{y}^{{k}}}^{\top} \boldsymbol{\bar{\mu}} \right)^{2}}\right),\\
  & \sup_{\mathrm{P} \in \mathcal{F}_{(\boldsymbol{\bar{\mu}}, \boldsymbol{\bar{\mathcal{C}}})}} \mathbb{E}_{\mathrm{P}}\left[\left(-u^k+\tau^k(\boldsymbol{y}^k,\boldsymbol{\tilde{t}})\right)\right]^+ = 1/2\left(-u^{k}+{\boldsymbol{y}^{{k}}}^{\top} \boldsymbol{\bar{\mu}} + \sqrt{{{\boldsymbol{y}^{k}} ^ {\top}} \boldsymbol{\bar{\mathcal{C}}} \boldsymbol{y}^{{k}} + \left( -u^{k}+{\boldsymbol{y}^{{k}}}^{\top} \boldsymbol{\bar{\mu}} \right)^{2}}\right).
\end{align*}

\noindent Replacing these equations in \eqref{main_subproblem} and ignoring the condition $0\leq \ell^k \leq u^k$ for a moment, we can derive the optimum $\bar{\ell}^k$ and $\bar{u}^k$ using the first order optimality condition as stated in the following proposition.

\begin{proposition}
\label{eq:optimal-time-window}
    For a given route decision $\boldsymbol{y}^{k}$ for customer $k \in V_{0}$, on-time measures \eqref{eq:h_l_robust} and \eqref{eq:h_u_robust} under the ambiguity set $\mathcal{F}_{(\boldsymbol{\bar{\mu}}, \boldsymbol{\bar{\mathcal{C}}})}$ results in the following optimal service time window $[\bar \ell^k, \bar u^k]$:
\end{proposition}
\begin{align}
    \label{eq:lbar_robust}
    &\bar{\ell}^k = \boldsymbol{y}^{k^{\top}} \boldsymbol{\bar{\mu}} - \frac{a_{\ell}^{k}-2 a_{w}^{k}}{\sqrt{1-\left(a_{\ell}^{k}-2 a_{w}^{k}\right)^{2}}} \sqrt{\boldsymbol{y}^{k^{\top}} \boldsymbol{\bar{\mathcal{C}}} \boldsymbol{y}^{k}}, \quad \emph{\text{and}} \\
    \label{eq:ubar_robust}
& \bar{u}^k = \boldsymbol{y}^{k^{\top}} \boldsymbol{\bar{\mu}} + \frac{a_{u}^{k}-2 a_{w}^{k}}{\sqrt{1-\left(a_{u}^{k}-2 a_{w}^{k}\right)^{2}}} \sqrt{\boldsymbol{y}^{k^{\top}} \boldsymbol{\bar{\mathcal{C}}} \boldsymbol{y}^{k}}.
\end{align}

From this proposition, one can observe that the service time window assigned to each customer $k \in V_0$ is built around the expected arrival time at its location, and each wing (the length of such a window on either side) is a positive multiple of the arrival time's standard deviation. It is clear that the optimal time window solution constructed in \eqref{eq:lbar_robust} and \eqref{eq:ubar_robust} satisfies the condition $ \bar{\ell}^k \leq \bar{u}^k$.

Both time window's wings depend on $\boldsymbol{a}^k$, the service provider's risk preference parameters. For a fixed $(a_\ell^k, a_u^k)$, as $a_w^k$ increases, the coefficients of the arrival time's standard deviation in \eqref{eq:lbar_robust} and \eqref{eq:ubar_robust} decrease, which will lead to shorter wings around the expected arrival time. This is in line with the definition of $a_w^k$ as the penalty of time window's length. For a fixed value of $a_w^k$, choosing different values to parameters $(a_\ell^k, a_u^k)$ will impact the time window differently. If $a_\ell^k = a_u^k$, the service provider has the same time window violation tolerance on either side of the window. Hence, we acquire a symmetric time window centered on the expected arrival time at the customer's location.
However, if  $a_\ell^k < a_u^k$, the service provider is more concerned about the tardiness than the earliness. Therefore, the arrival time's standard deviation on the right-hand side of the expected arrival will be multiplied by a larger coefficient, resulting in a longer right wing assigned to the customer. A similar argument can be developed for the case where $a_\ell^k > a_u^k$.
\section{Integrated Routing and Service Time Window Design}
\label{sec:integrated} 

The procedure described in Section \ref{sec:model} takes as input a routing decision (i.e., sequence of customers to be visited in last-mile delivery) and provides a time window to visit each customer with a certain level of confidence based on the service provider's risk tolerance. Even though the confidence levels are only functions of the service provider's risk tolerance parameters, the lengths of the resulting time windows and the corresponding percentage and amount of time window violations could differ for any two input routes, as demonstrated in Figure \ref{fig:two_routes}. 
Therefore, in this section, we develop a modeling framework that simultaneously optimizes the routing decision and the time windows design in last-mile delivery.

Let us show by $S_{xy}$ the set of feasible routes each of which is a Hamiltonian path that starts from the depot, visits each customer $k\in V_0$ exactly once, and ends again at the depot within a time budget TB. The time budget can be interpreted as the available driver's shift or the maximum duration the service provider is willing to allot to serving all the customers. Set $S_{xy}$ contains two main sets of binary decision variables: (i) $x_{ij}$ for each arc $(i,j)\in A$, which is equal to 1 if arc $(i,j)$ is in the route, and 0 otherwise; and (ii) $y_{ij}^k$ which becomes equal to 1 when going from deport to customer $k$ requires traversing arc $(i,j)\in A$. The mathematical description of $S_{xy}$ is give as:

\begin{equation}\label{eq:routes}
S_{xy}=\left\{\begin{array}{l|l}
\begin{array}{l}\boldsymbol{x} \in\{0,1\}^{|A|}\\
\boldsymbol{y} \in \mathbb{R}_{+}^{|A| \times\left|V_{0}\right|}\end{array} & \begin{array}{lr}
\sum_{(i, j) \in \delta^{+}(i)} x_{i j}=1 & \forall i \in V \quad (\mathrm{a}) \\
\sum_{(j, i) \in \delta^{-}(i)} x_{j i}=1 & \forall i \in V \quad (\mathrm{b}) \\
\sum_{(i, j) \in \delta^{+}(i)} y_{i j}^{k}-\sum_{(j, i) \in \delta^{-}(i)} y_{j i}^{k}=\left\{\begin{array}{ll}
1, & \text { if } i=0 \\
-1, & \text { if } i=k \\
0, & \text { otherwise }
\end{array}\right. & \forall k \in V_{0}, \forall i \in V \quad (\mathrm{c}) \\
y_{i j}^{k} \leq x_{i j} & \forall k \in V_{0}, \forall(i, j) \in A \quad (\mathrm{d})
\end{array}
\end{array}\right\}    
\end{equation}
Constraints (\ref{eq:routes}a) and (\ref{eq:routes}b) state that each node must be visited exactly once. Constraints (\ref{eq:routes}c) represent the flow conservation constraints, which ensure that the route must start from and end at origin 0. Constraints (\ref{eq:routes}d) guarantee that arc $(i,j) \in A$ can be used for reaching customer $k \in V_0$ only when it exists in the route. Note that the binary requirement of variables $\boldsymbol{y}$ is guaranteed by (\ref{eq:routes}c) and (\ref{eq:routes}d).

\begin{remark}
We can extend the route description to the multiple capacitated vehicles by incorporating the index $h\in H$ (with $H$ as the set of vehicles) into the variables $\boldsymbol{x}$ and $\boldsymbol{y}$. In particular, for each vehicle $h\in H$ and each arc $(i,j)\in A$, we can modify  the binary variable $x_{ij}$ to $x_{ij}^h$ which is equal to 1 if arc $(i,j)$ is traversed by vehicle $h$, and $0$ otherwise. We can then impose vehicle capacity on each route and ensure that all routes are connected to the depot by adding the following constraints to $S$
\begin{equation*}
    \textstyle\sum_{h\in H}\textstyle\sum_{i\notin \mathcal{V
    }}\textstyle\sum_{j \in \mathcal{V
    }:(i,j)\in A} x^h_{ij} \geq \gamma (\mathcal{V
    }) \quad \forall \mathcal{V} \subseteq V_0,
\end{equation*}
where $\gamma(\mathcal{V
    })$ shows the minimum number of vehicles required to serve the customers in subset $\mathcal{V
    }\subseteq V_0$ according to their demands.
\end{remark}

The goal of the integrated routing and service time window design  is to plan an optimal route (i.e., optimal $\boldsymbol{x}$ and $\boldsymbol{y}$) with tight service time windows (i.e., optimal $\boldsymbol{\ell}$ and $\boldsymbol{u}$) that minimize the overall service time window design cost. This can be achieved by solving the following optimization model (OM):
\begin{subequations}
    \label{eq:om}
    \begin{align} 
    \label{eq:om_obj}
    \text { \textbf{OM}: } \min_{\boldsymbol{x}, \boldsymbol{y}, \boldsymbol{\ell},\boldsymbol{u}} \quad & \sum_{k \in V_{0}} \Big(a_w^k(u^k-\ell^k) + a_{\ell}^k\mathcal{H}_{\ell}^k(\boldsymbol{y}^k, \ell^k) + a_u^k\mathcal{H}_u^k(\boldsymbol{y}^k, u^k) \Big)\\
    \mbox{s.t.}\quad&
    \label{eq:om_order}  0 \leq {\ell}^k \leq  {u}^k \quad \forall k\in V_0\\
    &\mathcal{H}_{TB}(\boldsymbol{x}) \leq \mathrm{TB} 
    \label{eq:om_budget} \\
    & (\boldsymbol{x}, \boldsymbol{y}) \in S_{xy}.  \nonumber
    \end{align}
\end{subequations}

The function $\mathcal{H}_{TB}(\boldsymbol{x})$ measures the total time needed to complete the route and constraint \eqref{eq:om_budget} ensures that the time budget is not violated. The OM model can be represented as stochastic programming or a DRO under the fully and partially known distribution of random travel times studied in Sections \ref{sec:stochastic} and \ref{sec:robust}, respectively. For the former case,  
$\mathcal{H}_{TB}(\boldsymbol{x}) = \mathbb{E}_{\mathrm{P}}\left(\tilde{\boldsymbol{t}}^{\top}\boldsymbol{x}\right)$ representing the expected completion time of the tour, while in the latter case, $\mathcal{H}_{TB}(\boldsymbol{x}) = \sup _{\mathrm{P} \in \mathbb{D}}\mathbb{E}_{\mathrm{P}}\left(\tilde{\boldsymbol{t}}^{\top}\boldsymbol{x}\right)$ indicating the tour's worst expected completion time in $\mathbb{D}$. In the following subsections, we show how to derive tractable deterministic models for the stochastic and DRO models. 


\subsection{OM under the Fully Known Distribution}
\label{sec:SAA-OM}

Under the SAA described in Section \ref{sec:stochastic}, the OM model can be expressed as the following sample-based optimization model SM: 
\begin{subequations}
\label{eq:SM}
\begin{align} 
    \label{eq:SM_obj}\text { \textbf{SM}: } \min_{\boldsymbol{x}, \boldsymbol{y}, \boldsymbol{\ell}, \boldsymbol{u}, \boldsymbol{v_{1}}, \boldsymbol{v_{2}}} & \sum_{k \in V_{0}}\left( a_w^k(u^k - \ell^k) + \frac{a_\ell^k}{Q} \sum_{q=1}^{Q} v_{1}^{k[q]} + \frac{a_u^k}{Q} \sum_{q=1}^{Q} v_{2}^{k[q]}\right)
    \\
    \mbox{s.t.} \quad
    \label{eq:SM_TB}&\frac{1}{Q} \sum_{q=1}^{Q} \sum_{(i, j) \in A} t_{i j}^{[q]} x_{i j} \leq \mathrm{TB} \\
    \label{eq:SM_v1}&\ell^{k}-\sum_{(i, j) \in A}\left(t_{i j}^{[q]} y_{i j}^{k}\right) \leq v_{1}^{k[q]} \qquad \forall k \in V_{0}, \forall q \in\{1,2, \ldots, Q\} \\
    \label{eq:SM_v2}&\sum_{(i, j) \in A}\left(t_{i j}^{[q]} y_{i j}^{k}\right)-u^{k} \leq v_{2}^{k[q]} \qquad \forall k \in V_{0}, \forall q \in\{1,2, \ldots, Q\} \\
    &(\boldsymbol{x}, \boldsymbol{y}) \in S_{xy} \\
    \label{eq:SM_sign}&0\leq \ell^{k}\leq u^{k}, \quad v_{1}^{k[q]}, v_{2}^{k[q]} \geq 0 \qquad \forall k \in V_{0}, \forall q \in\{1,2, \ldots, Q\}
\end{align}
\end{subequations}

The SM is a mixed integer linear program (MILP) that can be solved efficiently by  state-of-the-art optimization solvers. However, when the underlying network is dense and/or the number of samples is large, it becomes very challenging for the solvers to obtain the optimal solution or even a feasible solution in a reasonable time or memory usage. 

One way to deal with this difficulty is to partition the SM into an integer master problem (to find the best route) and linear subproblems (to obtain time windows for customers) that are more manageable in size and computationally easier to solve with respect to the original model OM. The routing decisions $(\boldsymbol{x}, \boldsymbol{y})$ are incorporated into the master problem, while variables $(\boldsymbol{\ell},\boldsymbol{u}, \boldsymbol{v_1},\boldsymbol{v}_2)$ associated with time windows and linearization are projected out and replaced by a  variable $\omega_k$. The resulting master problem, which we refer to as MP(SM), is then given by
\begin{subequations}
    \label{eq:Benders}
    \begin{align} 
    \label{eq:Benders_obj}\text { \textbf{MP(SM):}} ~\min_{\boldsymbol{x}, \boldsymbol{y},\omega} \quad &\sum_k\omega^{k}\\
    \mbox{s.t.} 
    \label{eq:Benders_omegak} \quad &\omega^{k} \geq \phi^{k}(\boldsymbol{y}^k) \quad \forall k \in V_{0}\\
    \nonumber&(\boldsymbol{x}, \boldsymbol{y}) \in S_{xy} \\
    &\omega^{k} \geq 0,\quad \forall k \in V_{0}
    \end{align}
\end{subequations}
where the convex (not necessarily differentiable everywhere) function $\phi^{k}(\boldsymbol{y}^k)$ appearing in \eqref{eq:Benders_omegak} gives the cost associated with the time window assignment for each customer $k \in V_0$ as defined in \eqref{eq:SP(SM)}. 

The decomposition idea is based on successively adding cuts in the $(\boldsymbol{x}, \boldsymbol{y}, \omega)$-space to approximate $\phi^k$ until an optimal solution $(\boldsymbol{x}^*, \boldsymbol{y}^*, \omega^*)$ with $\omega^*=\sum_{k}\phi^{k}(\boldsymbol{y}^k)$ is identified. Because of convexity, function $\phi^{k}(\boldsymbol{y})$ can be underestimated by a supporting hyperplane at $\hat{\boldsymbol{y}}$, so we can write the following linear inequality, known as a generalized Benders cut \citep[see][]{geoffrion1972generalized}:
\begin{align}
\label{eq:GBCut}
    \omega^{k} \geq\ \phi^{k}(\boldsymbol{y}) \geq \phi^{k}(\hat{\boldsymbol{y}})+\sum_{(i, j) \in A} \hat{s}_{i j}^{k}\left(y_{i j}^{k}-\hat{y}_{i j}^{k}\right) \quad \forall k \in V_{0},
\end{align}
where $\hat{s}_{i j}^{k} \in \partial \phi^{k}(\hat{\boldsymbol{y}})$ is any \emph{subgradient} of $\phi^{k}$ at $\hat{\boldsymbol{y}}$. 
The following proposition formally shows how to derive these subgradients without the need of solving any linear program. 

\begin{proposition}
\label{prop:subgradient}
    Given $(\hat{\boldsymbol{x}},\hat{\boldsymbol{y}})$, a subgradient $\hat{s}_{i j}^{k} \in \partial \phi^{k}(\hat{\boldsymbol{y}})$ for each $k\in V_0$ and $(i,j)\in A$ can be obtained as
\begin{equation*}
    \hat{s}_{i j}^{k} =\sum_{q=1}^{Q} t_{i j}^{[q]}\left( {\bar \rho_{2}^{k[q]}}- {\bar \rho_{1}^{k[q]}} \right),
\end{equation*}
where, ${\bar \rho_{1}^{k[q]}}$ and ${\bar \rho_{2}^{k[q]}}$ are the optimal values of the dual variables associated with constraints \eqref{eq:SM_v1} and \eqref{eq:SM_v2}, respectively that are computed as
\begin{align*}
  & {\bar \rho}_1^{k[q]} = \begin{cases}
 \frac{a_\ell^k}{Q}  & q=1,\ldots, P_1^{k}-1\\
  a_w^k - \sum_{q=1}^{P_1^{k}-1} \frac{a_\ell^k}{Q}  & q=P_1^{k}\\
0& q=P_1^{k}+1, \ldots, Q
 \end{cases}
\quad &   {\bar \rho}_2^{k[q]} = \begin{cases}
 0& q=1, \ldots, P_2^{k}-1\\
  a_w^k - \sum_{q=P_2^{k}+1}^{Q} \frac{a_u^k}{Q}  & q=P_2^{k}\\
  \frac{a_u^k}{Q}  & q=P_2^{k}+1, \ldots, Q.
 \end{cases}
\end{align*}

\end{proposition}

\proof{Proof.} See Appendix C.
\hfill $\square$

In our implementation, which is evaluated in Section \ref{sec:comput}, we solve MP(SM) using a branch-and-cut framework of a state-of-the-art optimization solver. The optimality cuts are incorporated into the master problem by using callbacks allowing to add the cutting planes \eqref{eq:GBCut} step-by-step. A callback is executed whenever an optimal solution of the LP-relaxation is found at the root node of the branch-and-bound tree or an incumbent solution at any node of the branch-and-bound tree is found. For the current choice of variables $(\boldsymbol{x}, \boldsymbol{y})$, the subgradients are computed and the resulting cuts \eqref{eq:GBCut} are added to the master problem if they are violated. This procedure continues until an incumbent solution is found where none of the corresponding cuts are violated.

\subsubsection{Extension 1: Fixed-Length Time Windows.}
\label{subsec:stochastic_fixed_sol}

The modeling approach and Benders decomposition method proposed here can be similarly applied to the SAA-based model \eqref{eq:SP(SM)_appx_2} developed for fixed-length time windows in Section \ref{subsec:stochastic_fixed}. That model can be extended to include the routing decision variables $(\boldsymbol{x}, \boldsymbol{y}) \in S_{xy}$ and the time budget constraint \eqref{eq:SM_TB} to simultaneously optimize the routing plans and fixed-length time window assignments with certain service guarantees. This enables us to evaluate how fixed-length constraints affect time window characteristics and solution costs compared to variable windows within a stochastic programming framework—an analysis we present in Section \ref{sec:comput}.

\subsection{OM under the Partially Known Distribution}
\label{sec:DRO-OM}

Here, using the time window characteristics described in Section \ref{sec:robust}, we present how to convert the OM under the DRO setting to a deterministic optimization model through the following proposition. This is accomplished by turning from a particular case of the ambiguity set $\mathbb{D}$ in Proposition \ref{eq:optimal-time-window} where we assumed the random travel time $\tilde{\boldsymbol{t}}$ has \emph{known} mean $\boldsymbol{\Bar{\mu}}$ and covariance $\boldsymbol{\Bar{\mathcal{C}}}$ to the general case of $\mathbb{D}$ defined in \eqref{eq:ambiguity_set1} with \emph{observed} $\boldsymbol{\hat{\mu}}$ and $\boldsymbol{\hat{\mathcal{C}}}$. We define  $\mathcal{U}_{(\boldsymbol{\hat{\mu}}, \boldsymbol{\hat{\mathcal{C}}})} = \mathcal{U}_{\boldsymbol{\hat{\mu}}}\times\mathcal{U}_{\boldsymbol{\hat{\mathcal{C}}}}$ to show the support set of all mean vectors $\boldsymbol{\mu}$ and covariance matrices $\boldsymbol{\mathcal{C}} > 0$ satisfying  (\ref{eq:ambiguity_set1}b) and  (\ref{eq:ambiguity_set1}c) with
\begin{align}
    &\mathcal{U}_{\boldsymbol{\hat{\mu}}} \triangleq \left\{ \boldsymbol{\mu} : \, \left( \boldsymbol{\mu}-\boldsymbol{\hat{\mu}}\right)^{\top} \boldsymbol{\hat{\mathcal{C}}}^{-1}\left(\boldsymbol{\mu}-\boldsymbol{\hat{\mu}}\right) \leq \alpha_1 \right\}, \label{eq:U_mu}\\
    &\mathcal{U}_{\boldsymbol{\hat{\mathcal{C}}}} \triangleq \left\{\boldsymbol{\mathcal{C}}  :\,\|\boldsymbol{\mathcal{C}}-\boldsymbol{\hat{\mathcal{C}}} \|_F \leq \alpha_2 \right\}. \label{eq:U_C}
\end{align}
This way, we factor in the size of the ambiguity set determined by the positive parameters $\alpha_1$ and $\alpha_2$ which provide means of quantifying the service provider's confidence in the observed values of the mean vector and covariance matrix, respectively.


\begin{proposition}
\label{theorem:Robust2-2}
The DRO reformulation of the OM under the ambiguity set $\mathbb{D}$ is equivalent to the following deterministic optimization model: 
\begin{subequations}
\label{eq:RM2}
    \begin{align}
    \text { \textbf{RM:}} ~    
    \label{eq:RM2_obj}
    \min_{\boldsymbol{x}, \boldsymbol{y}} \quad & \sum_{k \in V_{0}} (\Gamma_\ell^k + \Gamma_u^k) \sqrt{{\boldsymbol{y}^{k}}^{\top} (\boldsymbol{\hat{\mathcal{C}}} + \alpha_2 I_{|A|}) \boldsymbol{y}^{k}} \\
    \label{eq:RM2_budget} \mbox{s.t.} \quad
    &\boldsymbol{\hat{\mu}}^{\top} \boldsymbol{x} + \sqrt{\alpha_1 \,(\boldsymbol{x}^{\top} \boldsymbol{\hat{\mathcal{C}}} \boldsymbol{x})} \leq \mathrm{TB} \\
    & (\boldsymbol{x}, \boldsymbol{y}) \in S_{xy}, \nonumber
    \end{align}
\end{subequations}
where $I_{|A|}$ is the identity matrix of size ${|A|}$, and 
\begin{equation*}
    \Gamma_\ell^k = \frac{a_{\ell}^{k} - (a_{\ell}^{k}-2 a_{w}^{k})^2}{2\sqrt{1-\left(a_{\ell}^{k}-2 a_{w}^{k}\right)^{2}}} \quad \text{and} \quad
    \Gamma_u^k = \frac{a_{u}^{k} - (a_{u}^{k}-2 a_{w}^{k})^2}{2\sqrt{1-\left(a_{u}^{k}-2 a_{w}^{k}\right)^{2}}}.
\end{equation*}
\end{proposition}

\proof{Proof.} See Appendix D.
\hfill $\square$

Two immediate corollaries can be observed under the ambiguity set $\mathbb{D}$. First, the objective is to minimize the overall standard deviation (for all customers) of random arrival time at each customer
$\text{STD}_{\mathrm{P} \in \mathbb{D}}[\tau^k(\boldsymbol{y}^k,\boldsymbol{\tilde{t}})] = \sqrt{{\boldsymbol{y}^{k}}^{\top} (\boldsymbol{\hat{\mathcal{C}}} + \alpha_2 I_{|A|}) \boldsymbol{y}^{k}}$, which is the worst case standard deviation over $\mathbb{D}$, i.e., $\sup_{\boldsymbol{\bar{\mathcal{C}}} \in \mathcal{U}_{\boldsymbol{\hat{\mathcal{C}}}}} \sqrt{{\boldsymbol{y}^{k}}^{\top} \boldsymbol{\Bar{\mathcal{C}}} \boldsymbol{y}^{k}}$. Second, the expected value of the random arrival at each customer, $\mathbb{E}_{\mathrm{P} \in \mathbb{D}}[\tau^k(\boldsymbol{y}^k,\boldsymbol{\tilde{t}})]$, is $\boldsymbol{\hat{\mu}}^{\top} \boldsymbol{y}^k + \sqrt{\alpha_1} \sqrt{{\boldsymbol{y}^k}^{\top} \boldsymbol{\hat{\mathcal{C}}} \boldsymbol{y}^k}$, which is the worst case expected arrival over $\mathbb{D}$, i.e.,
$\sup_{\boldsymbol{\bar{\mu}} \in \mathcal{U}_{\boldsymbol{\hat{\mu}}}} \boldsymbol{\bar{\mu}}^{\top} \boldsymbol{y}^k$.

The RM in \eqref{eq:RM2} is a mixed integer non-linear program (MINLP) that can be reformulated as a mixed integer conic quadratic program (MICQP) by introducing a positive continuous variable $\vartheta^k$ for each customer $k\in V_0$ as follows:
\begin{subequations}
    \label{eq:RM-S'}
    \begin{align} 
    \label{eq:RM-S'_obj}\text { \textbf{RM'}: } \min_{\boldsymbol{x}, \boldsymbol{y},\boldsymbol{\vartheta}} &\quad  \sum_{k \in V_{0}} (\Gamma_{\ell}^k + \Gamma_{u}^k) \vartheta^k \\ \nonumber
     \mbox{s.t.} \quad &(\boldsymbol{x}, \boldsymbol{y}) \in S_{xy}, ~
    \eqref{eq:RM2_budget} \\
    \label{eq:RM-S'_conic}  &  \sum_{(i, j) \in A} \sum_{(r, s) \in A} \bar{\bar{\mathcal{C}}}_{i j r s} y_{i j}^{k} y_{r s}^{k} \leq {\left(\vartheta^k\right)}^2 \qquad \forall k \in V_0 \\
    \label{eq:RM-S'_set} & \vartheta^k \in \mathbb{R}_+ \qquad \forall k \in V_0,
    \end{align}
\end{subequations}
where $\bar{\bar{\mathcal{C}}}_{i j r s}$ is an entry of matrix $\boldsymbol{\hat{\mathcal{C}}} + \alpha_2 I_{|A|}$ in \eqref{eq:RM2_obj} representing the robust covariance between the two arcs $(i, j) \in A$ and $(r, s) \in A$.
If the binary restrictions on variables $\boldsymbol{x}$ are relaxed, the above formulation will become a second order cone program (SOCP), also known as a conic quadratic program. Due to their special structure, SOCP are computationally tractable and can be solved by interior-point algorithms in polynomial time. Therefore, the transformation of model \eqref{eq:RM2} to \eqref{eq:RM-S'} facilitates the solution through the embedded branch-and-cut algorithm in state-of-the-art solvers such as CPLEX and Gurobi. One can find the overview of the SOCP in \cite{ben2001polyhedral} and \cite{alizadeh2003second}.

However, the complexity of RM' increases as the  number of customers increases and/or the underlying network is dense. This is because  of decision variables $\boldsymbol{y}^k$, with index $k$ for each customer, and constraints \eqref{eq:RM-S'_conic}. To overcome this challenge, we develop a decomposition technique based on outer approximation (OA) described as follows. We consider the main model RM and introduce the convex function $\phi^{k}(\boldsymbol{y})$ defined as 
\begin{align*}
    \phi^{k}(\boldsymbol{y}) \triangleq \sqrt{\sum_{(i, j) \in A} \sum_{(r, s) \in A} \mathcal{\bar{\bar{C}}}_{i j r s} y_{i j}^{k} y_{r s}^{k}} \qquad \forall k \in V_{0},
\end{align*}
to result the following reformulation of RM:
\begin{subequations}
\label{eq:MP(RM-S)}
    \begin{align} 
    \label{eq:MP2_obj}\text {\textbf{MP(RM):}  } \min_{\boldsymbol{x}, \boldsymbol{y}, \boldsymbol{\omega}} & \quad  \sum_{k \in V_0} (\Gamma_{\ell}^k + \Gamma_{u}^k) \omega^k \\
    \nonumber \mbox{s.t.}  \quad &(\boldsymbol{x}, \boldsymbol{y}) \in S_{xy}, ~
    \eqref{eq:RM2_budget} \\
    \label{eq:MP2_phi}&\omega^k\geq \phi^k(\boldsymbol{y}) \qquad \forall k \in V_0 \\
    &\omega^k \geq 0 \qquad \forall k \in V_0.
\end{align}
\end{subequations}

Because of convexity, function $\phi^k (\boldsymbol{y})$ can be underestimated by a supporting hyperplane at any feasible solution $\hat{\boldsymbol{y}}$ according to \eqref{eq:GBCut} with $\hat{s}_{i j}^{k} \in \partial \phi^{k}(\hat{\boldsymbol{y}})$ being a subgradient of $\phi^{k}$ at $\hat{\boldsymbol{y}}$ computed as 

$$
\hat{s}_{i j}^{k}=\frac{\partial \phi^{k}(\hat{\boldsymbol{y}})}{\partial y_{i j}^{k}}=\frac{\sigma_{i j}^{2} \hat{y}_{i j}^{k}+\sum_{\acute{a} =(r, s) \in A \atop \acute{a} \neq a} \mathcal{\bar{\bar{C}}}_{i j r s} \hat{y}_{r s}^{k}}{\sqrt{\sum_{(i, j) \in A} \sum_{(r, s) \in A} \mathcal{\bar{\bar{C}}}_{i j r s} \hat{y}_{i j}^{k} \hat{y}_{r s}^{k}}} \qquad \forall k \in V_{0}, \forall a=(i, j) \in A.
$$

The overall branch-and-cut algorithm we have implemented to solve the  MP(RM) is similar to what we explained for MP(SM) in Section \ref{sec:SAA-OM}.  More precisely, for the current choice of variables $(\boldsymbol{x}, \boldsymbol{y})$, the subgradients are computed and the resulting cuts \eqref{eq:GBCut} are added to the master problem if they are violated. We have also implemented the single cut strategy by aggregating $\omega^k$ variables as a single variable $\omega =\sum_{k\in V_0}\omega_k$.  Note that one can follow the same idea to linearize constraints \eqref{eq:RM2_budget} using a supporting hyperplane. However, in our numerical experiments, we found that this has only a marginal impact on the computational efficiency, and the main bottleneck is \eqref{eq:MP2_phi}.
\section{Numerical Experiments} 
\label{sec:comput}

This section presents the computational study evaluating  our proposed models and the solution algorithms as well as discussing managerial implications to the reliable delivery operations management. Aligned with our formulations, all numerical experiments study the case of a single vehicle. We aim to address two main questions: First, whether the newly proposed on-time metrics and models can provide reasonable and reliable time window solutions as well as valuable managerial insights for the service provider under both  full and partial statistical information. Second, whether the proposed solution methodologies are capable of reducing the computational burden of solving the mathematical models. To address the first question, we solved the SM and RM' models in \eqref{eq:SM} and \eqref{eq:RM-S'}, respectively, by CPLEX based on the problem instances in \cite{adulyasak2016models}. Because these instances are defined on a sparse (incomplete) graph, they are not too challenging for the solver. The largest instances (50 customers) of the SM and RM' models were solved in an average of 40 minutes and 8 minutes, respectively. Therefore, we tested our decomposition algorithms on several  dense problem instances presented in  \cite{rostami2021branch} to address the second question.
 
All the models and decomposition algorithms were coded in Python, and all the instances were  run on a PC with an Intel Core i9 CPU processor @ 1.90GHz, 10 Cores, and 32GB RAM by calling CPLEX 22.1 as MILP and MINLP solver. CPLEX was set to exploit parallel computations (using 20 threads) while it solved the nodes of the branch-and-cut tree for all the models and algorithms. The generic callbacks were performed in CPLEX for the decomposition algorithms to separate integer feasible LP solutions in a context of lazy constraints.

\subsection{Datasets}
\label{subsec:problemInstances}

We consider six datasets introduced and used in \cite{jaillet2016routing} and \cite{adulyasak2016models} to evaluate the designed time windows for customers. These datasets are called IG-1 to IG-6, each of which consists of 20 problem instances. The IG-1, IG-3, IG-4, IG-5, and IG-6 are composed of the instances of  size $|V_0| =$ 10, 20, 30, 40, and 50, respectively, and $|A| = 3|V_0|$. IG-2 is the same as IG-1 with $|A|=50$. Since these datasets only provide $\hat{\boldsymbol{\mu}}$ and the time budget parameters, we needed to generate a positive semidefinite covariance matrix $\hat {\boldsymbol{\mathcal{C}}}$ for each. The details are provided in Appendix E. 
However, the datasets in \cite{rostami2021branch} contain the covariance matrices, which are used to assess the proposed decomposition algorithms' performance in a complete network.

Recall in the SM model, the assumption is the precise knowledge of the distribution of the random travel time vector $\Tilde{\boldsymbol{t}}$, i.e., P. In contrast, in the RM model, we adopt an ambiguous distribution of the travel times, where the travel time mean vector $\hat {\boldsymbol{\mu}}$ and covariance matrix $\hat {\boldsymbol{\mathcal{C}}}$ are known. Since the focus of our analyses is not on the size of the ambiguity set, we assume that $\alpha_1 = \alpha_2 = 0$. However, our results can be replicated for any other values of $\alpha_1$ and $\alpha_2$. To be able to compare the results of the analysis in the DRO setting with the case where P is known, we generated $Q$ sample travel times vectors $\boldsymbol{t}^{[1]}, \boldsymbol{t}^{[2]}, ...,\boldsymbol{t}^{[Q]}$ with known $\hat{\boldsymbol{\mu}}$ and $\hat{\boldsymbol{\mathcal{C}}}$.

\subsection{Models' Evaluation and Managerial Insights} 

In this section, we evaluate the capability of the proposed models in helping the service provider with providing reliable service time windows for its customers. 
For conducting analysis with the SAA method in SM, we generated $Q=1000$ sample travel times $t_{ij}^{[q]}, \forall (i,j) \in A$ and $\forall q \in \{1,2,..., Q\}$, from a Normal distribution with the known $\hat{\boldsymbol{\mu}}$ and $\hat{\boldsymbol{\mathcal{C}}}$.
To evaluate the performance of the routes and the time windows created by solving the above models, we also generated 1000 out-of-sample test instances from the same Normal distribution. That is, the previously designed routes and time windows are now tested on the graph with the new arcs' travel times to investigate the violation of the time windows on both sides (before and after the time windows). In our experiments, we examined our models for different choices of penalty parameters,  resulting in the different confidence levels with the violation rates $\beta_\ell, \beta_u \in \{0.025, 0.05, 0.075\}$. For example, a $\beta_\ell = 0.075, \beta_u = 0.025$ indicates that the service provider expects to arrive at the customers' locations before the earliest times and after the latest times assigned to them only at most $7.5\%$ and $2.5\%$ of times, respectively.  In other words, the service provider wants to be at least $90\%$ confident that all the services will start within the assigned time windows without early or late violations. Although, a higher level of service guarantee is desired for the late arrivals as more penalty is assigned to them.

\subsubsection{Models' Performance in Providing Reliable Delivery.} 
\begin{figure}[h!]
    \centering
    \includegraphics[scale=0.9]{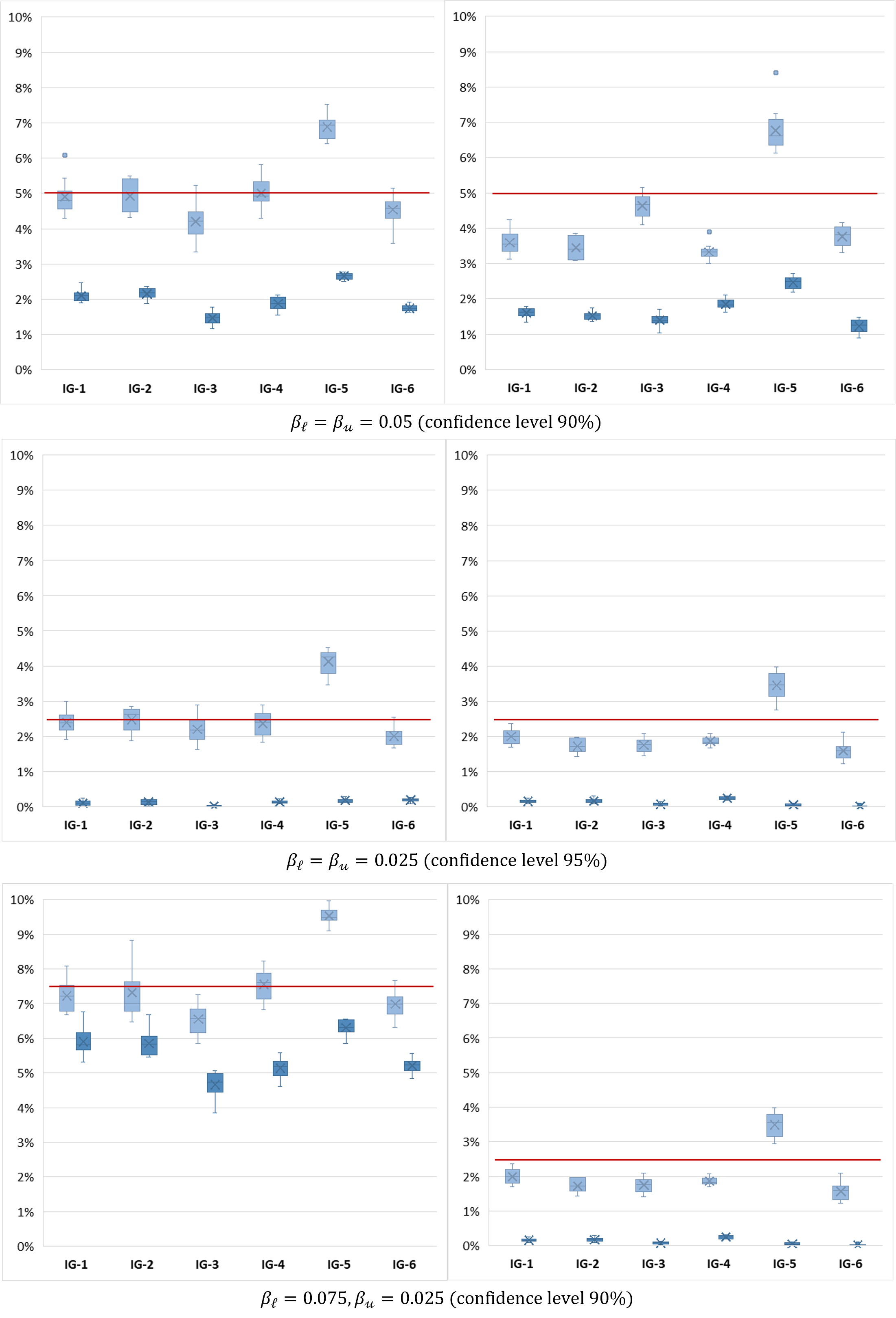}
    \caption{Percentage of out-of-sample time window violations before the earliest time (left) and after the latest time (right) in SM (lighter bars) and RM (darker bars).}
    \label{fig:2}
\end{figure}

Figure \ref{fig:2} presents the out-of-sample performance of both the SM and RM under different combinations of $\beta_\ell$ and $\beta_u$. Each diagram illustrates the percentage of arrivals at the customers in the test instances before or after the assigned time windows. The red line in each figure depicts the maximum acceptable violation rate specified by the service provider's desired confidence level on each side of the window. As can be seen, routes and time windows designed under the robust model consistently keep violations below the red line, whereas the stochastic model occasionally exceeds the acceptable threshold. When the confidence level increases from $90\%$ to $95\%$, the violation rates under both models decrease and are more significant under the robust setting. Moreover, with the same confidence level of $90\%$ in the first and the third diagrams, redistributing  $\beta_\ell$ and $\beta_u$ would result in different early arrival and late arrival violations which are aligned with our theoretical observations in Sections \ref{sec:stochastic} and \ref{sec:robust}. 

\subsubsection{Length of Delivery Time Windows.} 
\begin{figure}[t!]
    \centering
    \includegraphics[scale=0.9]{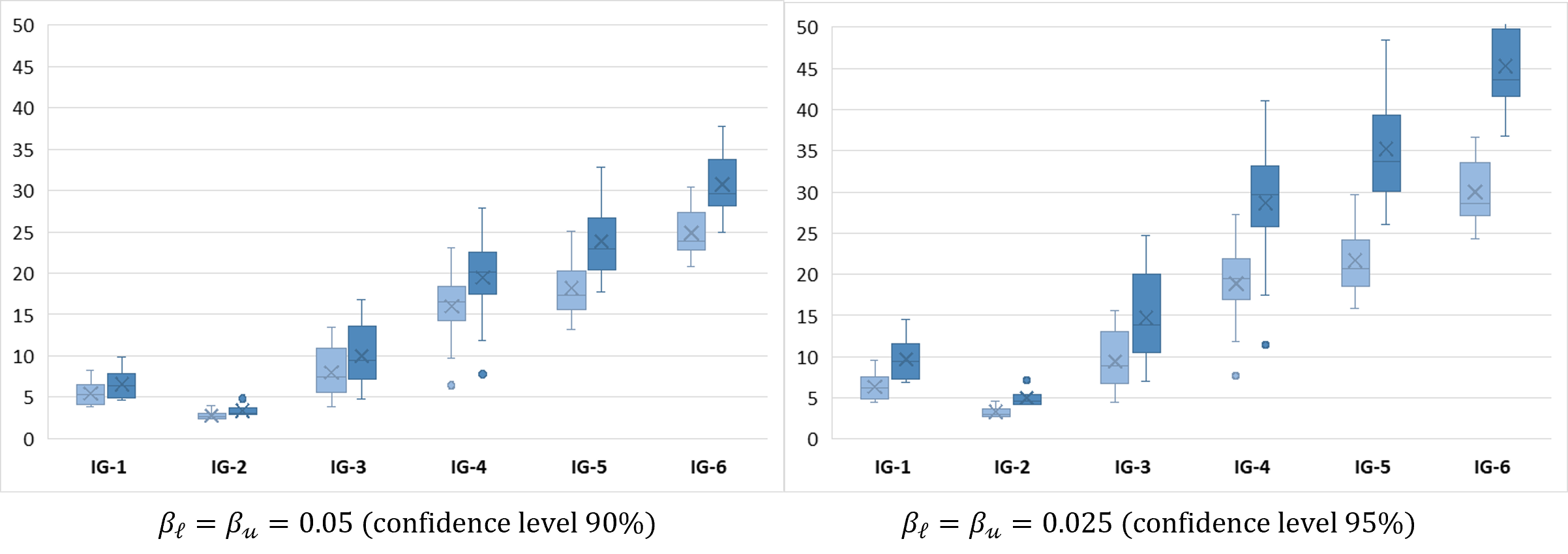}
    \caption{Time window lengths in SM (lighter bars) and RM (darker bars)}
    \label{fig:3}
\end{figure}

Figure \ref{fig:3} illustrates the lengths of the time windows generated under both the SM and RM across various instance groups and confidence levels, providing a direct comparison of how each method sizes its time windows. As it can be observed, the service time windows provided by the SM are tighter than those given by the RM for both confidence levels. Therefore, having a less time window violation rate under the RM (Figure \ref{fig:2}) comes at the expense of assigning longer time windows to the customers compared to the time windows designed under the SM. This highlights the fact that a trade-off exists between the time window length and the number of its violation. Overall, the results presented in  Figures \ref{fig:2} and \ref{fig:3} demonstrate how the RM model is able to design more reliable routes and time windows that are violated much less than the rate allowed by the confidence levels through making the time windows somewhat longer.

\subsubsection{Impact of Delivery's Guarantee Level.} 
\begin{figure}[htp]
    \centering
    \includegraphics[scale=0.8]{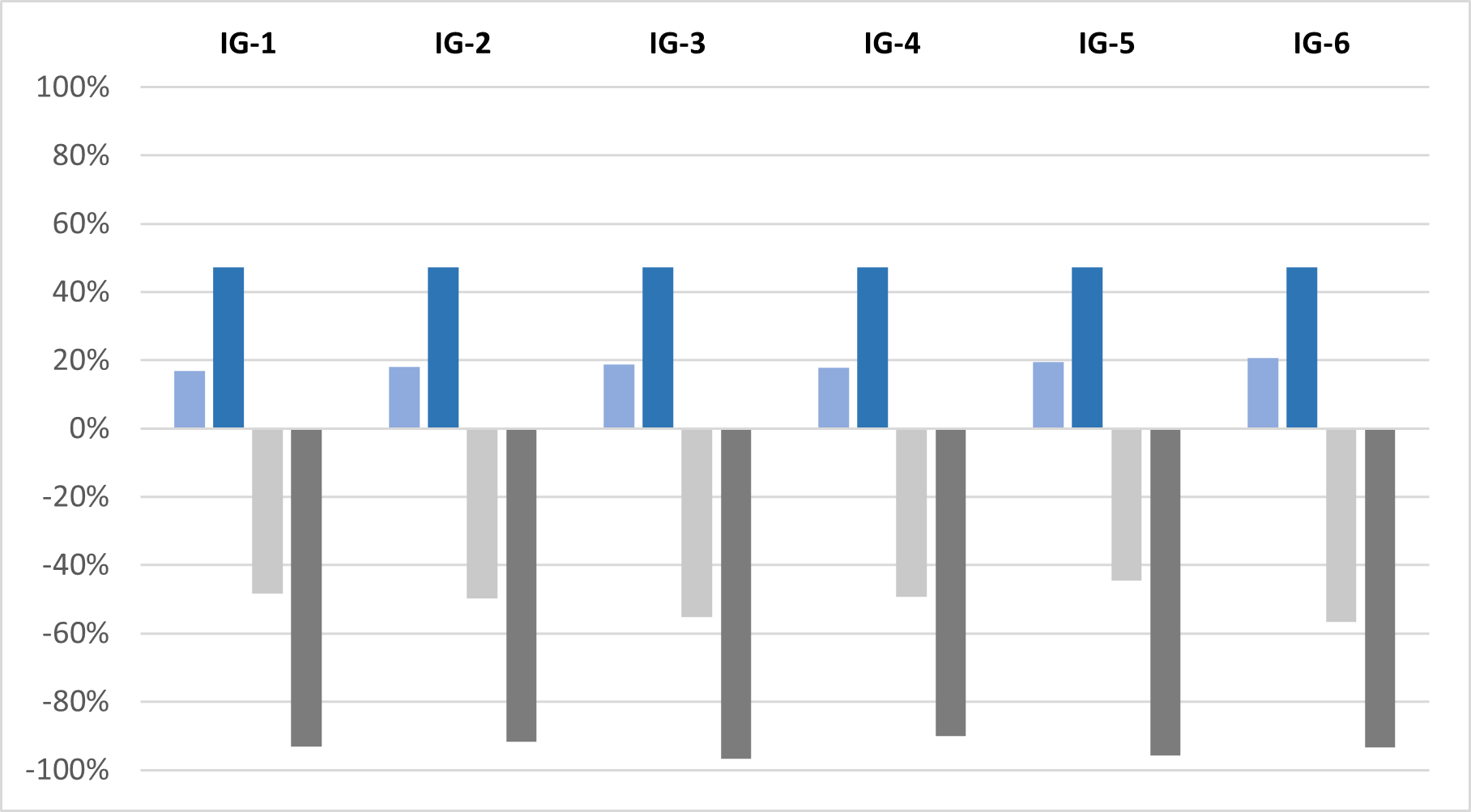}
    \caption{Increase in time window lengths (positive) and decrease in total number of time window violations (negative), by using the higher confidence level $95\%$ $(\beta_\ell = \beta_u = 0.025)$ instead of $90\%$ $(\beta_\ell = \beta_u = 0.05)$ in SM (lighter bars) and RM (darker bars)}
    \label{fig:5}
\end{figure}

Using $\beta_\ell = \beta_u = 0.025$ instead of $\beta_\ell = \beta_u = 0.05$ means that the service provider requires more confidence ($95\%$ instead of $90\%$), or tolerates lower risks, to ensure that the time windows assigned to the customers are not violated either way more than $2.5\%$ of times. Figure \ref{fig:5} displays how such an increase in the confidence level will impact the average length and number of violations of the assigned time windows. It is observed that increasing the confidence level by $5\%$ would result in almost 50$\%$ and 95$\%$  reduction in the number of violations under the stochastic and the robust settings, respectively. However, such a benefit comes at the cost of increased time window lengths by almost $19\%$ and $48\%$ under the SM and RM, respectively. 

\subsubsection{A Guideline to Assign Guaranteed Delivery Time Windows.} 
\begin{figure}[h!]
    \centering
    \includegraphics[scale=0.6]{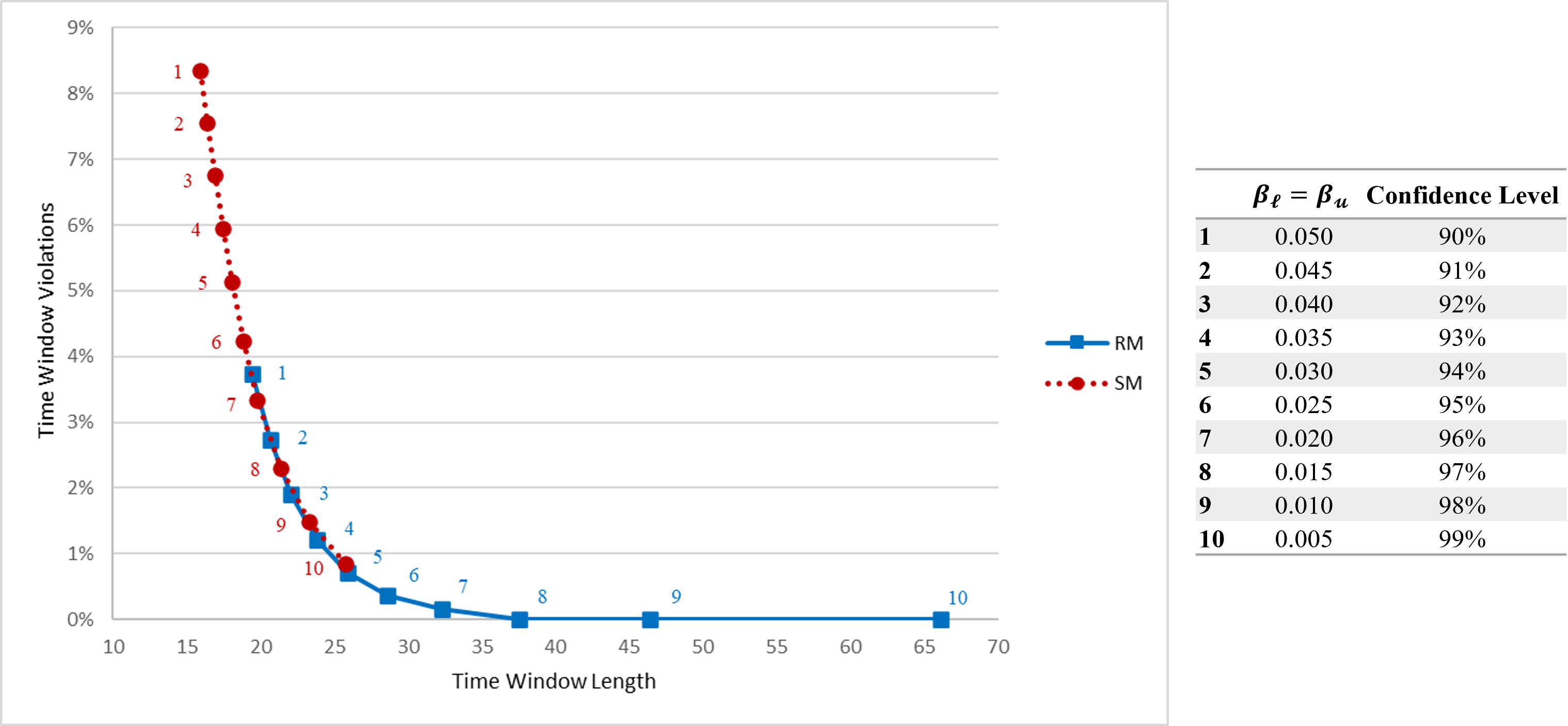}
    \caption{A guideline for selecting appropriate model and confidence level for IG-4}
    \label{fig:6}
\end{figure}

In order to provide a guideline for the service provider in dealing with each problem instance, a chart similar to what is displayed in Figure \ref{fig:6} can be generated. This chart helps the service provider to select the appropriate model and confidence level depending on the acceptable violation rate (or the risk tolerance) on either side of the time windows, as well as the length of the time windows that the service provider considers suitable for its customers. This chart clearly displays the trade-off between the time window violation rate and its length in designing an appropriate route for problem instances IG-4. Whereas customers prefer shorter time windows, service providers favor longer ones to reduce the frequency of early or late arrival violations.

\subsubsection{Impact of Fixed-Length Time windows.}
\begin{figure}[h!]
    \centering
    \includegraphics[scale=0.9]{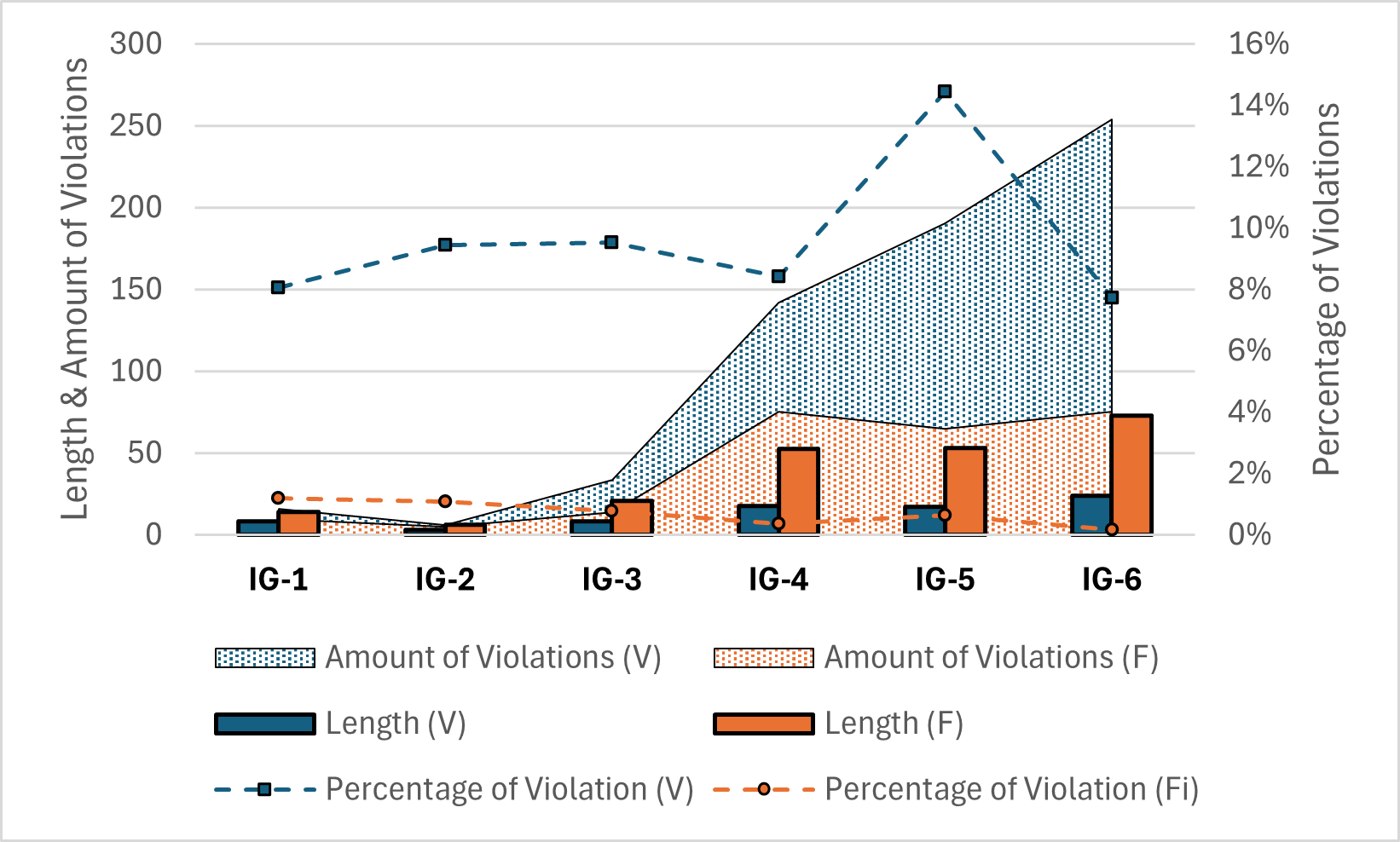}
    \caption{Length vs. percentage and amount of violations for variable-length (V) and fixed-length (F) time windows}
    \label{fig:f_vs_v}
\end{figure}

Figure \ref{fig:f_vs_v} compares the average performance of fixed‐length (Extension \ref{subsec:stochastic_fixed_sol}) versus variable‐length time windows under our stochastic programming approach. In each instance, the variable windows (blue bars) are consistently shorter than the fixed windows (orange bars). As a natural consequence of these narrower intervals, the amount and percentage of time window violations tend to be higher under the variable approach. Nevertheless, the violation rates remain below the $10\%$ risk‐tolerance threshold—equivalent to a $90\%$ confidence level in schedule reliability—in all but the out‐of‐sample IG‑5 case. Hence, while the variable‐window scheme experiences more violations, these remain within acceptable limits, and customers typically benefit from the reliable shorter intervals.

Extending the fixed‐window concept to a DRO setting entails considerable mathematical complexity, so the DRO model developed in this paper continued to focus on variable windows. Overall, both approaches remain viable under the stochastic programming approach, underscoring the trade‐off between tighter, tailored windows and a modestly increased risk of arriving outside those windows.

\subsection{Decomposition Algorithms' Performance}

\begin{figure}[htp]
  \centering
  \subfigure[]{\includegraphics[width=19.4pc]{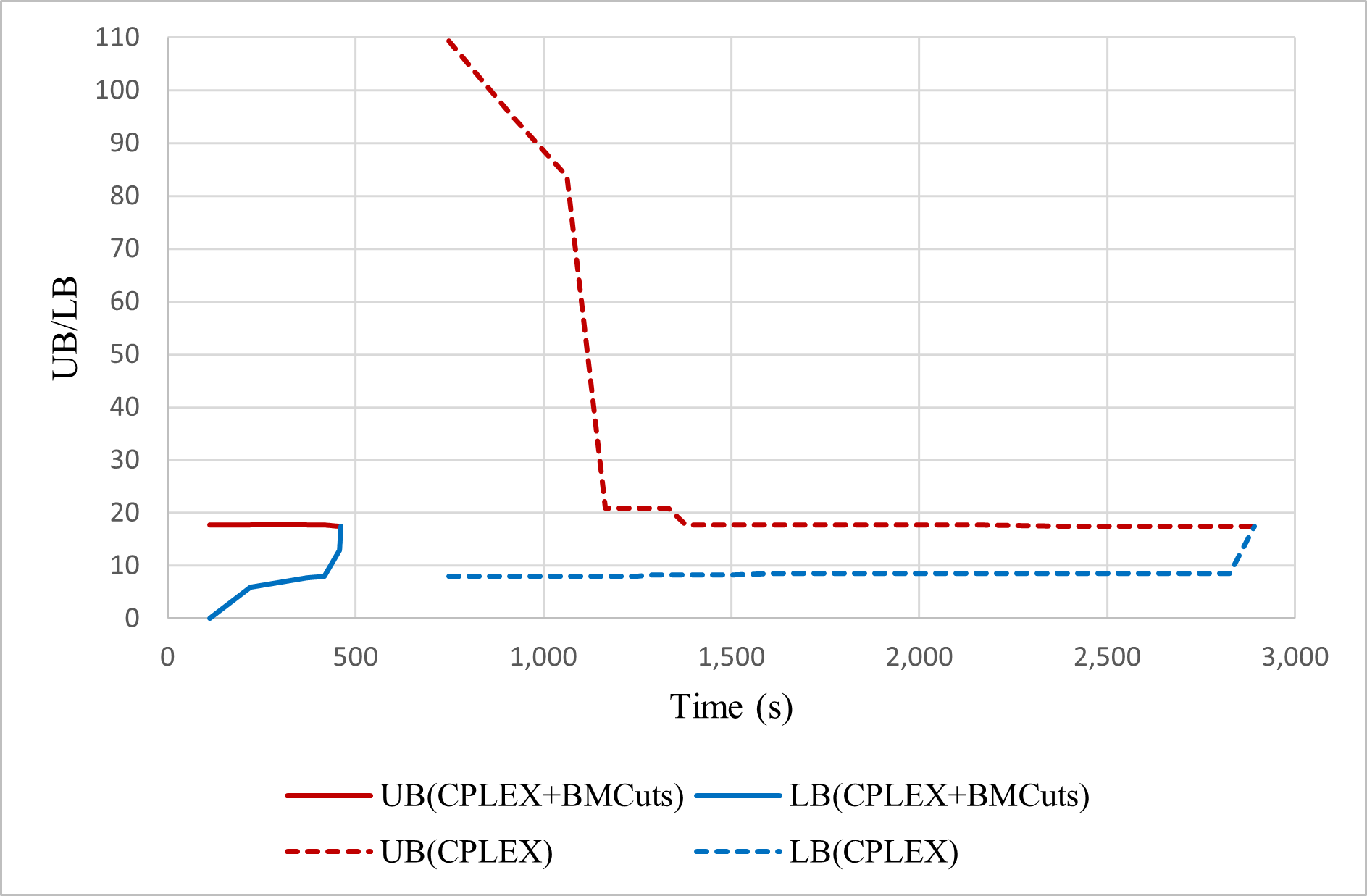}}\label{fig:7-a}
  \subfigure[]{\includegraphics[width=19.4pc]{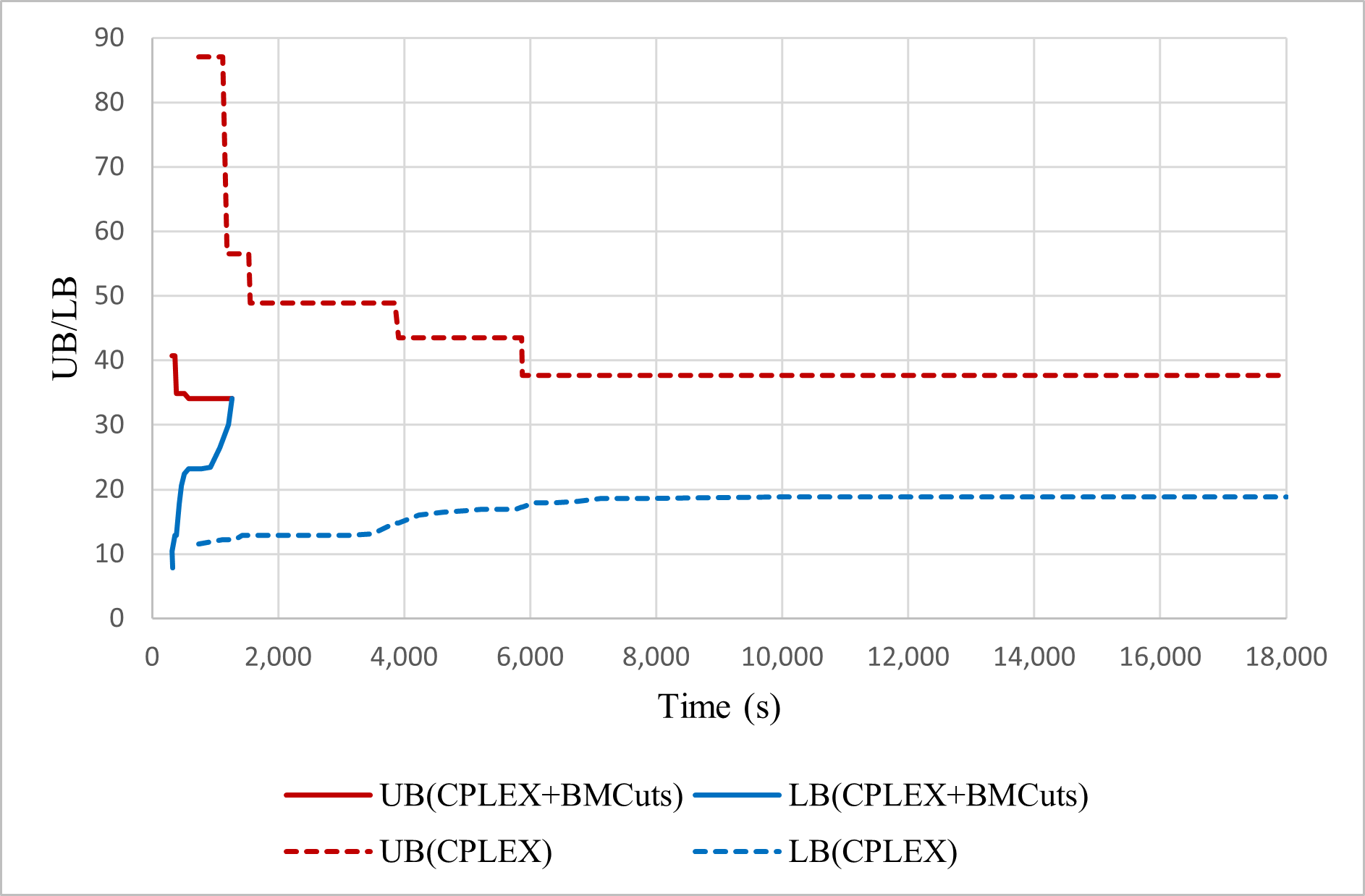}}\label{fig:7-b}
  \caption{Gap between the upper bound and lower bound by CPLEX with and without Benders cuts to solve the SM on complete graphs for (a) 20 customers and (b) 21 customers}
  \label{fig:7}
\end{figure}

\begin{figure}[htp]
  \centering
  \subfigure[]{\includegraphics[width=19.4pc]{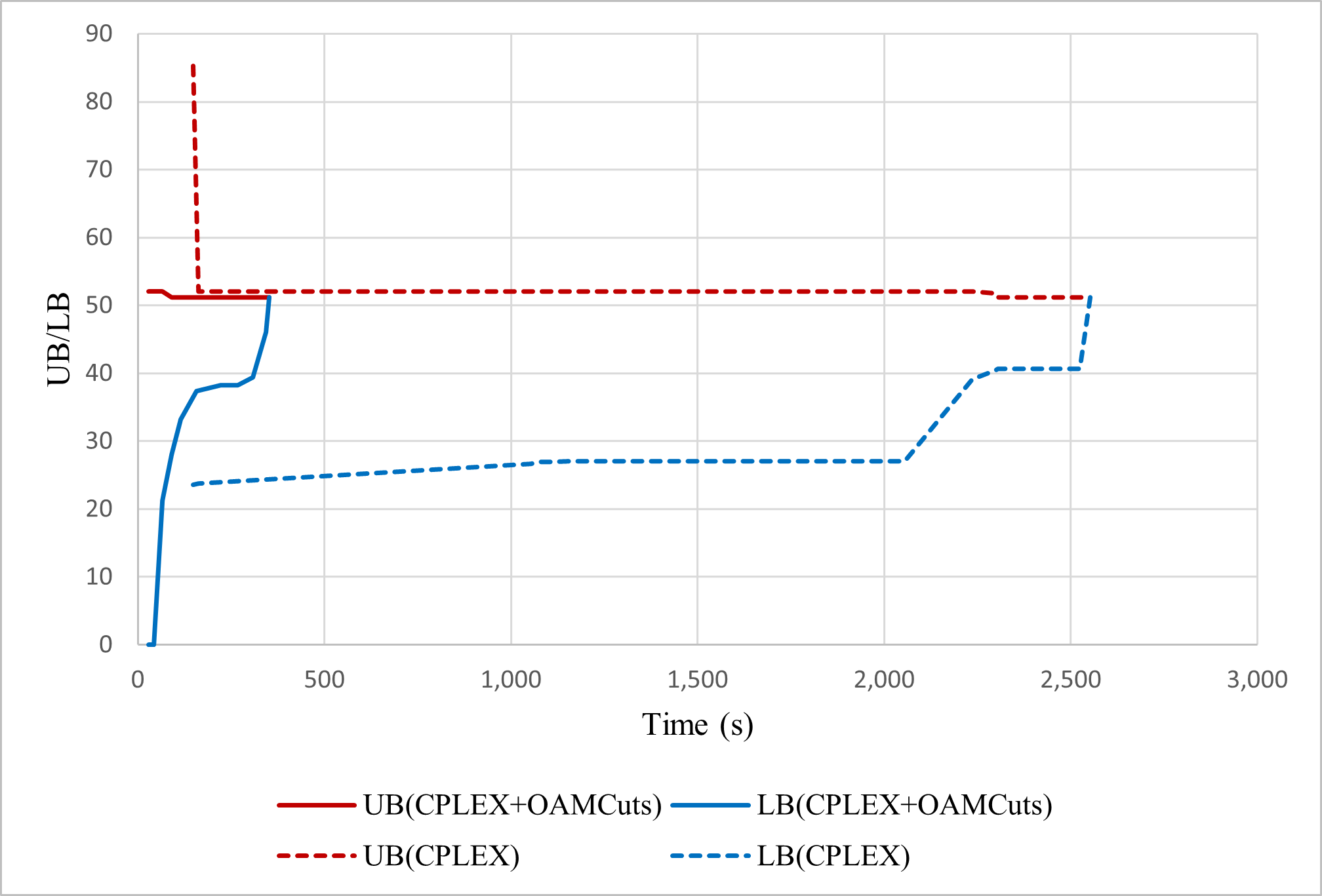}}\label{fig:8-a}
  \subfigure[]{\includegraphics[width=19.4pc]{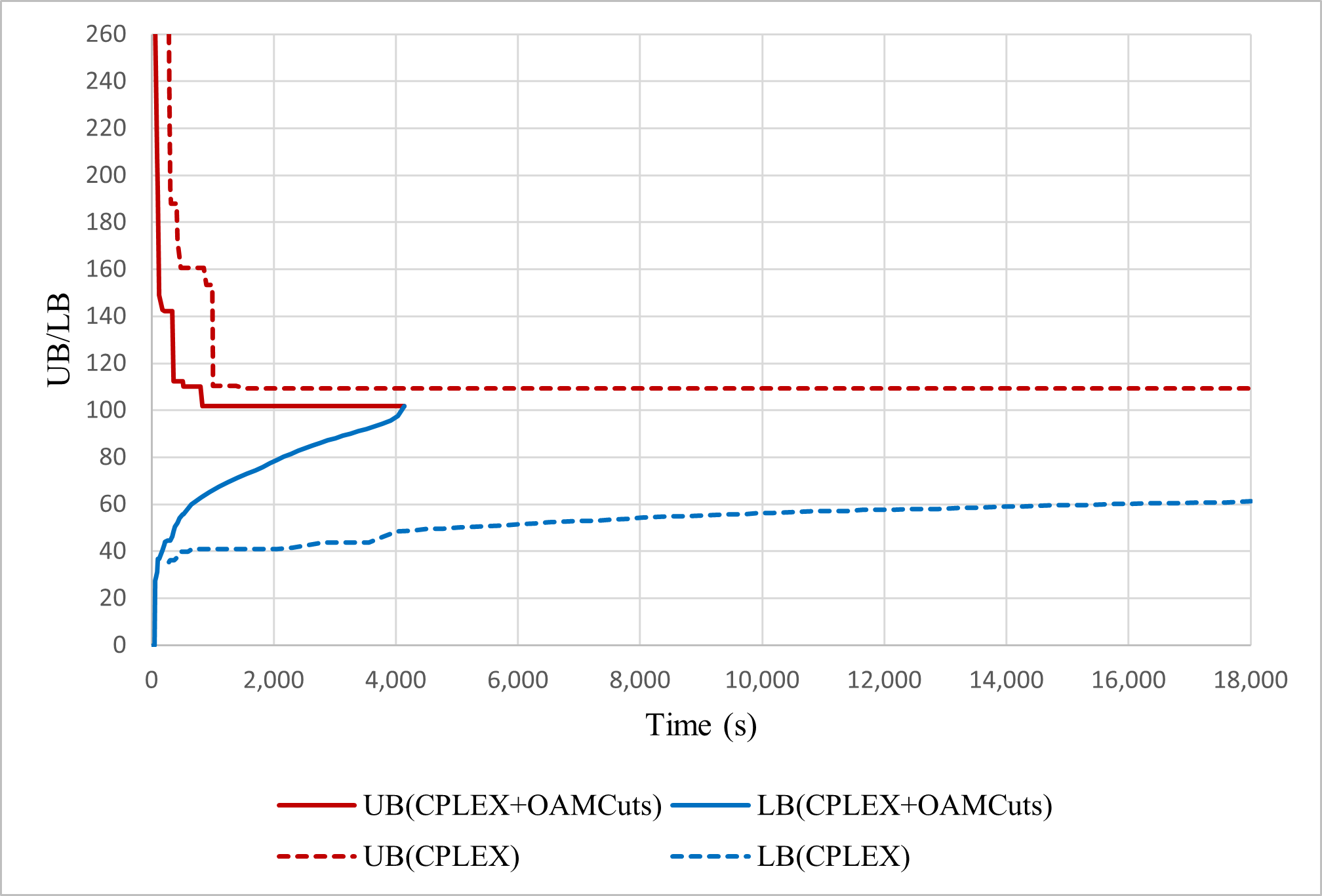}}\label{fig:8-b}
  \caption{Gap between the upper bound and lower bound by CPLEX with and without outer approximation cuts to solve the RM on complete graphs for (a) 20 customers and (b) 21 customers}
  \label{fig:8}
\end{figure}
In this section, we present our computational experiments to evaluate the performance of our proposed decomposition algorithms on dense graphs. The experiments used R101 instances introduced in \cite{rostami2021branch}. These instances are characterized by random customer geographical locations, a complete underlying network, and provided mean vector $\hat{\boldsymbol{\mu}}$ and covariance matrix $\hat{\boldsymbol{\mathcal{C}}}$. These instances were derived from the instances introduced by \cite{solomon1987algorithms} for the VRP with time windows, where the time windows were discarded for our experiments. To ensure statistical significance, and similar to the previous section, a sample size of 1000 was used in the stochastic model.

We evaluated three versions of the CPLEX branch-and-cut algorithm: CPLEX, CPLEX+BSCut, and CPLEX+BMCuts to solve the stochastic model (SM). The first version directly uses CPLEX to solve the Mixed-Integer Programming (MIP) model given in \eqref{eq:SM}. The second and third versions incorporate single Benders cuts and multiple Benders cuts, respectively. To compute the optimality cuts efficiently, we employed the closed-form solution presented in Proposition \ref{Proposition3}. This approach significantly reduced computational time. Similarly, we considered three versions of the CPLEX branch-and-cut algorithm to solve the robust model (RM): CPLEX, CPLEX+OASCut, and CPLEX+OAMCuts. These versions utilize direct CPLEX usage (with the reformulation given in RM'), OA single cuts, and OA multiple cuts, respectively. Through various settings, we found that separating Benders/OA cuts at integer solutions in the tree (as lazy constraints) and considering only fractional solutions at the root node yielded the best performance for the decomposition-based algorithms. A time limit (TL) of 5 hours was set for the experiments.

Detailed results of the algorithms' performance can be found in Appendix F 
through Tables 1 to 4.
To gain insights into the algorithms' behavior, we plotted the improvements of the lower bound (LB) and upper bound (UB) throughout the decision procedure (branch-and-bound) for two instances, one with 20 customers and another with 21 customers. Figure \ref{fig:7} showcases how the inclusion of Benders cuts in the CPLEX branch-and-cut algorithm reduces the computational time required to narrow the optimality gap when solving the SM on complete graphs. Similarly, Figure \ref{fig:8} illustrates the impact of adding outer approximation cuts to the CPLEX branch-and-cut algorithm for solving the RM. Although CPLEX effectively reduces the UB within a reasonable amount of time, the LB improvement is notably slower. In contrast, the addition of cuts noticeably facilitates LB growth, leading to faster convergence with the UB. Without the proposed Benders and OA cuts, CPLEX requires significantly more time to decrease the optimality gap and eventually prove optimality for instances with 20 customers. However, even after running CPLEX for 5 hours, a gap of zero is not achieved for the instance with 21 customers. This gap is closed much sooner when utilizing the cuts generated and added through our proposed decomposition algorithms.

The scalability of our proposed decomposition methods can be further enhanced by incorporating heuristic strategies. General heuristic frameworks like Adaptive Large Neighborhood Search (ALNS) \citep{pisinger2007general} or specialized route construction algorithms for VRP with time windows \citep{braysy2005vehicle} have proven effective in constructing high-quality initial feasible solutions within the master problems. Furthermore, local search methods, commonly used within broader metaheuristic contexts \citep{cordeau2002guide}, can refine incumbent solutions. While the concept of a restricted master problem (RMP) is rooted in \citep{magnanti1981accelerating}, its practical application is often integrated within recent Benders implementations. Finally, advanced cutting plane methods, including lift-and-project cuts \citep{balas1993lift}, can be employed for efficient cut generation and selection. These strategies, when embedded in our frameworks, balance solution quality and computational cost, particularly for large-scale instances.
\section{Conclusion}
\label{sec:conclusion}

For many businesses involved in last-mile operations, providing a high-quality delivery service in terms of reliability is critical for customer satisfaction and retention. In this paper, we proposed a new routing optimization approach with time window assignment using which a service provider can promise reliable goods/service delivery to a set of customers in a network with stochastic and possibly correlated arc travel times. To design such time windows, we have introduced two criteria that address the length of the time windows and the violation risk associated with early and late arrivals to the customers. We have provided two modeling frameworks based on stochastic and distributionally robust optimization and analytically demonstrated how these criteria provide certain levels of service guarantee for the customers. In particular, we have found the closed-form solutions for the optimal time windows in both settings with various risk tolerances when a route is obtained from any source (e.g., delivery routing software), and showed how to later exploit them in developing decomposition-based exact algorithms for solving the integrated routing and design problems. 
 
In our computational experiments, we show how both models are capable of finding routes with reliable time windows for the customers based on the service provider's risk preference. Moreover, the results show that while a small portion of the time windows designed by the stochastic model is violated on the out-of-sample test instances, the distributionally robust model generates more reliable routes and time windows whose violation rates never exceeded the risk tolerance of the service provider on either side. This, however, comes at the cost of assigning longer time windows to the customers. Solving the proposed models could become computationally expensive in a dense network. Thus, we developed two decomposition algorithms based on Benders decomposition and outer approximation to solve the stochastic and distributionally robust models on complete graphs, respectively. Our computational study validated the efficacy of these algorithms in reducing the required computational time to find the optimal solution within a time limit and generating higher quality solutions in the case of acceding to a good integer solution found in a limited time. 

In this study, we introduced two key extensions—incorporating fixed-length time windows and allowing waiting before time windows’ lower bounds—to enhance the practical applicability of our proposed approach. However, these extensions were not explored within a DRO framework. A natural direction for future research is to integrate these enhancements into a DRO setting, where uncertainty in travel times or demand distributions is explicitly accounted for. Investigating how fixed-length time windows and strategic waiting policies interact with distributional uncertainty could provide deeper insights into robust and adaptive decision-making, particularly in time-sensitive and risk-averse applications.

Several other future research avenues can also be considered. While our study aims to design an a priori route that can robustly accommodate uncertainties and variations that may arise in the a posteriori route, future studies may adapt our contributions to a setting that allows for dynamic adjustments of time windows and/or routes. Another one that we plan to consider in the near future is to approximate the proposed models using historical data through machine learning techniques integrated with optimization to predict the arc travel times. We believe this will improve the reliability of the designed time windows and that the proposed algorithms can be extended to deal with the new models.


\begin{thebibliography}{60}
\providecommand{\natexlab}[1]{#1}
\providecommand{\url}[1]{\texttt{#1}}
\providecommand{\urlprefix}{URL }

\bibitem[{Adulyasak \protect\BIBand{} Jaillet(2016)}]{adulyasak2016models}
Adulyasak Y, Jaillet P (2016) Models and algorithms for stochastic and robust vehicle routing with deadlines. \emph{Transportation Science} 50(2):608--626.

\bibitem[{Agrawal et~al.(2012)Agrawal, Ding, Saberi, \protect\BIBand{} Ye}]{agrawal2012price}
Agrawal S, Ding Y, Saberi A, Ye Y (2012) Price of correlations in stochastic optimization. \emph{Operations Research} 60(1):150--162.

\bibitem[{Alizadeh \protect\BIBand{} Goldfarb(2003)}]{alizadeh2003second}
Alizadeh F, Goldfarb D (2003) Second-order cone programming. \emph{Mathematical Programming} 95(1):3--51.

\bibitem[{Amazon(2023)}]{Amazon2023}
Amazon (2023) Estimated delivery windows. accessed: 01.05.2023. \urlprefix\url{https://www.amazon.com/gp/help/customer/display.html?nodeId=GK8CZJ8DR2J2WS5H}.

\bibitem[{Bakach et~al.(2021)Bakach, Campbell, Ehmke, \protect\BIBand{} Urban}]{bakach2021solving}
Bakach I, Campbell AM, Ehmke JF, Urban TL (2021) Solving vehicle routing problems with stochastic and correlated travel times and makespan objectives. \emph{EURO Journal on Transportation and Logistics} 10:100029.

\bibitem[{Balas et~al.(1993)Balas, Ceria, \protect\BIBand{} Cornu{\'e}jols}]{balas1993lift}
Balas E, Ceria S, Cornu{\'e}jols G (1993) A lift-and-project cutting plane algorithm for mixed 0--1 programs. \emph{Mathematical programming} 58(1):295--324.

\bibitem[{Ben-Tal \protect\BIBand{} Nemirovski(2001)}]{ben2001polyhedral}
Ben-Tal A, Nemirovski A (2001) On polyhedral approximations of the second-order cone. \emph{Mathematics of Operations Research} 26(2):193--205.

\bibitem[{Boysen et~al.(2021)Boysen, Fedtke, \protect\BIBand{} Schwerdfeger}]{boysen2021last}
Boysen N, Fedtke S, Schwerdfeger S (2021) Last-mile delivery concepts: a survey from an operational research perspective. \emph{Or Spectrum} 43:1--58.

\bibitem[{Br{\"a}ysy \protect\BIBand{} Gendreau(2005)}]{braysy2005vehicle}
Br{\"a}ysy O, Gendreau M (2005) Vehicle routing problem with time windows, part i: Route construction and local search algorithms. \emph{Transportation science} 39(1):104--118.

\bibitem[{Carlsson \protect\BIBand{} Delage(2013)}]{carlsson2013robust}
Carlsson JG, Delage E (2013) Robust partitioning for stochastic multivehicle routing. \emph{Operations Research} 61(3):727--744.

\bibitem[{Cordeau et~al.(2002)Cordeau, Gendreau, Laporte, Potvin, \protect\BIBand{} Semet}]{cordeau2002guide}
Cordeau JF, Gendreau M, Laporte G, Potvin JY, Semet F (2002) A guide to vehicle routing heuristics. \emph{Journal of the Operational Research society} 53(5):512--522.

\bibitem[{Cui et~al.(2020)Cui, Lu, Sun, \protect\BIBand{} Golden}]{cui2020sooner}
Cui R, Lu Z, Sun T, Golden J (2020) Sooner or later? promising delivery speed in online retail. \emph{March 29, https://dx.doi.org/10.2139/ssrn.3563404} .

\bibitem[{Dayarian \protect\BIBand{} Savelsbergh(2020)}]{dayarian2020crowdshipping}
Dayarian I, Savelsbergh M (2020) Crowdshipping and same-day delivery: Employing in-store customers to deliver online orders. \emph{Production and Operations Management} 29(9):2153--2174.

\bibitem[{Delage \protect\BIBand{} Ye(2010)}]{delage2010distributionally}
Delage E, Ye Y (2010) Distributionally robust optimization under moment uncertainty with application to data-driven problems. \emph{Operations Research} 58(3):595--612.

\bibitem[{Deloitte(2020)}]{Deloitte2020last}
Deloitte (2020) Last mile logistics, challenges and solutions in spain. \emph{Deloitte consulting, Department of Marketing \& Brand} .

\bibitem[{Deng et~al.(2021)Deng, Fang, \protect\BIBand{} Lim}]{deng2021urban}
Deng Q, Fang X, Lim YF (2021) Urban consolidation center or peer-to-peer platform? the solution to urban last-mile delivery. \emph{Production and Operations Management} 30(4):997--1013.

\bibitem[{DispatchTrack(2022)}]{DispatchTrack2022big}
DispatchTrack (2022) Big and bulky delivery. \emph{DispatchTrack Report} .

\bibitem[{European~Commission \protect\BIBand{} SMEs(2022)}]{Euro2022}
European~Commission IE Directorate-General for Internal~Market, SMEs (2022) Domestic postal traffic, letter mail and parcel services. \urlprefix\url{http://data.europa.eu/88u/dataset/6rFQHnqYW7HDFi1hIuDHg}.

\bibitem[{Fatehi \protect\BIBand{} Wagner(2022)}]{fatehi2022crowdsourcing}
Fatehi S, Wagner MR (2022) Crowdsourcing last-mile deliveries. \emph{Manufacturing \& Service Operations Management} 24(2):791--809.

\bibitem[{Gendreau et~al.(2014)Gendreau, Jabali, \protect\BIBand{} Rei}]{gendreau2014chapter}
Gendreau M, Jabali O, Rei W (2014) Chapter 8: Stochastic vehicle routing problems. \emph{Vehicle Routing: Problems, Methods, and Applications, Second Edition}, 213--239 (SIAM).

\bibitem[{Geoffrion(1972)}]{geoffrion1972generalized}
Geoffrion AM (1972) Generalized benders decomposition. \emph{Journal of Optimization Theory and Applications} 10(4):237--260.

\bibitem[{Hoogeboom et~al.(2021)Hoogeboom, Adulyasak, Dullaert, \protect\BIBand{} Jaillet}]{hoogeboom2021robust}
Hoogeboom M, Adulyasak Y, Dullaert W, Jaillet P (2021) The robust vehicle routing problem with time window assignments. \emph{Transportation Science} 55(2):395--413.

\bibitem[{Jabali et~al.(2015)Jabali, Leus, Van~Woensel, \protect\BIBand{} De~Kok}]{jabali2015self}
Jabali O, Leus R, Van~Woensel T, De~Kok T (2015) Self-imposed time windows in vehicle routing problems. \emph{OR Spectrum} 37(2):331--352.

\bibitem[{Jaillet et~al.(2016)Jaillet, Qi, \protect\BIBand{} Sim}]{jaillet2016routing}
Jaillet P, Qi J, Sim M (2016) Routing optimization under uncertainty. \emph{Operations Research} 64(1):186--200.

\bibitem[{Laporte(2010)}]{laporte2010concise}
Laporte G (2010) A concise guide to the traveling salesman problem. \emph{Journal of the Operational Research Society} 61(1):35--40.

\bibitem[{Lecluyse et~al.(2009)Lecluyse, Van~Woensel, \protect\BIBand{} Peremans}]{lecluyse2009vehicle}
Lecluyse C, Van~Woensel T, Peremans H (2009) Vehicle routing with stochastic time-dependent travel times. \emph{4OR} 7(4):363--377.

\bibitem[{Letchford \protect\BIBand{} Nasiri(2015)}]{Letchford2015steiner}
Letchford AN, Nasiri SD (2015) The steiner travelling salesman problem with correlated costs. \emph{European Journal of Operational Research} 245(1):62--69.

\bibitem[{Lim et~al.(2023)Lim, Wang, \protect\BIBand{} Webster}]{lim2023right}
Lim SFW, Wang Q, Webster S (2023) Do it right the first time: Vehicle routing with home delivery attempt predictors. \emph{Production and Operations Management} 32(4):1262--1284.

\bibitem[{Liu et~al.(2019)Liu, Wang, \protect\BIBand{} Susilo}]{liu2019assessing}
Liu C, Wang Q, Susilo YO (2019) Assessing the impacts of collection-delivery points to individual’s activity-travel patterns: A greener last mile alternative? \emph{Transportation Research Part E: Logistics and Transportation Review} 121:84--99.

\bibitem[{Liu et~al.(2021)Liu, He, \protect\BIBand{} Max~Shen}]{liu2021time}
Liu S, He L, Max~Shen ZJ (2021) On-time last-mile delivery: Order assignment with travel-time predictors. \emph{Management Science} 67(7):4095--4119.

\bibitem[{Loqate(2022)}]{loqate2022fixing}
Loqate (2022) Fixing failed deliveries, stamping out faulty fulfilment. \emph{Loqate GBG Report} .

\bibitem[{Lyu \protect\BIBand{} Teo(2022)}]{lyu2022last}
Lyu G, Teo CP (2022) Last mile innovation: The case of the locker alliance network. \emph{Manufacturing \& Service Operations Management} 24(5):2425--2443.

\bibitem[{Macioszek(2018)}]{macioszek2018first}
Macioszek E (2018) First and last mile delivery--problems and issues. \emph{Advanced Solutions of Transport Systems for Growing Mobility: 14th Scientific and Technical Conference" Transport Systems. Theory \& Practice 2017" Selected Papers}, 147--154 (Springer).

\bibitem[{Magnanti \protect\BIBand{} Wong(1981)}]{magnanti1981accelerating}
Magnanti TL, Wong RT (1981) Accelerating benders decomposition: Algorithmic enhancement and model selection criteria. \emph{Operations research} 29(3):464--484.

\bibitem[{Mangiaracina et~al.(2019)Mangiaracina, Perego, Seghezzi, \protect\BIBand{} Tumino}]{mangiaracina2019innovative}
Mangiaracina R, Perego A, Seghezzi A, Tumino A (2019) Innovative solutions to increase last-mile delivery efficiency in b2c e-commerce: a literature review. \emph{International Journal of Physical Distribution \& Logistics Management} .

\bibitem[{Martins et~al.(2019)Martins, Ostermeier, Amorim, H{\"u}bner, \protect\BIBand{} Almada-Lobo}]{martins2019product}
Martins S, Ostermeier M, Amorim P, H{\"u}bner A, Almada-Lobo B (2019) Product-oriented time window assignment for a multi-compartment vehicle routing problem. \emph{European Journal of Operational Research} 276(3):893--909.

\bibitem[{Merchan et~al.(2022)Merchan, Arora, Pachon, Konduri, Winkenbach, Parks, \protect\BIBand{} Noszek}]{merchan20222021}
Merchan D, Arora J, Pachon J, Konduri K, Winkenbach M, Parks S, Noszek J (2022) 2021 amazon last mile routing research challenge: Data set. \emph{Transportation Science} Articles in Advance.

\bibitem[{Mohajerin~Esfahani \protect\BIBand{} Kuhn(2018)}]{mohajerin2018data}
Mohajerin~Esfahani P, Kuhn D (2018) Data-driven distributionally robust optimization using the wasserstein metric: Performance guarantees and tractable reformulations. \emph{Mathematical Programming} 171(1):115--166.

\bibitem[{Nicholson(2015)}]{Nicholson201514}
Nicholson A (2015) Travel time reliability benefits: Allowing for correlation. \emph{Research in Transportation Economics} 49:14 -- 21, ISSN 0739-8859.

\bibitem[{Parent \protect\BIBand{} LeSage(2010)}]{parent2010spatial}
Parent O, LeSage JP (2010) A spatial dynamic panel model with random effects applied to commuting times. \emph{Transportation Research Part B: Methodological} 44(5):633--645.

\bibitem[{Pisinger \protect\BIBand{} Ropke(2007)}]{pisinger2007general}
Pisinger D, Ropke S (2007) A general heuristic for vehicle routing problems. \emph{Computers \& operations research} 34(8):2403--2435.

\bibitem[{Qi et~al.(2018)Qi, Li, Liu, \protect\BIBand{} Max~Shen}]{qi2018shared}
Qi W, Li L, Liu S, Max~Shen ZJ (2018) Shared mobility for last-mile delivery: Design, operational prescriptions, and environmental impact. \emph{Manufacturing \& Service Operations Management} 20(4):737--751.

\bibitem[{Rajabi-Bahaabadi et~al.(2019)Rajabi-Bahaabadi, Shariat-Mohaymany, Babaei, \protect\BIBand{} Vigo}]{rajabi2019reliable}
Rajabi-Bahaabadi M, Shariat-Mohaymany A, Babaei M, Vigo D (2019) Reliable vehicle routing problem in stochastic networks with correlated travel times. \emph{Operational Research} 1--32.

\bibitem[{Rostami et~al.(2021)Rostami, Desaulniers, Errico, \protect\BIBand{} Lodi}]{rostami2021branch}
Rostami B, Desaulniers G, Errico F, Lodi A (2021) Branch-price-and-cut algorithms for the vehicle routing problem with stochastic and correlated travel times. \emph{Operations Research} 69(2):436--455.

\bibitem[{Salari et~al.(2022)Salari, Liu, \protect\BIBand{} Shen}]{salari2022real}
Salari N, Liu S, Shen ZJM (2022) Real-time delivery time forecasting and promising in online retailing: when will your package arrive? \emph{Manufacturing \& Service Operations Management} 24(3):1421--1436.

\bibitem[{Savelsbergh \protect\BIBand{} Van~Woensel(2016)}]{Savelsbergh2016city}
Savelsbergh M, Van~Woensel T (2016) 50th anniversary invited article—city logistics: Challenges and opportunities. \emph{Transportation Science} 50(2):579–590.

\bibitem[{Scarf(1958)}]{scarf1958min}
Scarf H (1958) A min-max solution of an inventory problem. \emph{Studies in the Mathematical Theory of Inventory and Production} .

\bibitem[{Seshadri \protect\BIBand{} Srinivasan(2012)}]{seshadri2012algorithm}
Seshadri R, Srinivasan KK (2012) An algorithm for the minimum robust cost path on networks with random and correlated link travel times. Levinson DM, Liu HX, Bell M, eds., \emph{Network Reliability in Practice}, 171--208 (New York, NY: Springer), ISBN 978-1-4614-0947-2, \urlprefix\url{http://dx.doi.org/10.1007/978-1-4614-0947-2_11}.

\bibitem[{Solomon(1987)}]{solomon1987algorithms}
Solomon MM (1987) Algorithms for the vehicle routing and scheduling problems with time window constraints. \emph{Operations Research} 35(2):254--265.

\bibitem[{Spliet et~al.(2018)Spliet, Dabia, \protect\BIBand{} Van~Woensel}]{spliet2018time}
Spliet R, Dabia S, Van~Woensel T (2018) The time window assignment vehicle routing problem with time-dependent travel times. \emph{Transportation Science} 52(2):261--276.

\bibitem[{Spliet \protect\BIBand{} Desaulniers(2015)}]{spliet2015discrete}
Spliet R, Desaulniers G (2015) The discrete time window assignment vehicle routing problem. \emph{European Journal of Operational Research} 244(2):379--391.

\bibitem[{Spliet \protect\BIBand{} Gabor(2015)}]{spliet2015time}
Spliet R, Gabor AF (2015) The time window assignment vehicle routing problem. \emph{Transportation Science} 49(4):721--731.

\bibitem[{Subramanyam et~al.(2018)Subramanyam, Wang, \protect\BIBand{} Gounaris}]{subramanyam2018scenario}
Subramanyam A, Wang A, Gounaris CE (2018) A scenario decomposition algorithm for strategic time window assignment vehicle routing problems. \emph{Transportation Research Part B: Methodological} 117:296--317.

\bibitem[{Ulmer et~al.(2024)Ulmer, Goodson, \protect\BIBand{} Thomas}]{ulmer2024optimal}
Ulmer MW, Goodson JC, Thomas BW (2024) Optimal service time windows. \emph{Transportation Science} 58(2):394--411.

\bibitem[{Van~Loon et~al.(2015)Van~Loon, Deketele, Dewaele, McKinnon, \protect\BIBand{} Rutherford}]{van2015comparative}
Van~Loon P, Deketele L, Dewaele J, McKinnon A, Rutherford C (2015) A comparative analysis of carbon emissions from online retailing of fast moving consumer goods. \emph{Journal of Cleaner Production} 106:478--486.

\bibitem[{Vareias et~al.(2019)Vareias, Repoussis, \protect\BIBand{} Tarantilis}]{vareias2019assessing}
Vareias AD, Repoussis PP, Tarantilis CD (2019) Assessing customer service reliability in route planning with self-imposed time windows and stochastic travel times. \emph{Transportation Science} 53(1):256--281.

\bibitem[{WaterlooRegionalPolice(2023)}]{WaterlooPolice2023}
WaterlooRegionalPolice (2023) Increased parcel thefts. \urlprefix\url{https://x.com/i/web/status/1733884476269248873}.

\bibitem[{Wiesemann et~al.(2014)Wiesemann, Kuhn, \protect\BIBand{} Sim}]{wiesemann2014distributionally}
Wiesemann W, Kuhn D, Sim M (2014) Distributionally robust convex optimization. \emph{Operations Research} 62(6):1358--1376.

\bibitem[{Yu et~al.(2023)Yu, Shen, Badri-Koohi, \protect\BIBand{} Seada}]{yu2023time}
Yu X, Shen S, Badri-Koohi B, Seada H (2023) Time window optimization for attended home service delivery under multiple sources of uncertainties. \emph{Computers \& Operations Research} 150:106045.

\bibitem[{Zhang et~al.(2024)Zhang, Zhang, \protect\BIBand{} Baldacci}]{zhang2024generalized}
Zhang Z, Zhang Y, Baldacci R (2024) Generalized riskiness index in vehicle routing under uncertain travel times: Formulations, properties, and exact solution framework. \emph{Transportation Science} .

\end{thebibliography}


\newpage
\setcounter{page}{1}

\cleardoublepage
\phantomsection
\pdfbookmark[1]{Online Supplement}{supplement}

\begin{center}
\Huge \bfseries Online Supplement
\end{center}

\vspace{1cm}

\begin{APPENDICES}
    \phantomsection
\pdfbookmark[2]{Appendix A: Proof of Proposition \ref{Proposition1}}{appA}
\section*{Appendix A: Proof of Proposition \ref{Proposition1}}
\label{sec:AppB}
To prove the proposition, we use the results of the following lemma.
\begin{lemma}
    With $\boldsymbol y$ fixed, let $H(\epsilon) = \mathbb{E}_{\mathrm{P}} \Big[\big(\tau(\boldsymbol{y}, \boldsymbol{t})-\epsilon\big)^{+} \Big] = \int_{t \in \mathbb{R}^{|A|}} \big(\tau(\boldsymbol{y}, \boldsymbol{t})-\epsilon\big)^{+} p(\boldsymbol{t}) \mathrm{d} \boldsymbol{t}$. Then, assuming $\tau(\boldsymbol{y}, \boldsymbol{t})$ has a non-atomic (continuous) distribution, $H(\epsilon)$ is convex and continuously differentiable in $\epsilon$ with
    $$
    H'(\epsilon) = 
    F(\boldsymbol{y}, \epsilon) - 1.
    $$
\end{lemma}
 
\emph{Proof.}

\textbf{Convexity.} Fix any realization \(\boldsymbol{t}\). The function $\big(\tau(\boldsymbol{y},\boldsymbol{t}) \;-\;\epsilon\big)^+ = \max\bigl\{\,0,\;\tau(\boldsymbol{y},\boldsymbol{t})-\epsilon\bigr\}$ is a maximum of two affine (linear) functions in $\epsilon$, so it is convex. Then $H(\epsilon)$, being an expectation of this function over $\boldsymbol{t}$, remains convex (since integrals preserve convexity).

\textbf{Differentiability.}
For each fixed \(\boldsymbol{t}\),
\[
  \frac{\partial}{\partial \epsilon}\,
  \big(\tau(\boldsymbol{y},\boldsymbol{t})-\epsilon\big)^+
  \;=\;
  \begin{cases}
    -1, & \text{if } \tau(\boldsymbol{y},\boldsymbol{t}) > \epsilon,\\
     0, & \text{if } \tau(\boldsymbol{y},\boldsymbol{t}) < \epsilon,
  \end{cases}
\]
and is undefined (a kink) only when \(\epsilon = \tau(\boldsymbol{y},\boldsymbol{t})\). Because we assume the distribution of \(\tau(\boldsymbol{y},\boldsymbol{t})\) has no point masses (atoms), i.e.,
\[
  \Pr\big(\tau(\boldsymbol{y},\boldsymbol{t}) \;=\;\epsilon\big) 
  \;=\; 0
  \quad
  \text{for all }\epsilon,
\]
the set of \(\boldsymbol{t}\) that causes the kink has measure zero, so the above derivative exists \emph{almost everywhere} (a.e.) in \(\epsilon\). Next, by standard dominated convergence arguments, we can pass the (sub)derivative inside the expectation, so:
\[
  H'(\epsilon)
  \;=\;\mathbb{E}_{\mathrm{P}}\!\Bigl[
       \frac{\partial}{\partial \epsilon}\,
       \big(\tau(\boldsymbol{y},\boldsymbol{t})-\epsilon\big)^+
     \Bigr]
  \;=\;-\,\Pr\big(\tau(\boldsymbol{y},\boldsymbol{t})>\epsilon\big).
\]
Thus \(H(\epsilon)\) is differentiable almost everywhere.

\textbf{Continuity.} Finally, due to the non-atomic distribution of \(\tau(\boldsymbol{y},\boldsymbol{t})\), its cumulative distribution function $F(\boldsymbol{y}, \epsilon)$ is continuous, and hence the function $H'(\epsilon) = F(\boldsymbol{y}, \epsilon)-1$ is a continuous function of \(\epsilon\), making \(H(\epsilon)\) continuously differentiable in~\(\epsilon\).
\hfill $\square$

Likewise, under the usual non‐atomic assumptions, it can be shown that if $H(\epsilon) = \mathbb{E}_{\mathrm{P}} \Big[\big(\epsilon - \tau(\boldsymbol{y}, \boldsymbol{t})\big)^{+} \Big]$, then $H(\epsilon)$ will will be a convex continuously differentiable function with $ H'(\epsilon) = F(\boldsymbol{y}, \epsilon)$.

Consequently, $ \mathcal{H}_\ell^k(\boldsymbol{y}^k, \ell^k) = \mathbb{E}_{\mathrm{P}} \Big[\big(\ell^k - \tau^k(\boldsymbol{y}^k, \boldsymbol{\tilde{t}})\big)^+\Big] $ and $ \mathcal{H}_u^k(\boldsymbol{y}^k, u^k) = \mathbb{E}_{\mathrm{P}} \Big[\big(\tau^k(\boldsymbol{y}^k, \boldsymbol{\tilde{t}}) - u^k\big)^+\Big] $ are convex in \( \ell^k \) and \( u^k \), respectively. Since the linear term \( a_w^k (u^k - \ell^k) \) is affine (hence convex), summing convex terms with non-negative coefficients \( a_u^k, a_\ell^k \geq 0 \) preserves convexity. Thus, \( \mathcal{H}^k(\boldsymbol{y}^k, \ell^k, u^k) \) is jointly convex in \( (\ell^k, u^k) \).

Moreover, the full function \( \mathcal{H}^k \) is a sum:  
\[
\mathcal{H}^k = \underbrace{a_w^k (u^k - \ell^k)}_{\text{affine}} + \underbrace{a_\ell^k \mathcal{H}_\ell^k}_{\text{continuously differentiable in } \ell^k} + \underbrace{a_u^k \mathcal{H}_u^k}_{\text{continuously differentiable in } u^k},
\]   
whose partial derivatives are (by Lemma 1):  
\[
\frac{\partial}{\partial \ell^k} \mathcal{H}^k = -a_w^k + a_\ell^k \underbrace{F^k(\boldsymbol{y}^k, \ell^k)}_{\text{continuous in } \ell^k}, \quad
\frac{\partial}{\partial u^k} \mathcal{H}^k = a_w^k + a_u^k \underbrace{(F^k(\boldsymbol{y}^k, u^k) - 1)}_{\text{continuous in } u^k},
\]  
both of which exist, and since \( F(\boldsymbol{y}^k, \cdot) \) is continuous (no atoms), both partial derivatives are continuous (almost) everywhere in \( \ell^k \) and \( u^k \), respectively.

Furthermore, we investigate the Hessian matrix of \(\mathcal{H}^k\) as follows:  
\[
\nabla^2 \mathcal{H}^k = 
\begin{bmatrix}
a_\ell^k f^k(\boldsymbol{y}^k, \ell^k) & 0 \\
0 & a_u^k f^k(\boldsymbol{y}^k, u^k)
\end{bmatrix},
\]  
where \(f^k(\boldsymbol{y}^k, \ell^k)\) is the PDF of \(\tau^k\) evaluated at \(\ell^k\) and \(f^k(\boldsymbol{y}^k, u^k)\) is the PDF of \(\tau^k\) evaluated at \(u^k\). It is obvious that the Hessian has a diagonal structure as the cross-derivatives \(\frac{\partial^2}{\partial \ell^k \partial u^k} \mathcal{H}^k\) and \(\frac{\partial^2}{\partial u^k \partial \ell^k} \mathcal{H}^k\) are zero because \(\mathcal{H}^k\) is separable in \(u^k\) and \(\ell^k\).  Also, the diagonal entries \(a_\ell^k f^k(\boldsymbol{y}^k, \ell^k)\) and \(a_u^k f^k(\boldsymbol{y}^k, u^k)\) are non-negative since \(a_\ell^k \geq 0\), \(a_u^k \geq 0\), and \(f^k \geq 0\) (because \(f^k\) is a PDF).  This makes the Hessian a diagonal matrix with non-negative entries, which is positive semi-definite.  In conclusion, the Hessian also confirms \(\mathcal{H}^k\) is convex, which guarantees the global optimality of solutions derived from the first-order conditions.
\hfill $\square$
    \vspace{0.5cm}
    \phantomsection
\pdfbookmark[2]{Appendix B: Proof of Proposition \ref{Proposition3}}{appB}
\section*{Appendix B: Proof of Proposition \ref{Proposition3}}
\label{sec:AppC}

To find an optimal value of $[\ell^k,u^k]$ for each $k\in V_0$, we write the dual of $\tilde{\text{SP}}^k$ in \eqref{eq:SP(SM)} by ignoring the constraint $\ell^k \leq u^k$ for now, and construct the primal and dual solutions that satisfy the strong duality. The dual for fixed $\boldsymbol{y}^k$ reads as
\begin{subequations}
\label{eq:DSP(SM)}
\begin{align}
    \tilde{\text {DSP}}^k(\boldsymbol{y}^k): &\max _{\boldsymbol{\rho}_{1}, \boldsymbol{\rho}_{2}} \sum_{k \in V_{0}} \sum_{q=1}^{Q}\left(\sum_{(i, j) \in A} t_{i j}^{[q]} {y}_{i j}^{k}\right) \rho_{2}^{k[q]}-\sum_{k \in V_{0}} \sum_{q=1}^{Q}\left(\sum_{(i, j) \in A} t_{i j}^{[q]} {y}_{i j}^{k}\right) \rho_{1}^{k[q]} 
    \label{eq:DSP(SM)_of}\\
    \mbox{s.t.} \quad
    &\sum_{q=1}^{Q} \rho_{1}^{k[q]} \geq a_w^k \quad \forall k \in V_{0} 
    \label{eq:DSP(SM)_rho1} \\
    &\sum_{q=1}^{Q} \rho_{2}^{k[q]} \leq a_w^k \quad \forall k \in V_{0} 
    \label{eq:DSP(SM)_rho2} \\
    &0 \leq \rho_{1}^{k[q]} \leq \frac{a_\ell^k}{Q} \quad \forall k \in V_{0}, \forall q \in\{1,2, \ldots, Q\} 
    \label{eq:DSP(SM)-rho1-sign} \\
    &0 \leq \rho_{2}^{k[q]} \leq \frac{a_u^k}{Q} \quad \forall k \in V_{0}, \forall q \in\{1,2, \ldots, Q\}
    \label{eq:DSP(SM)_rho2-sign}
\end{align}
\end{subequations}
where dual variables $\rho_{1}^{k[q]}$ and $\rho_{2}^{k[q]}$ are associated with constraints \eqref{eq:Hl_v1}, and \eqref{eq:Hu_v2}, respectively, $\forall k \in V_{0}, \forall q \in\{1,2, \ldots, Q\}$. 

$\tilde{\text{DSP}}$ can be decomposed into two main subproblems for $\boldsymbol{\rho}_1$ and $\boldsymbol{\rho}_2$, each of which also decomposes into $|V_0|$ subproblems, one for each $k \in V_0$. For a given $k\in V_0$, let us assume the costs $c^{k[q]}=\sum_{(i, j) \in A} t_{i j}^{[q]} {y}_{i j}^{k}$ of variables $\boldsymbol{\rho}_1$ and $\boldsymbol{\rho}_2$  have been sorted to get $ c^{k[\Lambda_1]} \leq  c^{k[\Lambda_2]} \leq \ldots \leq   c^{k[\Lambda_Q]}$. For ease of presentation, we let $\Lambda=(1,2,\ldots Q)$. We then define a critical index $P_1^{k}\in \{1, 2, \ldots, Q\}$ such that
\begin{align}
\label{eq:Qorder1_proof}
    \sum_{q=1}^{P_1^{k}-1} \frac{a_\ell^k}{Q} < a_w^k \leq  \sum_{q=1}^{P_1^{k}} \frac{a_\ell^k}{Q}. 
\end{align}
Then an optimal value for variables $\boldsymbol{\rho}_1$ is obtained by setting
\begin{equation}
\label{eq:dualsol1}
   {\bar \rho}_1^{k[q]} = \begin{cases}
 \frac{a_\ell^k}{Q}  & q=1,\ldots, P_1^{k}-1;\\
  a_w^k - \sum_{q=1}^{P_1^{k}-1} \frac{a_\ell^k}{Q}  & q=P_1^{k};\\
0& q=P_1^{k}+1, \ldots, Q.
 \end{cases}
\end{equation}
This is due to the fact that we want to minimize the second term of the objective function in \eqref{eq:DSP(SM)_of}, which can be accomplished by assigning 0 to variables $\boldsymbol{\rho_1}$ with the highest costs while assigning as much as possible (at most $\frac{a_\ell^k}{Q}$ according to \eqref{eq:DSP(SM)-rho1-sign}) to variables $\boldsymbol{\rho_1}$ with the lowest cost (the first $P_1^k$ ones) to satisfy \eqref{eq:DSP(SM)_rho1}. In a similar fashion, we can define a critical index $P_2^{k}\in \{1,2, \ldots, Q\}$ such that
\begin{align}
\label{eq:Qorder2_proof}
    \sum_{q=P_2^{k}+1}^{Q} \frac{a_u^k}{Q} < a_w^k \leq  \sum_{q=P_2^{k}}^{Q} \frac{a_u^k}{Q}, 
\end{align}
and obtain an optimal value for variables $\boldsymbol{\rho}_2$ stated as
\begin{equation}
\label{eq:dualsol2}
   {\bar \rho}_2^{k[q]} = \begin{cases}
 0& q=1, \ldots, P_2^{k}-1;\\
  a_w^k - \sum_{q=P_2^{k}+1}^{Q} \frac{a_u^k}{Q}  & q=P_2^{k};\\
  \frac{a_u^k}{Q}  & q=P_2^{k}+1, \ldots, Q.
 \end{cases}
\end{equation}

We then accordingly construct the primal solutions for each customer $k\in V_0$ as follows:
\begin{align*}
  &{\bar \ell^k}=c^{k[P_1^{k}]}=\sum_{(i, j) \in A} t_{i j}^{[P_1^{k}]} {y}_{i j}^{k},
    &\bar v_1^{k[q]} = \begin{cases}
 c^{k[P_1^{k}]}-c^{k[q]}  & q=1,\ldots, P_1^{k}-1;\\
0& q=P_1^{k}, \ldots, Q.
 \end{cases}\\
&{\bar u^k}=c^{k[P_2^{k}]}=\sum_{(i, j) \in A} t_{i j}^{[P_2^{k}]} {y}_{i j}^{k},
    &\bar v_2^{k[q]} = \begin{cases}
0& q=1, \ldots, P_2^{k};\\
c^{k[q]}-c^{k[P_2^{k}]}  & q=P_2^{k}+1,\ldots, Q.
 \end{cases}
\end{align*}

Optimality comes from the feasibility of the primal and dual solutions for their problems and from the fact that the primal cost is equal to the dual cost. This can be evidently achieved by replacing the solutions in \eqref{eq:SP(SM)} and \eqref{eq:DSP(SM)}. It is, however, necessary to show that $P_1^{k} \leq P_2^{k}$ in order to satisfy the constraint $\ell^k \leq u^k$ for each customer $k \in V_0$. From \eqref{eq:Qorder1_proof} and \eqref{eq:Qorder2_proof}, we can gain $\frac{P_1^{k}-1}{Q} < \frac{a_w^k}{a_\ell^k}$ and $\frac{Q-P_2^{k}}{Q} < \frac{a_w^k}{a_u^k}$, respectively. Then considering the fact that $a_w^k/a_\ell^k + a_w^k/a_u^k \leq 1$, we have  $P_1^{k} \leq P_2^{k}$.
    \hfill $\square$
    \vspace{0.5cm}
    \phantomsection
\pdfbookmark[2]{Appendix C: Computing Subgradients in Proposition 5}{appC}
\section*{Appendix C: Computing Subgradients in Proposition \ref{prop:subgradient}}
\label{sec:subgradient}

Given $(\hat{\boldsymbol{x}},\hat{\boldsymbol{y}})$, one can efficiently solve the subproblem SP(SM) for fixed $(\boldsymbol{x},\boldsymbol{y})=(\hat{\boldsymbol{x}},\hat{\boldsymbol{y}})$. We can formulate the Lagrangian function for problem SP(SM) as (see $\tilde{\text{SP}}^k$ in \eqref{eq:SP(SM)} and Appendix \ref{sec:AppC})
\begin{equation}
\begin{split}
\nonumber
    & L\left(\boldsymbol{\ell}, \boldsymbol{u}, \boldsymbol{v_{1}}^{[q]}, \boldsymbol{v_{2}}^{[q]}, \boldsymbol{\rho_{1}}^{[q]}, \boldsymbol{\rho_{2}}^{[q]}, \boldsymbol{\lambda_{1}}, \boldsymbol{\lambda_{2}}, \boldsymbol{\delta_{1}}^{[q]}, \boldsymbol{\delta_{2}}^{[q]}\right)
    = \sum_{k \in V_{0}} L^{k}\left(\ell^{k}, u^{k}, v_{1}^{k[q]}, v_{2}^{k[q]}, \rho_{1}^{k[q]}, \rho_{2}^{k[q]}, \lambda_{1}^{k}, \lambda_{2}^{k}, \delta_{1}^{k[q]}, \delta_{1}^{k[q]}\right) = \\
    & \sum_{k \in V_{0}}\left( a_w^k(u^k - \ell^k) + \frac{a_\ell^k}{Q} \sum_{q=1}^{Q} v_{1}^{k[q]} + \frac{a_u^k}{Q} \sum_{q=1}^{Q} v_{2}^{k[q]}\right) + \sum_{k \in V_{0}} \sum_{q=1}^{Q} \rho_{1}^{k[q]}\left(\ell^{k}-\sum_{(i, j) \in A}\left(t_{i j}^{[q]} \hat{y}_{i j}^{k}\right)-v_{1}^{k[q]}\right) + \\ 
    & \sum_{k \in V_{0}} \sum_{q=1}^{Q} \rho_{2}^{k[q]}\left(\sum_{(i, j) \in A}\left(t_{i j}^{[q]} \hat{y}_{i j}^{k}\right)-u^{k}-v_{2}^{k[q]}\right) -\sum_{k \in V_{0}} \lambda_{1}^{k} \ell^{k}
    -\sum_{k \in V_{0}} \lambda_{2}^{k} u^{k}
    -\sum_{k \in V_{0}} \sum_{q=1}^{Q} \delta_{1}^{k[q]} v_{1}^{k[q]}-\sum_{k \in V_{0}} \sum_{q=1}^{Q} \delta_{2}^{k[q]} v_{2}^{k[q]},
\end{split}
\end{equation}
where $\lambda_{1}^{k}, \lambda_{2}^{k}, \delta_{1}^{k[q]}, \delta_{2}^{k[q]}$ are associated with the range constraints of primal variables $\ell^{k}, u^{k}, v_{1}^{k[q]}, v_{2}^{k[q]}$, respectively, in \eqref{eq:SM_sign}, $\forall k \in V_{0}, \forall q \in\{1,2, \ldots, Q\}$.

Let $\ell^{k^{*}}, u^{k^{*}}, v_{1}^{k[q]^{*}}, v_{2}^{k[q]^{*}}$ be the optimal primal solutions found, and let $\rho_{1}^{k[q]^{*}}, \rho_{2}^{k[q]^{*}}, \lambda_{1}^{k^{*}}, \lambda_{2}^{k^{*}}, \delta_{1}^{k[q]^{*}}, \delta_{1}^{k[q]^{*}}$ be the optimal dual variables. Using Lagrangian duality and Karush–Kuhn–Tucker (KKT) conditions, and assuming constraint qualifications hold, a subgradient $\forall k \in V_{0}, \forall(i, j) \in A$ can be obtained as (see Geoffrion 1972):


$$
\hat{s}_{i j}^{k} = \frac{\partial L^{k}\left(\ell^{k^{*}}, u^{k^{*}}, v_{1}^{k[q]^{*}}, v_{2}^{k[q]^{*}}, \rho_{1}^{k[q]^{*}}, \rho_{2}^{k[q]^{*}}, \lambda_{1}^{k^{*}}, \lambda_{2}^{k^{*}}, \delta_{1}^{k[q]^{*}}, \delta_{2}^{k[q]^{*}}\right)} {\partial y_{i j}^{k}} = \sum_{q=1}^{Q} t_{i j}^{[q]}\left(\rho_{2}^{k[q]^{*}}-\rho_{1}^{k[q]^{*}}\right), 
$$
where $\rho_{1}^{k[q]^{*}}$ and $\rho_{2}^{k[q]^{*}}$ are computed using \eqref{eq:dualsol1} and \eqref{eq:dualsol2} of the proof of the Proposition \ref{Proposition3}, respectively.
    \vspace{0.5cm}
\phantomsection
\pdfbookmark[2]{Appendix D: Proof of Proposition 6}{appD}
\section*{Appendix D: Proof of Proposition \ref{theorem:Robust2-2}} 
\label{sec:AppF}

For the convenience of analysis, let us define for each customer $k \in V_0$
\begin{equation*}
    H(\ell^k, u^k) = 
    a_w^k(u^k-\ell^k) + a_{\ell}^k \mathbb{E}_{\mathrm{P}}[\ell^k-\tau(\boldsymbol{y}^k,\boldsymbol{\tilde{t}})]^+ + 
    a_u^k \mathbb{E}_{\mathrm{P}}[\tau(\boldsymbol{y}^k,\boldsymbol{\tilde{t}})-u^k]^+.
\end{equation*}
Also, according to \eqref{eq:U_mu} and \eqref{eq:U_C}, we have the uncertainty set $\mathcal{U}_{(\boldsymbol{\hat{\mu}}, \boldsymbol{\hat{\mathcal{C}}})} $ as follows
\begin{equation*}
    \mathcal{U}_{(\boldsymbol{\hat{\mu}}, \boldsymbol{\hat{\mathcal{C}}})} = \left\{(\boldsymbol{\mu}, \boldsymbol{\mathcal{C}}) : \left( \boldsymbol{\mu}-\boldsymbol{\hat{\mu}}\right)^{\top} \boldsymbol{\hat{\mathcal{C}}}^{-1}\left(\boldsymbol{\mu}-\boldsymbol{\hat{\mu}}\right) \leq \alpha_1, 
    \|\boldsymbol{\mathcal{C}}-\boldsymbol{\hat{\mathcal{C}}} \|_F \leq \alpha_2 \right\}.
\end{equation*}

We first gain the supremum of $H(\ell^k, u^k)$ with respect to $\mathrm{P} \in \mathcal{F}_{(\boldsymbol{\bar{\mu}}, \boldsymbol{\bar{\mathcal{C}}})}$, and then with respect to $\left(\boldsymbol{\bar{\mu}}, \boldsymbol{\bar{\mathcal{C}}}  \right)$ within the uncertainty set $\mathcal{U}_{(\boldsymbol{\hat{\mu}}, \boldsymbol{\hat{\mathcal{C}}})}$. By plugging the optimal time window for each customer $k \in V_0$ found by \eqref{eq:lbar_robust} and \eqref{eq:ubar_robust} in Proposition \ref{eq:optimal-time-window} into $H(\ell^k, u^k)$ and from the definition of $\mathbb{D}$ in \eqref{eq:ambiguity_set1}, we have

\begin{equation*}
\begin{aligned}
    \sup _{\mathrm{P} \in \mathbb{D}} H(\ell^k, u^k) = \sup_{(\boldsymbol{\bar{\mu}}, \boldsymbol{\bar{\mathcal{C}}}) \in \mathcal{U}_{(\boldsymbol{\hat{\mu}}, \boldsymbol{\hat{\mathcal{C}}})}} \sup_{\mathrm{P} \in \mathcal{F}_{(\boldsymbol{\bar{\mu}}, \boldsymbol{\bar{\mathcal{C}}})}} H(\ell^k, u^k) 
    = \sup_{(\boldsymbol{\bar{\mu}}, \boldsymbol{\bar{\mathcal{C}}}) \in \mathcal{U}_{(\boldsymbol{\hat{\mu}}, \boldsymbol{\hat{\mathcal{C}}})}} \{ (\Gamma_\ell^k + \Gamma_u^k) \sqrt{\boldsymbol{y}^{k^{\top}} \boldsymbol{\Bar{\mathcal{C}}} \boldsymbol{y}^{k}} \} \\
    = (\Gamma_\ell^k + \Gamma_u^k) \sup_{\boldsymbol{\bar{\mathcal{C}}} \in \mathcal{U}_{\boldsymbol{\hat{\mathcal{C}}}}} \sqrt{\boldsymbol{y}^{k^{\top}} \boldsymbol{\Bar{\mathcal{C}}} \boldsymbol{y}^{k}}.
\end{aligned}
\end{equation*}

Now, let $\boldsymbol{\tilde{\mathcal{C}}} \triangleq \boldsymbol{\bar{\mathcal{C}}} - \boldsymbol{\hat{\mathcal{C}}}$. Then the problem $\sup_{\boldsymbol{\bar{\mathcal{C}}} \in \mathcal{U}_{\boldsymbol{\hat{\mathcal{C}}}}} \sqrt{\boldsymbol{y}^{k^{\top}} \boldsymbol{\bar{\mathcal{C}}} \boldsymbol{y}^{k}}$ can be formulated as 
\begin{align*}
    \sup_{\boldsymbol{\tilde{\mathcal{C}}}} \sqrt{\boldsymbol{y}^{k^{\top}} \boldsymbol{\tilde{\mathcal{C}}} \boldsymbol{y}^{k} + \boldsymbol{y}^{k^{\top}} \boldsymbol{\hat{\mathcal{C}}} \boldsymbol{y}^{k}} \\
    \mbox{s.t.} \quad
    \|\boldsymbol{\mathcal{\tilde{C}}} \|_F \leq \alpha_2.
\end{align*}

For this problem, we rewrite the supremum as follows
\[
\sup_{\boldsymbol{\tilde{\mathcal{C}}} : \|\boldsymbol{\tilde{\mathcal{C}}}\|_F \leq \alpha_2} \sqrt{\boldsymbol{y}^{k \top} \boldsymbol{\tilde{\mathcal{C}}} \boldsymbol{y}^k + \boldsymbol{y}^{k \top} \boldsymbol{\hat{\mathcal{C}}} \boldsymbol{y}^k} 
= \sup_{\boldsymbol{\tilde{\mathcal{C}}} : \|\boldsymbol{\tilde{\mathcal{C}}}\|_F \leq \alpha_2} \sqrt{\boldsymbol{\tilde{\mathcal{C}}} \circ \boldsymbol{y}^k \boldsymbol{y}^{k \top} + \boldsymbol{y}^{k \top} \boldsymbol{\hat{\mathcal{C}}} \boldsymbol{y}^k},
\]
where \(\circ\) denotes the Frobenius inner product, which satisfies the following inequality according to the Cauchy-Schwarz inequality
\[
\boldsymbol{\tilde{\mathcal{C}}} \circ (\boldsymbol{y}^k \boldsymbol{y}^{k \top}) \leq \|\boldsymbol{\tilde{\mathcal{C}}}\|_F \cdot \|\boldsymbol{y}^k \boldsymbol{y}^{k \top}\|_F,
\]
where the equality holds if and only if \(\boldsymbol{\tilde{\mathcal{C}}}\) is proportional to \(\boldsymbol{y}^k \boldsymbol{y}^{k \top}\), i.e., there exists a scalar \(\lambda\) such that
\[
\boldsymbol{\tilde{\mathcal{C}}} = \lambda (\boldsymbol{y}^k \boldsymbol{y}^{k \top}).
\]
This alignment (that matrix $\boldsymbol{\tilde{\mathcal{C}}}$ is perfectly aligned with $\boldsymbol{y}^k \boldsymbol{y}^{k \top}$) ensures that the Frobenius inner product achieves its maximum possible value.
Substituting $\boldsymbol{\tilde{\mathcal{C}}} = \lambda (\boldsymbol{y}^k \boldsymbol{y}^{k \top})$ into the Frobenius norm gives
\[
\|\boldsymbol{\tilde{\mathcal{C}}}\|_F = |\lambda| \cdot \|\boldsymbol{y}^k \boldsymbol{y}^{k \top}\|_F.
\]
To satisfy the constraint \(\|\boldsymbol{\tilde{\mathcal{C}}}\|_F = \alpha_2\), we solve for \(\lambda\)
\[
|\lambda| = \frac{\alpha_2}{\|\boldsymbol{y}^k \boldsymbol{y}^{k \top}\|_F},
\]
using which the aligned \(\boldsymbol{\tilde{\mathcal{C}}}\) becomes
\[
\boldsymbol{\tilde{\mathcal{C}}} = \alpha_2 \frac{\boldsymbol{y}^k \boldsymbol{y}^{k \top}}{\|\boldsymbol{y}^k \boldsymbol{y}^{k \top}\|_F},
\]
whose substitution into the Frobenius inner product gives
\[
\boldsymbol{\tilde{\mathcal{C}}} \circ (\boldsymbol{y}^k \boldsymbol{y}^{k \top}) = \left(\alpha_2 \frac{\boldsymbol{y}^k \boldsymbol{y}^{k \top}}{\|\boldsymbol{y}^k \boldsymbol{y}^{k \top}\|_F}\right) \circ (\boldsymbol{y}^k \boldsymbol{y}^{k \top}) = \alpha_2 \frac{\|\boldsymbol{y}^k \boldsymbol{y}^{k \top}\|_F^2}{\|\boldsymbol{y}^k \boldsymbol{y}^{k \top}\|_F} = \alpha_2 \|\boldsymbol{y}^k \boldsymbol{y}^{k \top}\|_F.
\]
This achieves the maximum possible value of the inner product under the Frobenius norm constraint. Thus, the supremum becomes
\[
\sup_{\boldsymbol{\tilde{\mathcal{C}}} : \|\boldsymbol{\tilde{\mathcal{C}}}\|_F \leq \alpha_2} \sqrt{\boldsymbol{y}^{k \top} \boldsymbol{\tilde{\mathcal{C}}} \boldsymbol{y}^k + \boldsymbol{y}^{k \top} \boldsymbol{\hat{\mathcal{C}}} \boldsymbol{y}^k} = \sqrt{\alpha_2 \|\boldsymbol{y}^k\|_2^2 + \boldsymbol{y}^{k \top} \boldsymbol{\hat{\mathcal{C}}} y^k},
\]
which results in
\begin{equation*}
    \sup_{\boldsymbol{\bar{\mathcal{C}}} \in \mathcal{U}_{\boldsymbol{\hat{\mathcal{C}}}}} \sqrt{{\boldsymbol{y}^{k}}^{\top} \boldsymbol{\bar{\mathcal{C}}} \boldsymbol{y}^{k}} = 
    \sqrt{{\boldsymbol{y}^{k}}^{\top} (\boldsymbol{\hat{\mathcal{C}}} + \alpha_2 I_{|A|}) \boldsymbol{y}^{k}},
\end{equation*}
where $I_{|A|}$ is the identity matrix of size ${|A|}$. Hence, we finaly have
\begin{equation*}
\label{eq:sup_objective}
     \sup _{\mathrm{P} \in \mathbb{D}} H(\ell^k, u^k) = (\Gamma_\ell^k + \Gamma_u^k) \sqrt{{\boldsymbol{y}^{k}}^{\top} (\boldsymbol{\hat{\mathcal{C}}} + \alpha_2 I_{|A|}) \boldsymbol{y}^{k}}.
\end{equation*}

In a similar fashion, to handle constraint \eqref{eq:om_budget}, we have
\begin{align*}
    \sup _{\mathrm{P} \in \mathbb{D}} \mathbb{E}_{\mathrm{P}}\left(\tilde{\boldsymbol{t}}\right)^{\top}\boldsymbol{x} = 
    \sup_{(\boldsymbol{\bar{\mu}}, \boldsymbol{\bar{\mathcal{C}}}) \in \mathcal{U}_{(\boldsymbol{\hat{\mu}}, \boldsymbol{\hat{\mathcal{C}}})}} \sup_{\mathrm{P} \in \mathcal{F}_{(\boldsymbol{\bar{\mu}}, \boldsymbol{\bar{\mathcal{C}}})}} \mathbb{E}_{\mathrm{P}}\left(\tilde{\boldsymbol{t}}\right)^{\top}\boldsymbol{x} =
    \sup_{(\boldsymbol{\bar{\mu}}, \boldsymbol{\bar{\mathcal{C}}}) \in \mathcal{U}_{(\boldsymbol{\hat{\mu}}, \boldsymbol{\hat{\mathcal{C}}})}} \left\{ \boldsymbol{\bar{\mu}}^{\top} \boldsymbol{x} \right\} =
    \sup_{\boldsymbol{\bar{\mu}} \in \mathcal{U}_{\boldsymbol{\hat{\mu}}}} \boldsymbol{\bar{\mu}}^{\top} \boldsymbol{x}.
\end{align*}
This way, the optimal solution to $\sup_{\boldsymbol{\bar{\mu}} \in \mathcal{U}_{\boldsymbol{\hat{\mu}}}} \boldsymbol{\bar{\mu}}^{\top} \boldsymbol{x}$ clearly can be shown to be
\begin{equation*}
\label{eq:sup_constraint}
    \sup_{\boldsymbol{\bar{\mu}} \in \mathcal{U}_{\boldsymbol{\hat{\mu}}}} \boldsymbol{\bar{\mu}}^{\top} \boldsymbol{x} =
    \boldsymbol{\hat{\mu}}^{\top} \boldsymbol{x} + \sqrt{\alpha_1} \sqrt{\boldsymbol{x}^{\top} \boldsymbol{\hat{\mathcal{C}}} \boldsymbol{x}}.
\end{equation*}
\hfill $\square$

    \vspace{0.5cm}
    \phantomsection
\pdfbookmark[2]{Appendix E: Generating a Positive Semidefinite Covariance Matrix}{appE}
\section*{Appendix E: Generating a Positive Semidefinite Covariance Matrix}
\label{appendix:CovMatrix}

In this section, we explain how to generate a positive semidefinite covariance matrix $\hat {\boldsymbol{\mathcal{C}}}$ used for evaluating our proposed frameworks to design service time windows.
Let $\varrho_{ijrs}$ represent the correlation coefficient between the travel times on arcs $(i,j) \in A$ and $(r,s) \in A$, and $\sigma_{ij}$ be the standard deviation of travel time for traversing arc $(i,j)$. This way, the entry of the covariance matrix  $\hat {\boldsymbol{\mathcal{C}}} \in \mathbb{R}^{|A| \times |A|}$ for this pair of arcs is given by
\begin{equation}
    \label{eq:covmatrix}
    \hat {\mathcal{C}}_{ijrs}= \begin{cases}
    \sigma_{ij}^2 & \text{if } (i,j)=(r,s) \\
    \varrho_{ijrs}\sigma_{ij}\sigma_{rs} & \text{if } (i,j)\neq (r,s).
    \end{cases}
\end{equation}

According to \eqref{eq:covmatrix}, $\boldsymbol{\hat{\mathcal{C}}}=\boldsymbol{R} \circ (\boldsymbol{D}\boldsymbol{D}^\top)$, where $\boldsymbol{R}\in\left[-1,1\right]^{|A|\times|A|}$ is the correlation matrix, $\boldsymbol{D}\in\mathbb{R}^{|A|}$ is the standard deviation vector, and ``$\circ$'' is the Hadamard (element-wise) product operation. Note that according to the Schur product theorem, the Hadamard product of two positive semidefinite matrices is also a positive semidefinite matrix. $(\boldsymbol{D}\boldsymbol{D}^\top)$ is a rank-one matrix, and thus positive semidefinite. Therefore, if we can generate a positive semidefinite matrix $\boldsymbol{R}$, the resulting covariance matrix $\hat{\boldsymbol{\mathcal{C}}}$ will be positive semidefinite as well. To generate the appropriate matrix $\boldsymbol{R}$ and vector $\boldsymbol{D}$, we modified the method used in \cite{rostami2021branch} to reflect the fact that as the distance between arcs increases, their travel times' correlation decreases. The procedure involves the following steps: 
\begin{itemize}
    \item Since the product of any matrix and its transpose is semidefinite, we set $\boldsymbol{R}=\boldsymbol{E}\boldsymbol{E}^{\top}$, where $\boldsymbol{E} \in \mathbb{R}^{|A|\times|V|}$, and hence $\boldsymbol{R}$ will be semidefinite, each element of which represents $\varrho_{ijrs}$, the correlation coefficient between arcs $(i,j) \in A$ and $(r,s) \in A$. Matrix $\boldsymbol{E}$ is gained from a $|A|\times|V|$ matrix $\boldsymbol{\Tilde{E}}=[\frac{1}{1+d_{ijk}}]$, where $d_{ijk}$ is the minimum distance (number of edges) between arc $(i,j)\in A$ and node $k\in V$. This way, smaller travel time correlations will be assigned to the arcs that are more distant from each other. $\boldsymbol{E}$ is actually the matrix resulting from normalizing the rows of $\boldsymbol{\Tilde{E}}$ to have length one. Therefore, matrix $\boldsymbol{R}$'s entries lie in [0, 1]. To generate negative correlations between the arcs, we multiply each element of matrix $\boldsymbol{\Tilde{E}}$ by $-1$ with the probability of $5\%$ resulting in a matrix $\boldsymbol{R}$ with all entries in $[-1, 1]$.
    
    \item Given the expected travel time $\mu_{ij}$, we generate vector $\boldsymbol{D}=[\sigma_{ij}]$ with $\sigma_{ij}=CV_{ij} \times \mu_{ij}$, where $CV_{ij}$ is the coefficient of variation for arc $(i,j)\in A$ and is drawn from a uniform distribution in the range [0.01, 0.2].
\end{itemize}  
    \vspace{0.5cm}
    \phantomsection
\pdfbookmark[2]{Appendix F: Decomposition Algorithms' Performance}{appF}
\section*{Appendix F: Decomposition Algorithms' Performance}
\label{sec:AppG}

The results are presented in Tables \ref{tab:Benders_SM1} to \ref{tab:OA_RM2}. For each SM and RM, we report the results in two tables. Tables \ref{tab:Benders_SM1} and \ref{tab:OA_RM1} present the results for the instances that were solved to optimality by all the algorithms within the time limit. The objective is to compare the performance of the algorithms in terms of computational time. Tables \ref{tab:Benders_SM2} and \ref{tab:OA_RM2} have been divided into two parts. In each of them, the upper part reports the results for the instances where at least one of our decomposition-based algorithms solved the instance within the time limit, while the lower part presents the results for instances where all the algorithms reached the time limit. The objective of these tables is to compare the algorithms in terms of computational time when applicable and optimality gap when algorithms reach the time limit.

In all the tables, we use $\#Customers$ to represent the number of customers in each instance. For each algorithm, we use $\# BBnodes$ to show the number of branch-and-bound nodes explored in the decision tree, $Time$ to show the computational time (in seconds) to solve each instance, and $Gap$ to display the optimality gap between the upper bound (UB; best integer solution found) and the lower bound (LB). For each decomposition-based algorithm, $\# Cuts$ indicates the number of Benders/OA user cuts added to the master problem within the tree. Moreover, in each table, the last two columns display the percentage improvements achieved by the decomposition-based algorithms compared to the CPLEX base algorithm in terms of either the computational time or the optimality gap. We used $((x_0 - x_d)/x_0) \times 100 $ formula to compute such an improvement quantity, where $x_d$ is the time/gap by the decomposition-based algorithms and $x_0$ stands for those of the base algorithm. In Tables \ref{tab:Benders_SM2} and \ref{tab:OA_RM2}, if the decomposition algorithm was able to reach the optimal solution within the time limit (a gap of zero), ++ presents the gap improvement instead of 100. We show the time by TL when an algorithm reached the time limit and could not solve the instance to optimality. For each instance, the results of the best algorithm in terms of time/gap are presented in bold.

As it can be seen from Table \ref{tab:Benders_SM1}, CPLEX was able to solve instances with up to 20 (except 17) customers in reasonable times and also solved an instance with 27 customers within the time limit. For these instances, both the CPLEX+BSCut and CPLEX+BMCuts outperform the CPLEX in terms of computational time except for two small instances with 10 and 12 customers where the base algorithm performs better. Comparing CPLEX+BSCut and CPLEX+BMCuts, we can observe that overall the latter outperforms the former, which indicates how adding the cuts shrinks the feasible set of LP relaxation more efficiently and hence allows it to explore more branch-and-bound nodes more effectively. From Table \ref{tab:Benders_SM2}, we can see that  CPLEX+BMCuts outperforms the others in terms of the computational time whenever instances were solved to optimality within the time limit, and in terms of the optimality gap when all the algorithms reached the time limit.  

The results in Tables \ref{tab:OA_RM1} and \ref{tab:OA_RM2} follow the similar patterns as those for solving the SM in  Tables \ref{tab:Benders_SM1} and \ref{tab:Benders_SM2}. More precisely, RM's instances with up to 20 customers and with 23 customers can be solved by CPLEX in reasonable times but get more difficult as the size increases. Where all the algorithms solve the instances to optimality, the CPLEX+OASCut and CPLEX+OAMCuts are, on average, 72.19\% and 76.57\% faster than CPLEX, respectively. As observed from Table \ref{tab:OA_RM2}, more instances remain unsolved within the time limit compared to the SM model, which indicates the difficulty of solving the RM model when the number of customers increases.   

\begin{table}[h]
\caption{Evaluating CPLEX  with and without Benders cuts to solve the SM on complete graphs when all algorithms reach the optimal solution within the time limit}
\label{tab:Benders_SM1}
\resizebox{\textwidth}{!}{%
\begin{tabular}{lllllllllllllll}
\toprule
\multirow{2}{*}{$\#Customers$} &  & \multicolumn{2}{c}{CPLEX} &  & \multicolumn{3}{c}{CPLEX+BSCut}                 &  & \multicolumn{3}{c}{CPLEX+BMCuts}                    &  & \multicolumn{2}{c}{\%Improvement (time)} \\ \cmidrule{3-4} \cmidrule{6-8} \cmidrule{10-12} \cmidrule{14-15}
                           &  & $\#BBnodes$   & $Time (s)$           &  & $\#BBnodes$  & $\#Cuts$    & $Time (s)$            &  & $\#BBnodes$ &  $\#Cuts$ & $Time (s)$                       &  & $Single$ & $Multiple$                              \\\midrule 
10                         &  & 127         & \textbf{47.03}     &  & 2          & 4         & 71.09               &  & 5          & 36        & 69.91                           &  & N/A & N/A                          \\
11                         &  & 2,160       & 218.44             &  & 1,825      & 12        & \textbf{88.19}      &  & 825        & 82        & 156.45                          &  & 59.63 & 28.38                       \\
12                         &  & 74          & \textbf{49.41}     &  & 0          & 3         & 136.89              &  & 2          & 28        & 117.26                          &  & N/A & N/A                          \\
13                         &  & 859         & 313.17             &  & 1,061      & 9         & \textbf{182.86}     &  & 922        & 149       & 252.03                          &  & 41.61 & 19.52                      \\
14                         &  & 647         & 402.97             &  & 280        & 7         & 194.03              &  & 278        & 59        & \textbf{132.52}                 &  & 51.85 & 67.11                      \\
15                         &  & 29          & 298.14             &  & 441        & 4         & \textbf{112.36}     &  & 260        & 25        & 132.02                          &  & 62.31 & 55.72                     \\
16                         &  & 37,273      & 12,378.94          &  & 32,304     & 16        & 464.53              &  & 7,460      & 105       & \textbf{252.41}                 &  & 96.25 & 97.96                      \\

18                         &  & 24,524      & 15,116.70          &  & 4,641      & 18        & 347.81              &  & 2,840      & 134       & \textbf{321.47}                 &  & 97.70 & 97.87                      \\
19                         &  & 2,566       & 2,557.95           &  & 2,698      & 25        & 516.94              &  & 1,591      & 213       & \textbf{310.47}                 &  & 79.79 & 87.86                      \\
20                         &  & 2,218       & 2,897.50           &  & 5,783      & 66        & 893.55              &  & 2,436      & 481       & \textbf{465.42}                 &  & 69.16 & 83.94                      \\

27                         &  & 6,148       & 16,732.23          &  & 33,001     & 63        & 5,763.08            &  & 8,077      & 592       & \textbf{1,912.48}               &  & 65.56 & 88.57 \\
 \bottomrule
\end{tabular}%
}
\end{table}

\begin{table}[htp]
\caption{Evaluating CPLEX  with and without Benders cuts to solve the SM on complete graphs when at least one algorithm cannot reach the optimal solution within the time limit}
\label{tab:Benders_SM2}
\resizebox{\textwidth}{!}{%
\begin{tabular}{llllllllllllllllll}
\toprule
\multirow{2}{*}{$\#Customers$} &  & \multicolumn{3}{c}{CPLEX} &  & \multicolumn{4}{c}{CPLEX+BSCut}           &  & \multicolumn{4}{c}{CPLEX+BMCuts}                    &  & \multicolumn{2}{c}{\%Improvement (gap)} \\ \cmidrule{3-5} \cmidrule{7-10} \cmidrule{12-15}\cmidrule{17-18}
                           &  & $\#BBnodes$  & $Time (s)$       &$Gap (\%)$    &  & $\#BBnodes$ & $\#Cuts$ & $Time (s)$        & $Gap (\%)$    &  & $\#BBnodes$ & $\#Cuts$ & $Time (s)$           & $Gap (\%)$             &  & $Single$ & $Multiple$                              \\  \midrule 
17                         &  & 33,885      & TL      & 26.68 &  & 43,815     & 15        & 455.67          & 0.00     &  & 30,042     & 136       & \textbf{402.53}    & 0.00              &  & ++ & ++                      \\
21                         &  & 9,935       & TL      & 50.12 &  & 66,110     & 23        & 1,802.88        & 0.00     &  & 58,779     & 129       & \textbf{1,341.78}  & 0.00              &  & ++ & ++                         \\
22                         &  & 5,955       & TL      & 42.84 &  & 84,629     & 27        & 2,658.14        & 0.00     &  & 55,562     & 276       & \textbf{1,434.02}  & 0.00              &  & ++ & ++                           \\
23                         &  & 6,991       & TL      & 20.27 &  & 34,253     & 75        & 3,380.72        & 0.00     &  & 26,800     & 441       & \textbf{1,772.28}  & 0.00              &  & ++ & ++                           \\

25                         &  & 6,397       & TL      & 15.13 &  & 34,195     & 107       & 6,224.64        & 0.00     &  & 28,995     & 768       & \textbf{3,341.00}  & 0.00              &  & ++ & ++                           \\
26                         &  & 4,804       & TL      & 53.18 &  & 156,920    & 29        & 15,793.13       & 0.00     &  & 119,985    & 206       & \textbf{9,874.73}  & 0.00              &  & ++ & ++                          \\
28                         &  & 3,975       & TL      & 44.38 &  & 22,244     & 46        & 5,287.77        & 0.00     &  & 21,574     & 168       & \textbf{2,847.98}  & 0.00              &  & ++ & ++                           \\
29                         &  & 1,501       & TL      & 77.25 &  & 84,389     & 57        & TL              & 23.89    &  & 80,205     & 483       & \textbf{10,922.41} & 0.00              &  & 69.07 & ++                        \\
\midrule
24                         &  & 5,804       & TL      & 57.46 &  & 184,114    & 34        & TL              & 42.56    &  & 343,806    & 303       & TL & \textbf{23.39}                    &  & 25.93 & 59.29                      \\
30                         &  & 1,172       & TL      & 78.77 &  & 53,077     & 126       & TL              & 47.28    &  & 65,156     & 1,308  & TL & \textbf{41.99}                    &  & 39.98 &46.69                      \\ \bottomrule
\end{tabular}%
}
\end{table}

\begin{table}[t]
\caption{Evaluating CPLEX  with and without outer approximation cuts to solve the RM on complete graphs when all algorithms reach the optimal solution within the time limit}
\label{tab:OA_RM1}
\resizebox{\textwidth}{!}{%
\begin{tabular}{lllllllllllllll}
\toprule
\multirow{2}{*}{$\#Customers$} &  & \multicolumn{2}{c}{CPLEX} &  & \multicolumn{3}{c}{CPLEX+OASCut}           &  & \multicolumn{3}{c}{CPLEX+OAMCuts}                    &  & \multicolumn{2}{c}{\%Improvement (time)} \\ \cmidrule{3-4} \cmidrule{6-8} \cmidrule{10-12} \cmidrule{14-15}
                           &  & $\#BBnodes$  & $Time (s)$       &  & $\#BBnodes$ & $\#Cuts$ & $Time (s)$                     &  & $\#BBnodes$ &  $\#Cuts$ & $Time (s)$                    &  & $Single$ & $Multiple$                              \\\midrule 
10                         &  & 216         & \textbf{3.75}     &  & 1,017      & 10        & 4.25                           &  & 505        & 47        & 4.20                           &  & N/A & N/A                          \\
11                         &  & 4,116       & 22.20             &  & 16,794     & 13        & 16.63                          &  & 12,398     & 49        & \textbf{16.06}                 &  & 25.09 & 27.66                      \\
12                         &  & 308         & 35.95             &  & 1,281      & 15        & 12.36                          &  & 596        & 48        & \textbf{10.06}                 &  & 65.62 & 72.02                      \\
13                         &  & 1,809       & 37.84             &  & 9,954      & 6         & 30.55                          &  & 4,025      & 147       & \textbf{23.00}                 &  & 19.27 & 39.22                      \\
14                         &  & 2,707       & 117.22            &  & 10,943     & 10        & 28.28                          &  & 6,866      & 159       & \textbf{25.06}                 &  & 75.87 & 78.62                      \\
15                         &  & 1,098       & 61.95             &  & 4,546      & 2         & \textbf{15.31}                 &  & 4,219      & 30        & 16.42                          &  & 75.29 & 73.49                      \\
16                         &  & 99,871      & 7,807.39          &  & 148,705    & 18        & 325.14                         &  & 84,329     & 107       & \textbf{178.69}                &  & 95.84 & 97.71                      \\
17                         &  & 127,521     & 13,250.91         &  & 313,606    & 19        & 1,309.55                       &  & 154,106    & 171       & \textbf{484.36}                &  & 90.12 & 96.34                      \\
18                         &  & 25,092      & 4,002.02          &  & 52,710     & 24        & 269.78                         &  & 43,390     & 119       & \textbf{193.38}                &  & 93.26 & 95.17                      \\
19                         &  & 22,927      & 3,568.84          &  & 25,182     & 22        & 231.52                         &  & 23,511     & 114       & \textbf{166.94}                &  & 93.51 & 95.32                      \\
20                         &  & 7,566       & 2,564.73          &  & 31,268     & 65        & 550.03                         &  & 21,639     & 443       & \textbf{331.05}                &  & 78.55 & 87.09                      \\ 
23                         &  & 32,511      & 15,410.19         &  & 161,321    & 58        & \textbf{2,828.98}              &  & 138,820    & 425       & 3,140.36                       &  & 81.64 & 79.62                      \\
\bottomrule
\end{tabular}%
}
\end{table}
\begin{table}[p]
\caption{Evaluating CPLEX  with and without outer approximation cuts to solve the RM on complete graphs when at least one algorithm cannot reach the optimal solution within the time limit}
\label{tab:OA_RM2}
\resizebox{\textwidth}{!}{%
\begin{tabular}{llllllllllllllllll}
\toprule
\multirow{2}{*}{$\#Customers$} &  & \multicolumn{3}{c}{CPLEX} &  & \multicolumn{4}{c}{CPLEX+OASCut}           &  & \multicolumn{4}{c}{CPLEX+OAMCuts}                    &  & \multicolumn{2}{c}{\%Improvement (gap)} \\ \cmidrule{3-5} \cmidrule{7-10} \cmidrule{12-15}\cmidrule{17-18}
                           &  & $\#BBnodes$ & $Time (s)$ & $Gap (\%)$ &  & $\#BBnodes$ & $\#Cuts$ & $Time (s)$        & $Gap (\%)$    &  & $\#BBnodes$ & $\#Cuts$ & $Time (s)$ & $Gap (\%)$ &  & $Single$ & $Multiple$                              \\  \midrule  
21                         &  & 20,952      & TL      & 43.92 &  & 559,787    & 31        & 6,799.45           & 0.00   &  & 391,931    & 213       & \textbf{4,273.69}  & 0.00           &  & ++ & ++                          \\
22                         &  & 26,893      & TL      & 44.52 &  & 1,024,765  & 21        & TL                 & 8.42   &  & 723,933    & 360       & \textbf{14,457.58} & 0.00           &  & 81.09 & ++                           \\
25                         &  & 14,649      & TL      & 30.90 &  & 157,558    & 113       & \textbf{8,110.94}  & 0.00   &  & 313,882    & 880       & 13,268.84          & 0.00           &  & ++ & ++                           \\
27                         &  & 7,258       & TL      & 36.02 &  & 172,658    & 66        & \textbf{11,621.41} & 0.00   &  & 209,664    & 505       & 14,519.22          & 0.00           &  & ++ & ++                           \\
\midrule
24                         &  & 14,290      & TL      & 70.92 &  & 314,644    & 22        & TL         & 47.48          &  & 351,075    & 287       & TL                 & \textbf{45.62} &  & 33.05 & 35.67                      \\
26                         &  & 8,178       & TL      & 48.13 &  & 189,768    & 14        & TL         & \textbf{36.67} &  & 190,883    & 191       & TL                 & 41.36          &  & 23.81 & 14.07                      \\
28                         &  & 6,039       & TL      & 36.96 &  & 217,883    & 42        & TL         & \textbf{22.24} &  & 150,124    & 260       & TL                 & 27.07          &  & 39.83 & 26.76                      \\
29                         &  & 1,722       & TL      & 61.74 &  & 126,462    & 38        & TL         & 51.05          &  & 119,493    & 349       & TL                 & \textbf{43.34} &  & 17.31 & 29.80                     \\
30                         &  & 3,570       & TL      & 62.22 &  & 95,557     & 22        & TL         & 51.85          &  & 87,653     & 1,113      & TL                 & \textbf{48.23} &  & 16.67 & 22.48                      \\ \bottomrule
\end{tabular}%
}
\end{table}  
    \vspace{0.5cm}
\end{APPENDICES}

\end{document}